\documentclass[11pt, psamsfonts, reqno]{amsart}

\usepackage[T1]{fontenc}
\usepackage[latin1]{inputenc}
\usepackage[frenchb, english]{babel}
\usepackage{amsthm}
\usepackage{amscd}
\usepackage{tensor}
\usepackage{mathtools}

\usepackage{amssymb,amsmath,bbm}
\usepackage{tikz-cd}
\usepackage[all]{xy}
\usepackage{mathrsfs}
\usepackage[left=0.0cm,right=0.0cm,top=1.7cm,bottom=1.7cm]{geometry}
\usepackage{hyperref}
\usepackage[capitalise]{cleveref}
\usepackage{xcolor}
\usepackage{relsize}
\usepackage{comment}

\numberwithin{equation}{subsection}

\theoremstyle{plain}
\newtheorem{theorem}{Theorem}[subsection]
\newtheorem{Theo}{Theorem}[section]
\newtheorem{Coro}{Corollary}[section]

\newtheorem{lemma}[theorem]{Lemma}
\newtheorem{proposition}[theorem]{Proposition}
\newtheorem{prop}[theorem]{Proposition}
\newtheorem{corollary}[theorem]{Corollary}

\newtheorem*{theorem*}{Theorem}

\theoremstyle{definition}
\newtheorem{definition}[theorem]{Definition}

\newtheorem{question}[Theo]{Question}

\theoremstyle{remark}
\newtheorem{Rem}{Remark}[section]
\newtheorem{remark}[theorem]{Remark}

\newtheorem*{convention}{Convention}

\usepackage{tikz}
\usetikzlibrary{calc}

\def\Q{{\bb Q}}

\def\Z{{\bb Z}}

\def\N{{\bb N}}

\def\zp{{\Z_p}}

\def\qp{{\Q_p}}

\def\Fil{\mathrm{Fil}}
\def\Hom{\mathrm{Hom}}

\def\Sym{\mathrm{Sym}}
\def\Hom{\mathrm{Hom}}

\def\epsilon{\varepsilon}

\def\et{{\rm \acute{e}t}}

\def\GL{\mathbf{GL}}

\def\Solid{\mathrm{Solid}}

\def\sol{\mathsmaller{\square}}

\DeclareMathOperator{\Spa}{Spa}
\DeclareMathOperator{\Spec}{Spec}
\DeclareMathOperator{\Spf}{Spf}
\DeclareMathOperator{\Spd}{Spd}
\DeclareMathOperator{\gr}{gr}

\DeclareMathOperator{\Rep}{Rep}
\DeclareMathOperator{\Lie}{Lie}

\DeclareMathOperator{\Cont}{{Cont}}

\DeclareMathOperator{\Sol}{Solid}
\DeclareMathOperator{\id}{id}

\DeclareMathOperator{\dR}{dR}

\DeclareMathOperator{\AnSpec}{AnSpec}

\DeclareMathOperator{\Fl}{\n{F}\ell}

\DeclareMathOperator{\ad}{ad}

\DeclareMathOperator{\Perf}{Perf}
\DeclareMathOperator{\Perfd}{Perfd}

\DeclareMathOperator{\Aut}{Aut}

\DeclareMathOperator{\FL}{FL}
\DeclareMathOperator{\rad}{rad}

\DeclareMathOperator{\Gr}{Gr}
\DeclareMathOperator{\KS}{KS}

\newcommand{\n}[1]{\mathcal{#1}}

\newcommand{\bb}[1]{\mathbb{#1}}
\newcommand{\ob}[1]{\mathrm{#1}}
\newcommand{\bbf}[1]{\mathbf{#1}}
\newcommand{\f}[1]{\mathfrak{#1}}
\newcommand{\s}[1]{\mathscr{#1}}

\newcommand{\OBdr}{\s{O}\!\bb{B}_{\dR}}
\newcommand{\OBdR}[1]{\s{O}\!\bb{B}_{\dR, {#1}}}
\newcommand{\OC}{\s{O}\!\bb{C}}
\newcommand{\OCC}[1]{\s{O}\!\bb{C}_{\log,{#1}}}

\DeclareMathOperator{\an}{\begin{scriptsize}
an
\end{scriptsize}}

\DeclareMathOperator{\proet}{\begin{scriptsize}
pro\acute{e}t
\end{scriptsize}}

\DeclareMathOperator{\GM}{\begin{scriptsize}
GM
\end{scriptsize}}
\DeclareMathOperator{\HT}{\begin{scriptsize}
HT
\end{scriptsize}}

\title{A {J}acquet-{L}anglands functor for $p$-adic locally analytic representations}

\author{ {Gabriel} {Dospinescu}}
\author{ {Juan Esteban}   {Rodr\'iguez Camargo}}

\makeatletter
\def\@tocline#1#2#3#4#5#6#7{\relax
  \ifnum #1>\c@tocdepth 
  \else
    \par \addpenalty\@secpenalty\addvspace{#2}%
    \begingroup \hyphenpenalty\@M
    \@ifempty{#4}{%
      \@tempdima\csname r@tocindent\number#1\endcsname\relax
    }{%
      \@tempdima#4\relax
    }%
    \parindent\z@ \leftskip#3\relax \advance\leftskip\@tempdima\relax
    \rightskip\@pnumwidth plus4em \parfillskip-\@pnumwidth
    #5\leavevmode\hskip-\@tempdima
      \ifcase #1
       \or\or \hskip 1em \or \hskip 2em \else \hskip 3em \fi%
      #6\nobreak\relax
    \dotfill\hbox to\@pnumwidth{\@tocpagenum{#7}}\par
    \nobreak
    \endgroup
  \fi}
\makeatother
\usepackage{hyperref}


\begin{document}

\footnote{ UMR $6620$ du CNRS-Laboratoire de Math\'ematiques Blaise Pascal-Universit\'e de Clermont-Auvergne, Campus des C\'ezeaux
3, place Vasarely, 63178 Aubi\`ere cedex, France and Institute of Mathematics "Simion Stoilow" of the Romanian Academy, 21 Calea Grivitei Street, 010702 Bucharest, Romania, 
 gabriel.dospinescu@uca.fr
 
 Max Planck Institute for Mathematics, Room 208, Vivatsgasse 7, 53111 Bonn, Germany, rodriguez@mpim-bonn.mpg.de} 

\maketitle

\begin{abstract}
We study the locally analytic theory of infinite level local Shimura varieties. As a main result, we prove that in the case of a duality of local Shimura varieties, the locally analytic vectors of different period sheaves at infinite level are independent of the actions of the $p$-adic Lie groups $G$ and $G_b$ of the two towers; this  generalizes a result of Pan for the Lubin-Tate and Drinfeld spaces for $\GL_2$. We apply this theory to show that Scholze's $p$-adic  Jacquet-Langlands functor  commutes with the passage to locally analytic vectors, and  is compatible with central characters of Lie algebras. We also prove that the compactly supported de Rham cohomology of the two towers are isomorphic as smooth representations of $G\times G_b$.
\end{abstract}

\tableofcontents

\section{Introduction}
\label{Intro}

 The main objective of this work is to give some new insights concerning the locally analytic part of the $p$-adic local Langlands and Jacquet-Langlands correspondences. We are interested in the infinite level Lubin-Tate and Drinfeld spaces (or more generally local Shimura varieties) and the locally analytic vectors for the action of the associated $p$-adic Lie groups on various period sheaves that are classical in $p$-adic Hodge theory.
 This study is used to establish a compatibility between the Banach and locally analytic $p$-adic Jacquet-Langlands correspondence, as well as 
 comparison results for the compactly supported de Rham cohomology of dual local Shimura varieties (such results are established in various degrees of generality in \cite{MR4073863}, \cite{FlipFlop}, \cite{BenchaoSudeRham} and are probably best understood using the new framework introduced in \cite{deRhamStacksFF}). Along the way, we establish analogs of the 
 results in \cite{RCLocAnCompleted} for local Shimura varieties and we answer a question of Lue Pan.

  To motivate our main results let us briefly recall the perfectoid geometry of the Lubin-Tate and Drinfeld spaces, which are the most classical and well-studied examples of local Shimura varieties \cite{MR3271247}, \cite{ScholzeWeinspadicgeometry}. 
 Let $F$ be a finite extension of $\bb{Q}_p$, with ring of integers $\n{O}$ and residue field $\bb{F}$, and let $\breve{F}$ be the completion of the maximal unramified extension of 
 $F$. Let $G={\rm GL}_n(F)$ and let $D$ be the central division algebra over $F$ with invariant $1/n$. The Lubin-Tate space at level $\GL_n(\n{O})$ is the rigid generic fiber $\n{LT}_{\breve{F}}$ of the formal scheme classifying deformations (by quasi-isogenies) of a formal $\n{O}$-module $\bb{X}_b$
 over $\overline{\bb{F}}$ (algebraic closure of $\bb{F}$) of dimension $1$ and $\n{O}$-height $n$. It is a countable disjoint union of open unit discs of dimension $n-1$ and has a natural (but very complicated) action of $D^{\times}$, as well as an \'etale Grothendieck-Messing period map $\n{LT}_{\breve{F}}\to \bb{P}^{n-1}_{\breve{F}}$. 
 The space $\n{LT}_{\breve{F}}$ has a \textit{dual} given by the Drinfeld space $\Omega_{\breve{F}}\subset \bb{P}^{n-1}_{\breve{F}}$, complement  of the $F$-rational hyperplanes of $\bb{P}^{n-1}_{\breve{F}}$. It has a natural (and fairly simple) action of $G$.

The spaces $\n{LT}_{\breve{F}}$ and $\Omega_{\breve{F}}$ are intimately related via perfectoid geometry in a very clean way: by work of Faltings, Fargues and Scholze-Weinstein (see  \cite{MR1936369},  \cite{MR2441312}, \cite{Scholze2013moduli}) 
 there is a perfectoid space $\n{LT}_{\infty,\breve{F}}$ endowed with an action of $G\times D^{\times}$ and a diagram
\begin{equation}\label{eqDiagramDualityLTIntro}
\begin{tikzcd}
 & \n{LT}_{\infty,\breve{F}}  \ar[rd,"{\pi_{\HT}}"] \ar[ld, "{\pi_{\GM}}"']& \\
\bb{P}^{n-1}_{\breve{F}} & & \Omega_{\breve{F}} 
\end{tikzcd}
\end{equation}
such that: 

\begin{itemize}

\item $\n{LT}_{\infty, \breve{F}}$ is the Lubin-Tate space at infinite level, obtained by trivializing the Tate module of the universal deformation of $\bb{X}_b$. It is a pro\'etale $\GL_n(\n{O})$-torsor over $\n{LT}_{\breve{F}}$, and $\pi_{\GM}: \n{LT}_{\infty,\breve{F}}\to \bb{P}^{n-1}_{\breve{F}}$ is a pro\'etale $G$-torsor, obtained by composing $\n{LT}_{\infty, \breve{F}}\to \n{LT}_{\breve{F}}$ with the Grothendieck-Messing period map $\n{LT}_{\breve{F}}\to \bb{P}^{n-1}_{\breve{F}}$.

\item $\pi_{\GM}$ is a $D^{\times}$-equivariant $G$-torsor, $D^{\times}$ acting on $\bb{P}^{n-1}_{\breve{F}}$ via its embedding into $\GL_n(\breve{F})$. 

\end{itemize}

 The above diagram carries a suitable Weil descent datum to $F$. Thus, its base change to $\bb{C}_p$ carries an action of the Weil group $W_F$. 
In particular, there is an action of $G\times D^{\times} \times W_F$ on $\n{LT}_{\infty, \bb{C}_p}$, which makes natural the expectation that both the Jacquet-Langlands correspondence, relating $G$ and $D^{\times}$-representations,  and the Langlands correspondence, relating $G$ and $W_F$-representations, can be realised in various cohomologies attached to $\n{LT}_{\infty, \bb{C}_p}$, see for example \cite{MR2441312}, \cite{MR1876802},  \cite{MR3702670},  \cite{MR4073863},  \cite{MR4255044}, \cite{MR4595286}.

In \cite{ScholzeLubinTate}, Scholze used the diagram \eqref{eqDiagramDualityLTIntro} to construct a \textit{$p$-adic Jacquet-Langlands functor} sending smooth admissible $p$-torsion representations of $G$ to smooth admissible $p$-torsion representations of $D^{\times}$: if $\pi$ is an admissible smooth representation of $G$ over an artinian ring $A$ which is $p$-power torsion (or a complex of such representations), we obtain a $D^{\times}$-equivariant \'etale sheaf $\n{F}_{\pi}$ on the rigid space $\bb{P}^{n-1}_{\bb{C}_p}$ using that 
 $\pi_{\GM}$ is a pro\'etale $G$-torsor, and a complex of smooth  $D^{\times}$-representations over $A$
\[
\n{JL}(\pi)= R\Gamma_{\et}(\bb{P}^{n-1}_{\bb{C}_p}, \n{F}_{\pi}).
\]
By \cite[Theorem 1.1]{ScholzeLubinTate} the cohomology groups $\n{JL}^i(\pi)$ of $\n{JL}(\pi)$ are smooth admissible representations of  $D^{\times}$ over $A$ and $\n{JL}^i(\pi)=0$ for $i>2(n-1)$. Moreover, this construction satisfies a local-global compatibility for $n=2$ (see 
\cite[Theorem 1.3]{ScholzeLubinTate}; similar results should hold for general $n$), justifying the compatibility with a more classical Jacquet-Langlands correspondence.

One can naturally extend the Jacquet-Langlands functor to unitary Banach representations and it is not hard to see that it also preserves admissible Banach representations, see \cref{CoroBanachAdmissible} (the situation is more delicate for non unitary admissible Banach representations, as then such a naive statement is not true anymore).  On the other hand, Schneider-Teitelbaum introduced a class of admissible locally analytic representations for $p$-adic Lie groups in  \cite{SchTeitDist}. A natural question arises:
\begin{question}\label{eqQuestion1Intro}
Is there a Jacquet-Langlands functor $\pi\mapsto \n{JL}(\pi)$ for admissible locally analytic representations of $\GL_n(F)$?   If so, is it compatible with the Jacquet-Langlands functor for Banach representations? 
\end{question}
In this paper we give a partial answer to this question, namely, that the Jacquet-Langlands functor for admissible unitary Banach representations is compatible with the passage to locally analytic vectors, see \cref{Theo8Intro} for a more precise statement.

In a different direction,  the works of Lue Pan \cite{LuePanI,LuePanII} studying the locally analytic vectors of perfectoid  modular curves use some special sheaves of locally analytic functions at infinite level. These sheaves encode, via the localization theory of Beilinson-Bernstein  on the flag variety \cite{BeilinsonBernstein} and the Hodge-Tate period map, many aspects of the $p$-adic Hodge theory of Shimura varieties. In  \cite{RCGeoSen,RCLocAnCompleted} some of these features have been generalized to arbitrary global Shimura varieties. Another aspect of this paper is to extend the results in geometric Sen theory from the global to the local setup. It is then natural to ask what additional properties local Shimura varieties acquire after taking locally analytic vectors, in particular one can ask the following question:

\begin{question}\label{Question2Intro}
   Are the $G$-locally analytic vectors of the structure sheaf $\widehat{\s{O}}_{\n{LT}}$ of $\n{LT}_{\infty, \bb{C}_p}$ equal to its $D^{\times}$-locally analytic vectors? More precisely, if $U_{\infty}\subset \n{LT}_{\infty, \bb{C}_p}$ is an affinoid perfectoid with (open) stabilizers 
   $K_U\subset G$ and $K_{U,D}\subset D^{\times}$, do we have 
\[
\widehat{\s{O}}_{\n{LT}}(U_{\infty})^{K_U-{\rm la}}=\widehat{\s{O}}_{\n{LT}}(U_{\infty})^{K_{U,D}-{\rm la}}
\]
as subspaces of $\widehat{\s{O}}_{\n{LT}}(U_{\infty})$?
\end{question}

For  $n=2$ this is proved by Pan in \cite[Corollary 5.3.9]{LuePanII} via explicit power series expansions. He added (see Remark 5.3.10 in loc.cit.) that "it will be very interesting and conceptually satisfying to have a more intrinsic proof". 
We answer this question and prove a much more general result that holds for an arbitrary duality of local Shimura varieties and arbitrary period sheaves appearing in the affinoid charts of relative Fargues-Fontaine curves, see \cref{CoroLocAnVectorInfiniteShimura} for a precise statement. Then, the partial result towards \cref{eqQuestion1Intro} mentioned above will be a rather formal consequence of this independence of locally analytic vectors at infinite level, after applying enough technology coming from the theory of solid locally analytic representations \cite{RRLocallyAnalytic,RJRCSolidLocAn2}. We also apply this independence of locally analytic vectors to  construct an equivariant isomorphism for the compactly supported de Rham cohomology between the two towers of a duality of local Shimura varieties, see \cref{Theo7Intro}.

 In order to present the main results of this paper we have separated the introduction in different paragraphs, going from the  general results on towers of rigid spaces, passing to the applications to local Shimura varieties, and finishing with the most specific applications to  the Lubin-Tate and Drinfeld  towers.

\subsubsection*{Cohomology of towers of rigid spaces} In this paragraph we explain the results of \cref{SubsecLocAnTower} about locally analytic vectors of period sheaves in towers of rigid spaces. Let $\Perfd$ be the category of perfectoid spaces and let $\Perf\subset \Perfd$ be the full subcategory of perfectoid spaces in characteristic $p$, seen as sites for the $v$-topology introduced in 
 \cite{scholze2022etale}. There are structural sheaves $\widehat{\s{O}}$ and $\widehat{\s{O}}^{+}$, as well as a Fontaine period sheaf $\bb{A}_{\inf}$
 mapping an affinoid perfectoid $\Spa(R,R^+)\in \Perfd$ to $R$, $R^{+}$ and $W(R^{\flat,+})$ respectively.  
 Let $\Perf_{\varpi}$ be the category of perfectoid spaces in characteristic $p$ endowed with a fixed pseudo-uniformizer $\varpi$, and with maps preserving the pseudo-uniformizer and let $I=[s,r]\subset (0,\infty)$ be a compact interval with rational ends. One obtains a period sheaf 
 $\bb{B}_I$ on $\Perf_{\varpi}$ by sending an affinoid perfectoid $\Spa(R,R^+)\in \Perf_{\varpi}$ in characteristic $p$ with fixed pseudo-uniformizer to the rational localization 
\[
\bb{B}_I(S)= W(R^+) \left\langle \frac{p}{[\varpi]^{1/r}}, \frac{[\varpi]^{1/s}}{p} \right\rangle \left[\frac{1}{[\varpi]}\right],
\]
in other words $\bb{B}_I(S)$ is the ring of analytic functions on the rational subspace of ${\rm Spa}(W(R^+))$ given by $$|p|^r\leq |[\varpi]|\leq |p|^s\ne 0.$$

     Let $K/\qp(\mu_{p^{\infty}})$ be a perfectoid field and let $X$ be a qcqs smooth rigid space over $K$ endowed with a continuous action of a compact
     $p$-adic Lie group $H$. Let $G$ be another compact $p$-adic Lie group and consider an $H$-equivariant pro-finite-\'etale $G$-torsor $\widetilde{X}\to X^{\lozenge}$ of the diamond attached to $X$. In particular, $\widetilde{X}$ is endowed with an action of the $p$-adic Lie group $G\times H$. The following theorem relates the locally analytic vectors of the $v$-cohomologies of period sheaves at infinite level. It is really the key new idea of this article.

\begin{Theo}[\cref{theo:LocAnTwoTowers}]\label{Theo1Intro}
Let $I\subset (0,\infty)$ be a compact interval with rational ends.  Then the $G$-locally analytic  vectors of the solid $\bb{Q}_p$-linear representation $R\Gamma_v(\widetilde{X}, \bb{B}_I)$ are $H$-locally analytic. More precisely,  the natural map of solid $G\times H$-representations
\[
R\Gamma_v(\widetilde{X}, \bb{B}_I)^{RG\times H-{\rm la}}\xrightarrow{\sim} R\Gamma_v(\widetilde{X}, \bb{B}_I)^{RG-{\rm la}}
\]
is an equivalence. 
\end{Theo}

 Taking $G$ trivial, it follows that pro\'etale cohomology groups of qcqs smooth rigid varieties $X$ endowed with continuous actions of $p$-adic Lie groups $H$ are locally analytic: 

\begin{Coro}[\cref{CohoEquivRigidSpaces}]\label{CoroIntro2}
Let $X$ be a  qcqs smooth rigid space over $K$, endowed with a continuous action of a $p$-adic Lie group $H$. Then for $I\subset (0,\infty)$ a compact interval with rational ends the solid $H$-representation $R\Gamma_v(X, \bb{B}_I)$ is $H$-locally analytic. Moreover, the solid 
$H$-representation $R\Gamma_{\rm et}(X, \qp)$ is locally analytic.
\end{Coro}

\begin{Rem}\label{RemIntro2}
In particular, the cohomology groups of period sheaves on $X$ admit an action of the Lie algebra of $H$ obtained by derivations. We found this fact surprising since there is no finiteness or Hausdorff assumptions in the cohomology groups (for instance, it is not known, even when $X$ is affinoid, whether the \'etale cohomology groups $H^i(X, \qp)$ are Banach spaces). This also suggests that there is a deeper structure in the period sheaves of rigid spaces that witnesses the locally analytic properties of their cohomologies. In a work in progress of Johannes Ansch\"utz, Arthur-C\'esar le Bras, Peter Scholze and the second author we give a conceptual explanation of these facts via the analytic prismatization. 
\end{Rem}

\subsubsection*{Geometric Sen theory over local Shimura varieties} 
 We state now the main results of \cref{Section:SenOperators} about the Sen operators of local Shimura varieties (see \cite{LuePanI} and \cite{RCLocAnCompleted} for the case of global Shimura varieties). Let $(\bbf{G},[b],\{\mu\})$ be a local Shimura datum (as in Lecture XXIV of \cite{ScholzeWeinspadicgeometry}), so $\bbf{G}$ is a connected reductive group over $\qp$, $\{\mu\}$ is a conjugacy class of geometric minuscule cocharacters of $\bbf{G}$ and $[b]$ is an isomorphism class of $G$-isocrystals. Let $E$ be the local reflex field, i.e. 
 the field of definition of $\{\mu\}$. For a compact open subgroup $K\subset \bbf{G}(\bb{Q}_p)$ let $$\n{M}_{K}=\n{M}_{\bbf{G},b,\mu,K,\breve{E}}$$ be the corresponding local Shimura variety of level $K$, a smooth rigid analytic space over $\breve{E}$, the completion of the maximal unramified extension of $E$.

 Let $\FL_{\bbf{G},\mu,E}$ and $\FL_{\bbf{G},\mu^{-1},E}$ be the algebraic flag varieties classifying decreasing and increasing $\mu$-filtrations of the trivial $\bbf{G}$-torsor respectively, and let $\Fl_{\bbf{G},\mu,E}$ and $\Fl_{\bbf{G},\mu^{-1},E}$ be their analytification, seen as adic spaces over ${\rm Spa}(E)$. 
 The corresponding infinite level local Shimura variety 
 $$\n{M}_{\infty}=\n{M}_{\bbf{G},b,\mu,\infty,\breve{E} }:=\varprojlim_{K} \n{M}_{\bbf{G},b,\mu,K,\breve{E}}^{\lozenge}$$
 is a diamond over $\breve{E}$ endowed with 
  Grothendieck-Messing and Hodge-Tate period maps
\[
\begin{tikzcd}
 & \n{M}_{ \infty } \ar[rd,"\pi_{\HT}"] \ar[ld,"\pi_{\GM}"']&  \\ 
  \Fl_{\bbf{G},\mu,\breve{E}}^{\lozenge}& & \Fl_{\bbf{G},\mu^{-1}, \breve{E}}^{\lozenge}. 
\end{tikzcd}
\]
and with an action of $\bbf{G}(\bb{Q}_p)\times \widetilde{G}_b$ for which both maps $\pi_{\GM}$ and $\pi_{\HT}$ are equivariant in a suitable sense 
(see \cref{SectionLocalShimura} for more details). Here 
 $\widetilde{G}_b$ is the group of automorphisms of the $\bbf{G}$-torsor $\n{E}_b$ over the Fargues-Fontaine curve attached to $[b]$, see \cite[III.5.1]{FarguesScholze}. 
By \cite[Corollary 23.3.2]{ScholzeWeinspadicgeometry}, when $b$ is basic there is a dual local Shimura datum $(\Check{\bbf{G}}, \Check{b}, \Check{\mu})$, an isomorphism $G_b=\widetilde{G}_b= \Check{\bbf{G}}(\bb{Q}_p)$  and a $\bbf{G}(\bb{Q}_p)\times \Check{\bbf{G}}(\bb{Q}_p)$-equivariant isomorphism of infinite level local Shimura varieties
\[
\n{M}_{\bbf{G},b,\mu,\infty,\breve{E}} \cong \n{M}_{\Check{\bbf{G}},\Check{b},\Check{\mu},\infty,\breve{E}} 
\]
that exchanges the Grothendieck-Messing and Hodge-Tate period maps.

 Fix a compact open subgroup $K$ of $\bbf{G}(\qp)$. 
The map $$\pi_K:\n{M}_{\bbf{G},b,\mu,\infty,\breve{E}}\to \n{M}_{\bbf{G},b,\mu,K,\breve{E}}^{\lozenge}$$ is a pro\'etale  $K$-torsor, in particular:

\begin{itemize}

\item  For any Ind-system $V=``\varinjlim_i" V_i$ of $p$-adically complete continuous representations of $K$ we can construct a $v$-sheaf
$$\n{F}_{V}=\varinjlim_{i} \n{F}_{V_{i}}$$
on $\n{M}_{K}^{\lozenge}$, where 
 $\n{F}_{V_i}$ is the $p$-complete $v$-sheaf obtained by descent along $\pi_K$ of the "constant" $v$-sheaf $\underline{V_i}$.  In particular, for $V$ an algebraic representation of $\bbf{G}$ (for instance the dual $\f{g}^{\vee}$ of the adjoint representation of $\bbf{G}$)  we have automorphic local systems $\n{F}_{V}$, and for $\pi$ a smooth admissible representation of $\bbf{G}(\bb{Q}_p)$ over a $p$-power torsion ring $A$ the sheaf $\n{F}_{\pi}$ is the \'etale local system considered in \cite{ScholzeLubinTate} for the $p$-adic Jacquet-Langlands functor (in the special case of the infinite level Lubin-Tate space).

\item There is a canonical geometric Sen operator (\cite[Theorem 3.3.4]{RCGeoSen}) 
$$\theta_{\n{M}_K}: \n{F}_{\f{g}^{\vee}}\otimes_{\bb{Q}_p} \widehat{\s{O}} \to \Omega^1_{\n{M}_K} \otimes_{\s{O}_{\n{M}_K}} \widehat{\s{O}}(-1),$$
where $\Omega^1_{\n{M}_K}$ is the cotangent bundle of the smooth rigid space $\n{M}_K$ over $\breve{E}$ and $\s{F}(n)$ is the $n$-th Hodge-Tate twist of $\s{F}$ by the $n$-th power of the cyclotomic character. 

\end{itemize}

  Before stating the main theorem of this section we need to introduce more notation, coming from the localization theory of Beilinson-Bernstein \cite{BeilinsonBernstein}. Let $E'$ be a finite extension of $E$ that splits $\bbf{G}$ and fix once and for all a 
 cocharacter $\mu: \bb{G}_{m,E'}\to \bbf{G}_{E'}$ in the conjugacy class $\{\mu\}$.
   Let $\bbf{P}_{\mu}$ and $\bbf{P}_{\mu^{-1}}$ be the (opposite) parabolic subgroups classifying decreasing and increasing filtrations of type $\mu$ on algebraic representations of $\bbf{G}_{E'}$, let $\bbf{N}_{\mu}\subset \bbf{P}_{\mu}$ and $\bbf{N}_{\mu^{-1}}\subset \bbf{P}_{\mu^{-1}}$ be their unipotent radicals respectively, and let $\bbf{M}=\bbf{M}_{\mu}=\bbf{M}_{\mu^{-1}}$ be their common Levi subgroup, the centralizer of $\mu$ in $\bbf{G}_{E'}$. The isomorphisms $$\FL_{\bbf{G},\mu,E'}\simeq \bbf{G}_{E'}/\bbf{P}_{\mu}, \quad \FL_{\bbf{G},\mu^{-1},E'}\simeq \bbf{G}_{E'}/\bbf{P}_{\mu^{-1}}$$ 
    induce equivalences between $\bbf{G}$-equivariant quasi-coherent sheaves on $\FL_{\bbf{G},\mu,E'}$ and $\FL_{\bbf{G},\mu^{-1},E'}$ and algebraic representations of $\bbf{P}_{\mu}$ and $\bbf{P}_{\mu^{-1}}$ respectively. In particular we obtain 
    $\bbf{G}$-equivariant Lie algebroids $$\f{n}^0_{\mu}\subset \f{p}^0_{\mu}\subset \f{g}^0_{\mu}=  \s{O}_{\FL_{\bbf{G},\mu,E'}}\otimes_{\bb{Q}_p} \f{g}$$  and $\f{m}^0_{\mu}= \f{p}^0_{\mu}/\f{n}^0_{\mu}$ attached to the adjoint action of $\bbf{P}_{\mu}$ on the Lie algebras 
 $\f{n}_{\mu}\subset  \f{p}_{\mu} \subset \f{g}$ and $\f{m}_{\mu}= \f{p}_{\mu}/\f{n}_{\mu}$ 
 of $\bbf{N}_{\mu}\subset \bbf{P}_{\mu}\subset \bbf{G}_{E'}$ and $\bbf{M}_{\mu}$. We denote in the same way their pullbacks to vector bundles over the analytic flag varieties, and use similar notation for the Lie algebras of the opposite parabolic and their associated Lie algebroids in $\FL_{\bbf{G},\mu^{-1}, E'}$.  

The next theorem is an analogue for local Shimura varieties of \cite[Theorem 5.2.5]{RCLocAnCompleted}.

\begin{Theo}[\cref{TheoremComputationGeoSenOperator}]\label{TheoIntro3}
  The  Kodaira-Spencer isomorphism (up to a sign)  induces an identification 
$$\pi^*_{\HT}(\f{n}^{0,\vee}_{\mu^{-1}})\cong \Omega^1_{\n{M}_K} \otimes_{\s{O}_{\n{M}_K}} \widehat{\s{O}}(-1),$$
which in turn induces a canonical identification 
$$\pi_K^*(\theta_{\n{M}_K})=\pi_{\HT}^*(\f{g}^{0,\vee}\to \f{n}^{0,\vee}_{\mu^{-1}}),$$
the map $\f{g}^{0,\vee}\to \f{n}^{0,\vee}_{\mu^{-1}}$ being dual to the inclusion $\f{n}^0_{\mu^{-1}}\subset \f{g}^0$.
 \end{Theo}

\begin{Rem}\label{RemarkIntro3}
  The flag varieties and $\pi_{\HT}$ are defined over $E$, and so are the Lie algebroids $\f{n}^0_{\mu^{-1}}$, $\f{p}^0_{\mu^{-1}}$ and $\f{m}^0_{\mu^{-1}}$. 
  In the previous theorem one can avoid to base change the flag varieties and local Shimura varieties to $E'$, by avoiding to fix a Hodge cocharacter $\mu$.
\end{Rem}

A first consequence of the computation of the geometric Sen operator is the vanishing of the higher locally analytic vectors of the structural sheaf $\widehat{\s{O}}$ at infinite level, as well as the computation of the arithmetic Sen operator in terms of geometric representation theory. Let $C/E'$ be a complete algebraically closed extension and consider the $C$-base change $\n{M}_{K,C}$ of $\n{M}_K$. Let $\n{V}_K=C^{{\rm la}}(K,\bb{Q}_p)$ be the space of locally analytic functions on $K$ endowed with the left regular action and let $\n{F}_{\n{V}_K}$ be the associated $v$-sheaf on $\n{M}_{K,C}$. The following theorem is an analogue for local Shimura varieties of Proposition 6.2.8, Corollary 6.2.13 and Theorem 7.2.1 of \cite{RCLocAnCompleted}.

\begin{Theo}[{\cref{TheoVanishingLocAnVectors}}]\label{Theo4Intro}
Let $U\subset \n{M}_{K,C}$ be an open affinoid subspace with toric coordinates\footnote{As usual, this means that 
$U$ has an \'etale map to a torus $\bb{T}^d_C$, which is a finite composition of rational localizations and finite \'etale maps.}, with pullback $U_{\infty}\subset \n{M}_{\infty,C}$. Then\footnote{Here the completed tensor product is a filtered colimit of $p$-completed tensor products obtained by writing $\n{F}_{\n{V}_K}$ as a colimit of Banach sheaves (equivalently a solid tensor product as in \cite{anschutz2024descent}).}

\begin{equation}\label{equationVanishingLocAnVectors}
R\Gamma_{v}(U, \widehat{\s{O}}\widehat{\otimes}_{\bb{Q}_p}\n{F}_{\n{V}_K})=\widehat{\s{O}}(U_{\infty})^{G-{\rm la}},
\end{equation}
in particular this complex is concentrated in degree $0$. 

  Moreover, the subsheaf $\f{n}^0_{\mu^{-1}}$ of 
$\f{g}^0_{\mu^{-1}}$ kills $\widehat{\s{O}}(U_{\infty})^{G-{\rm la}}$, inducing a horizontal action of  $\f{m}^0_{\mu^{-1}}$ on $\widehat{\s{O}}(U_{\infty})^{G-{\rm la}}$, and $\widehat{\s{O}}(U_{\infty})^{G-{\rm la}}$  has an arithmetic Sen operator as in \cite[Definition 7.1.2]{RCLocAnCompleted} given by the negative of the derivative of the Hodge cocharacter $-\theta_{\mu}=\theta_{\mu^{-1}}\in \f{m}^0_{\mu^{-1}}$.
\end{Theo}

\subsubsection*{Locally analytic vectors of local Shimura varieties}

In this paragraph we apply \cref{Theo1Intro} to prove the independence of locally analytic vectors for a duality of local Shimura varieties, generalizing a theorem of Pan  for the Lubin-Tate tower of $\GL_2$ \cite[Corollary 5.3.9]{LuePanII}.

 Let $C/E$ be  a complete algebraically closed field, let $G=\bbf{G}(\bb{Q}_p)$ and let $G_b$ be the profinite quotient of $\widetilde{G}_b$. 
 Let $\widehat{\s{O}}_{\n{M}}$ be the restriction of the structural sheaf $\widehat{\s{O}}$ in the $v$-site of $\n{M}_{\infty,C}$ to the underlying topological space $| \n{M}_{\infty,C}|$. Let $\s{O}^{G-{\rm la}}_{\n{M}}\subset \widehat{\s{O}}_{\n{M}}$ be the subsheaf whose values in a qcqs open subspace $U_{\infty}\subset  \n{M}_{\infty,C}$ are given by  the $G=\bbf{G}(\bb{Q}_p)$-locally analytic sections of $\widehat{\s{O}}_{\n{M}}$, more precisely 
\[
\s{O}^{G-{\rm la}}_{\n{M}}(U_{\infty}) = \widehat{\s{O}}(U_{\infty})^{K_{U_{\infty}}-{\rm la}}
\]
where $K_{U_{\infty}}\subset G$ is the stabilizer of $U_{\infty}$. 

\begin{Coro}[{\cref{TheoMainComparisonLocAn}}]\label{CoroLocAnVectorInfiniteShimura}
For any $p$-adic Lie group $H\subset \widetilde{G}_b$ and any qcqs open subspace $U_{\infty}\subset \n{M}_{\infty,C}$ with stabilizers $G_U\subset G$ and $H_U\subset H$, the natural map 
\[
\s{O}^{G-{\rm la}}_{\n{M}}(U_{\infty})^{RH_U-{\rm la}}\xrightarrow{\sim} \s{O}^{G-{\rm la}}_{\n{M}}(U_{\infty})
\] 
is an equivalence. In particular, if $b$ is basic  we have an equality of subsheaves of $\widehat{\s{O}}_{\n{M}}$
\[
\s{O}^{G-{\rm la}}_{\n{M}}= \s{O}^{G_b-{\rm la}}_{\n{M}}. 
\]
More generally, for $b$ basic and   $I\subset (0,\infty)$ a compact interval with rational ends,  we have  natural equivalences of derived solid locally analytic representations of $G\_Utimes G_{b,U}$
\[
R\Gamma_{v}(U_{\infty}, \bb{B}_{I})^{RG_{b,U}-{\rm la}} \xleftarrow{\sim} R\Gamma_{v}(U_{\infty}, \bb{B}_{I})^{R G_U\times G_{b,U}-{\rm la}}  \xrightarrow{\sim } R\Gamma_{v}(U_{\infty}, \bb{B}_{I})^{RG_U-{\rm la}}. 
\]
\end{Coro}

From now on we shall focus in the case when $b$ is basic. We will identify the $G$ and $G_b$-locally analytic vectors of the structural  sheaf $\widehat{\s{O}}_{\n{M}}$ at infinite level and simply write $\s{O}^{la}_{\n{M}}$. We will explain now how to identify the horizontal actions of the Levi Lie algebras of \cref{Theo4Intro}. By \cref{CorollaryIsoTorsorsShimuraVariety} we have a natural $G\times G_b$-equivariant isomorphism of $\widehat{\s{O}}_{\n{M}}$-vector bundles on $\n{M}_{\infty,C}$
 \[
 \f{m}^0_{\mu^{-1}}\otimes_{\s{O}_{\Fl_{\bbf{G},\mu^{-1},C}}} \widehat{\s{O}}_{\n{M}} \cong  \widehat{\s{O}}_{\n{M}} \otimes_{\s{O}_{\Fl_{\bbf{G},\mu,C}}} \f{m}^0_{\mu}.
 \]  
By taking locally analytic vectors we obtain an $\s{O}^{{\rm la}}_{\n{M}}$-vector bundle which we shall denote as $\f{m}^{0, {\rm la}}$, endowed with $G\times G_b$-equivariant isomorphisms
\begin{equation}\label{eqIsoSheavesLeviIntro}
\f{m}^0_{\mu^{-1}} \otimes_{\s{O}_{\Fl_{\bbf{G},\mu^{-1},C}}} \s{O}^{{\rm la}}_{\n{M}} \cong \f{m}^{0, {\rm la}} \cong  \s{O}^{ {\rm la}}_{\n{M}}\otimes_{\s{O}_{\Fl_{\bbf{G},\mu,C}}} \f{m}^0_{\mu}.
\end{equation}
We have the following theorem:

 \begin{Theo}[{\cref{TheoCenterBothSides}}]\label{Theo6Intro}
The actions of $\f{n}_{\mu}^0$ and $\f{n}_{\mu^{-1}}^0$ on $\s{O}^{ {\rm la}}_{\n{M}}$ vanish, and the actions of $\f{m}^0_{\mu}$ and $\f{m}^0_{\mu^{-1}}$ on $\s{O}^{{\rm la}}_{\n{M}}$ by derivations are identified via \eqref{eqIsoSheavesLeviIntro}. In particular, the central character of the actions of $\f{m}^0_{\mu}$ and $\f{m}^0_{\mu^{-1}}$ on $\s{O}^{{\rm la}}_{\n{M}}$ agree under the natural isomorphism of the center of the enveloping algebras $\n{Z}(\f{m}_{\mu,C}) \cong \n{Z}(\f{m}_{\mu^{-1},C})$, where $\f{m}_{\mu,C}$ and $\f{m}_{\mu^{-1},C}$ are the Levi subalgebras of  $\Lie G_b\otimes_{\bb{Q}_p} C $ and $\Lie G\otimes_{\bb{Q}_p} C$  respectively.
\end{Theo}

\subsubsection*{De Rham cohomology of towers of local Shimura varieties}

Our next result is the comparison between compactly supported de Rham cohomologies of the two towers in  a duality of local Shimura varieties. This  theorem has been also independently obtained by Wies\l{}awa Nizio\l{} and the first author \cite{FlipFlop} (with weaker functoriality properties) and Benchao Su \cite{BenchaoSudeRham} (for the Lubin-Tate and Drinfeld spaces). Such results 
are probably best understood now in terms of de Rham stacks and descent properties of them, as developed in \cite{camargo2024analytic} and \cite{deRhamStacksFF}. As pointed out by Arthur-C\'esar le Bras,  purely motivic techniques as those appearing in \cite{MR3880025} (see Proposition 4.5 in \textit{loc. cit.}) also yield such a result, concretely, it is a corollary of the de Rham realization of Berkovich motives obtained from the thesis of Aoki \cite{AokiBerkovichMotives}.

  Let $(\bbf{G},b,\mu)$ be a local Shimura datum with $b$ basic and let $(\Check{\bbf{G}},\Check{b},\Check{\mu})$ be the dual local Shimura datum, giving rise to towers of rigid spaces $(\n{M}_{\bbf{G}, b,\mu,K, \breve{E}})_{K\subset \bbf{G}(\bb{Q}_p)}$ and $(\n{M}_{\Check{\bbf{G}}, \Check{b},\Check{\mu},\Check{K}, \breve{E}})_{K\subset \Check{\bbf{G}}(\bb{Q}_p)}$

\begin{Theo}[\cref{TheoComparisonCoho}]\label{Theo7Intro}
There is a natural $\bbf{G}(\bb{Q}_p)\times \Check{\bbf{G}}(\bb{Q}_p)$-equivariant isomorphism of compactly supported de Rham cohomology groups
\[
\varinjlim_{K\subset \bbf{G}(\bb{Q}_p)} H^i_{dR,c}(\n{M}_{\bbf{G},b,\mu,K,C})\cong  \varinjlim_{\Check{K}\subset \Check{\bbf{G}}(\bb{Q}_p)} H^i_{dR,c}(\n{M}_{\Check{\bbf{G}}, \Check{b},\Check{\mu}, \Check{K},C}).
\]
\end{Theo}

The strategy is to relate the de Rham complexes of each tower with a suitable de Rham complex of the sheaf $\s{O}^{la}_{\n{M}}$ arising from the derivations of both groups $G$ and $G_b$.

\subsubsection*{Locally analytic Jacquet-Langlands functor in the Lubin-Tate case}

We finish the presentation of the main results with the principal motivation that initiated this project, that is, the  $p$-adic Jacquet-Langlands functor of the  Lubin-Tate tower treated in  \cite{ScholzeLubinTate}. We keep the notation of the beginning of the introduction regarding the Lubin-Tate and Drinfeld towers.  We have the following compatibility with the passage to locally analytic vectors:

\begin{Theo}[\cref{TheoJLLocAn}]\label{Theo8Intro}
Let $\pi$ be an admissible unitary Banach representation of $\GL_n(F)$ and let $\Pi=\pi^{\GL_n(F)-{\rm la}}$, with associated  
 pro\'etale sheaf $\n{F}_{\Pi}$ over $\bb{P}^{n-1}_{\bb{C}_p}$ (constructed by descent along $\pi_{\GM}$ from $\Pi$). There is a natural equivalence of solid locally analytic $D^{\times}$-representations
 \[
 R\Gamma_{\proet}(\bb{P}^{n-1}_{\bb{C}_p}, \n{F}_{\pi})^{RD^{\times}-{\rm la}} \cong  R\Gamma_{\proet}(\bb{P}^{n-1}_{\bb{C}_p}, \n{F}_{\Pi}), 
 \]
which induces an isomorphism of cohomology groups: 
 \[
H^i_{\proet}(\bb{P}^{n-1}_{\bb{C}_p}, \n{F}_{\pi})^{D^{\times}-{\rm la}} \cong H^i_{\proet}(\bb{P}^{n-1}_{\bb{C}_p}, \n{F}_{\Pi}).
 \]
\end{Theo}

The key strategy is to rewrite the pro\'etale cohomology of the sheaf $\n{F}_{\pi}$ in terms of period sheaves $\bb{B}_{I}$ and then to exploit the independence of locally analytic vectors at infinite level in \cref{CoroLocAnVectorInfiniteShimura} in order to jump between towers.

As a corollary of  Theorems \ref{Theo6Intro} and \ref{Theo8Intro} we can show that the Jacquet-Langlands functor preserves central characters for the locally analytic vectors of admissible Banach representations. Such results were first obtained in \cite{MR4940362} and \cite{MR4621880}, by completely different means (using $p$-adic analytic continuation and global methods). In our approach, we rather use the proof of \cref{Theo8Intro} and the compatibility of the horizontal characters for the sheaf $\s{O}^{la}_{\n{M}}$ of \cref{Theo6Intro}. 

\begin{Coro}[{\cref{CoroCentralChar}}]\label{CoroIntroCentralCharJL}
Let $\pi$ be an admissible unitary Banach representation of $\GL_n(F)$ such that $\pi^{\GL_n(F)-{\rm la}}$ has infinitesimal character $\chi$. Then, for all $i\in \bb{Z}$, the locally analytic $D^{\times}$-representation $H^i_{\proet}(\bb{P}^{n-1}_{\bb{C}_p}, \n{F}_{\pi})^{D^{\times}-{\rm la}} $ has infinitesimal character $\chi$ under the natural identification $\n{Z}(\Lie D^{\times})\otimes C \cong \n{Z}(\Lie G)\otimes C$. 
\end{Coro}

\subsection*{Outline of the paper}

 \cref{SectionPreliminaries} is a preliminary section devoted to recollections about the main players and tools in this paper: period sheaves, solid locally analytic representations and geometric Sen theory. The only new result is  \Cref{bplus}, which we found interesting despite the fact that one could prove the main results of this paper without such a strong result, and the crucial local analyticity criterion in  \Cref{PropLocAnCriterion}. Combined with geometric Sen theory, this is the major input in the proof of  \Cref{theo:LocAnTwoTowers}, to which the third section is devoted.  \Cref{SectionLocalShimura}  contains no original result, we briefly present the construction and main properties of local Shimura varieties and their period maps from the point of view of the stack ${\rm Bun}_G$ of Fargues-Scholze. \Cref{Section:SenOperators} is devoted to the computation of the geometric and arithmetic Sen operator of local Shimura varieties. The proof follows fully faithfully that of 
 \cite[Theorems 5.2.5 and 7.2.1]{RCLocAnCompleted}, once a $p$-adic Riemann-Hilbert correspondence for automorphic local systems is established in this setting (which is much easier than the one for global Shimura varieties) and once we have a 
purely representation theoretic construction of the Kodaira-Spencer isomorphism (essentially a reinterpretation of the anchor map of the reductive group acting on the flag variety). This point of view on the Kodaira-Spencer map allows us to compute the pullback of equivariant sheaves on flag varieties via the Hodge-Tate period maps in terms of automorphic vector bundles and the Faltings extension.
  Finally, \Cref{SectionLocAnInfiniteLevel,ssJLFunctor}   use the previously established results to prove the main theorems of the paper, as stated in the introduction.

\subsection*{Conventions}

In this paper we use the $v$-site of perfectoid spaces as introduced  in \cite{scholze2022etale}. We use  the theory of solid almost quasi-coherent sheaves of \cite{MannSix}, \cite{anschutz2024descent} and \cite{Anschutz-LeBras-Mann}; the use of these cohomology theories is important in order to properly keep track to the condensed or topological structure of cohomology complexes.  In particular, this work heavily depends on the theory of condensed mathematics of Clausen and Scholze \cite{ClausenScholzeCondensed2019,ClauseScholzeAnalyticGeometry}, and on higher category theory, for which we refer to \cite{HigherTopos,HigherAlgebra}. A different reason to use condensed mathematics    is to have access to the theory of solid locally analytic representations of \cite{RRLocallyAnalytic,RJRCSolidLocAn2}. This is important since, even though most of the main theorems  involve classical topological representations, the proofs will make appear very general solid representations which are not classical. 

\subsection*{Acknowledgements} 

We thank Johannes Ansch\"utz, Laurent Berger, Guido Bosco,  Arthur-C\'esar le Bras, Pierre Colmez,  Wies\l{}awa Nizio\l{}, Lue Pan  and Peter Scholze for enlightening conversations in different stages of this work. We also thank  Johannes Ansch\"utz, Guido Bosco, Arthur-C\'esar le Bras, Zhenghui Li, Wies\l{}awa Nizio\l{} and Benchao Su for corrections and comments concerning a first draft of the paper. Part of this project was done during the trimester program in Bonn: ``The Arithmetic of The Langlands Program'' during the summer of 2023. We heartily thank the organizers and the Hausdorff Research Institute for Mathematics for the excellent environment for mathematical discussions and  exchanges. The second author wants to thank Columbia University and the Simons Society of Fellows for the wonderful working conditions and  support as a postdoc and Junior Fellow during 2023-2024. He also wants to thank the Max Planck Institute for Mathematics for the support during the correction stage of this paper.  The first author was supported by the project "Group schemes, root systems, and related representations" founded by the European Union - NextGenerationEU through Romania's National Recovery and Resilience Plan (PNRR) call no. PNRR-III- C9-2023-I8, Project CF159/31.07.2023, and coordinated by the Ministry of Research, 
Innovation and Digitalization (MCID) of Romania. He would like to thank the Institute of Mathematics "Simion Stoilow" of the Romanian Academy for the wonderful working conditions.

\section{Preliminaries}
\label{SectionPreliminaries}

In this section we introduce the main objects and techniques used in the paper. In \cref{ss:GeneralNotations} we give some general definitions and notations.  In \cref{SubsecFFCurve} we recall the definition of families of Fargues-Fontaine curves following \cite[\S 11.2]{ScholzeWeinspadicgeometry} and \cite[\S II.1]{FarguesScholze}. Some period sheaves associated to affinoid charts of the curves will be explicitly introduced for reference in later sections, and a careful study is made in  \Cref{bplus}, one of the few new and rather delicate results in this chapter.  In \cref{SubsecLocAnRep}  we briefly introduce the theory of solid locally analytic representations of  \cite{RRLocallyAnalytic,RJRCSolidLocAn2}, in particular we extend the criterion of \cite[Proposition 3.4.1]{RJRCSolidLocAn2} for a solid representation of a $p$-adic Lie group to be locally analytic. We end this preliminary chapter with a brief review of geometric Sen theory, as developed in \cite{RCGeoSen}, which will play a crucial role in this paper.

\subsection{General notations}\label{ss:GeneralNotations}

   Fix once and for all a prime number $p$.
   Let $\Perfd$ be the category of perfectoid spaces (over $\bb{Z}_p$)
    and let 
   $\Perf$ be the full subcategory of perfectoid spaces of characteristic $p$. 
   We let $(-)^{\flat}: \Perfd\to \Perf$ be the tilting functor, which induces equivalences
   $\Perfd/X\simeq \Perf/X^{\flat}$ and $X_{\rm et}\simeq X^{\flat}_{\rm et}$ for any $X\in \Perfd$ (\cite[Corollary 3.20 and Theorem 6.3]{scholze2022etale}). 
   If $S\in \Perf$, an untilt of $S$ is a pair $(S^{\sharp}, \iota)$ with $S^{\sharp}\in \Perfd$ and $\iota: (S^{\sharp})^{\flat}\simeq S$. We will simply denote an untilt by $S^{\sharp}$ from now on.   
 
     We endow $\Perfd$ with the $v$-topology, for which a family of morphisms $\{f_i: X_i\to X\}_{i\in I}$ is a covering 
if for any quasi-compact open $U\subset X$ there is a finite subset $J\subset I$ and quasi-compact open subsets 
$U_j\subset X_j$ ($j\in J$) such that $U=\cup_{j\in J} f_j(U_j)$. This topology is subcanonical and the functors 
$\mathcal{O}: X\mapsto \mathcal{O}_X(X)$, $\mathcal{O}^+: X\mapsto \mathcal{O}^+(X)$ are sheaves, such that 
$H^i(X, \mathcal{O})=0$ and $H^i(X, \mathcal{O}^+)$ is almost zero for $i>0$ and $X$ affinoid perfectoid (see \cite[Corollary 8.6, Theorem 8.7 and  Proposition 8.8]{scholze2022etale}). 

 Recall that a small $v$-sheaf is a $v$-sheaf $X$ on $\Perf$ which is a quotient of (the sheaf represented by) a perfectoid space and that a small $v$-stack is a $v$-stack (of grupoids) $X$ on $\Perf$ for which there is a surjection of $v$-stacks $S\to X$ with 
 $S\in \Perfd$ such that $S\times_X S$ is a small $v$-sheaf. Let ${\rm Stack}_v$ be the category of small $v$-stacks.
 There is a functor $$|-|: {\rm Stack}_v\to {\rm Top}$$ to the category of topological spaces, sending surjections of small $v$-stacks to surjective quotient maps of topological spaces (and conversely, if $f: X\to Y$ is a \emph{qcqs} map in 
 ${\rm Stack}_v$ for which $|f|$ is surjective, then $f$ is surjective). 
  See \cite[Section 12]{scholze2022etale} for these results. A $v$-sheaf $X$ is called \emph{spatial} if $X$ is qcqs (and then $X$ is a small $v$-sheaf) and the $|U|$ form a basis for the topology of $|X|$ when $U$ runs through the qcqs open subfunctors of $X$. More generally, a small $v$-sheaf $X$ is called locally spatial if it has an open cover by spatial $v$-sheaves. 

 The presheaf $\bb{Z}_p^{\lozenge}={\rm Spd}(\bb{Z}_p)$ sending 
   $S\in \Perf$ to the set of isomorphism classes of untilts of $S$ is a $v$-sheaf \cite[Lemma 15.1]{scholze2022etale}. For any analytic adic space $X$ over ${\rm Spa}(\bb{Z}_p)$ we let $X^{\lozenge}$ be the presheaf sending $S\in \Perfd$ to the set of isomorphism classes of untilts of $S$ over $X$ (i.e. pairs $(S^{\sharp}, f)$ where $S^{\sharp}$ is an untilt of $S$ and $f: S^{\sharp}\to X$ is a map of analytic adic spaces over $\bb{Z}_p$). By \cite[Lemma 15.6]{scholze2022etale} this defines a $v$-sheaf, more precisely a locally spatial diamond, and there are natural identifications $|X^{\lozenge}|=|X|$, $X_{\rm et}=X^{\lozenge}_{\rm et}$, $X_{\rm fet}=X^{\lozenge}_{\rm fet}$. If $X={\rm Spa}(A, A^{\circ})$ is an affinoid, we simply write $A^{\lozenge}$ instead of $X^{\lozenge}$.

\subsection{The Fargues-Fontaine curve and sheaves of periods}
\label{SubsecFFCurve}
We refer to Section II.1 in \cite{FarguesScholze} for the proofs of the facts recalled in this section.

\subsubsection{Various relative curves} 

   For any $S\in \Perf$ there are natural sous-perfectoid analytic adic spaces 
   $\n{Y}_S, Y_S, X_S$ endowed with natural identifications of $v$-sheaves on $\Perf$
    $$\n{Y}_S^{\lozenge}\simeq \bb{Z}_p^{\lozenge}\times S,\,\, Y_S^{\lozenge}\simeq \bb{Q}_p^{\lozenge}\times S, \,\, X_S^{\lozenge}\simeq (\bb{Q}_p^{\lozenge}\times S)/(
{\rm id}\times \varphi_S^{\bb{Z}}),$$
where $\varphi_S$ is the absolute Frobenius automorphism of $S$.
The space $X_S$ is called the \emph{relative Fargues-Fontaine curve over $S$}
  and by construction we have $$Y_S=\n{Y}_S\times_{{\rm Spa}(\bb{Z}_p)} {\rm Spa}(\bb{Q}_p), \,\, X_S=Y_S/\varphi_S^{\bb{Z}},$$
the Frobenius automorphism $\varphi_S$ on $Y_S$ being totally discontinuous and extending to an automorphism $\varphi_S$ of $\n{Y}_S$. The construction of $\n{Y}_S$ is done in two stages: for 
$S$ affinoid we will recall it below, and it is compatible with 
open immersions of affinoid perfectoid spaces $S\to S'$, allowing to glue the various $\n{Y}_U$ for 
affinoid perfectoid subspaces $U$ of $S$.

Suppose that 
 $S=\Spa(R,R^+)\in \Perf$ is affinoid and let $\varpi\in R$ be a pseudo-uniformizer. 
Endow the ring of $p$-typical Witt vectors $W(R^+)$ with the $(p, [\varpi])$-adic topology (where $[-]: R^+\to W(R^+)$ is the Teichm\"uller lift) and set 
$$\Spa(W(R^+)):={\rm Spa}(W(R^+), W(R^+)).$$
Then  
$$\n{Y}_{S}=\Spa(W(R^+))\setminus \{|[\varpi]|=0\}$$
is an open subspace of $\Spa(W(R^+))$, independent of 
$\varpi$, the automorphism $\varphi_S$ being 
induced by the natural Frobenius lift on $W(R^+)$.

\subsubsection{Closed Cartier divisors}

   If $X$ is a uniform analytic adic space we denote by ${\rm CCD}(X)$ the set of closed Cartier divisors of 
   $X$ (see \cite[Definition 5.3.7]{ScholzeWeinspadicgeometry} for this notion). Let 
  $${\rm Div}^1=\bb{Q}_p^{\lozenge}/\varphi^{\bb{Z}}.$$
  By \cite[Propositions II.1.4 and II.1.18]{FarguesScholze} there are natural injections 
  $$\bb{Z}_p^{\lozenge}(S)\to {\rm CCD}(\n{Y}_S),\,\, \bb{Q}_p^{\lozenge}\to  {\rm CCD}(Y_S), \,\, 
  {\rm Div}^1(S)\to {\rm CCD}(X_S)$$
  for any $S\in \Perf$. More precisely a map $S\to \bb{Z}_p^{\lozenge}$ (resp. $S\to \bb{Q}_p^{\lozenge}$) corresponds to an untilt $S^{\sharp}$ of $S$ (resp. an untilt over $\bb{Q}_p$) and the corresponding closed Cartier divisor is $S^{\sharp}$: locally on $S={\rm Spa}(R,R^+)$ we have 
  $S^{\sharp}={\rm Spa}(R^{\sharp}, R^{\sharp, +})$ with 
  $R^{\sharp,+}=W(R^+)/(\xi)$ for some nonzero divisor 
  $\xi$, and the corresponding invertible ideal sheaf is $\xi \mathcal{O}_{\n{Y}_S}$ (resp. $\xi \mathcal{O}_{Y_S}$). Moreover, for a map $S\to \bb{Q}_p^{\lozenge}$ the composition $S^{\sharp}\to Y_S\to X_S$ is still an element of ${\rm CDD}(X_S)$, depending only on the composite $S\to \bb{Q}_p^{\lozenge}\to {\rm Div}^1$.

\subsubsection{Period sheaves}

 We introduce now some of the main players in this paper.
We make the convention that all compact intervals $I=[s,r]\subset [0,\infty)$ 
 considered below have rational endpoints (we allow $s=r$). 
 
 It will be convenient to work with 
 test objects $S$ endowed with a pseudo-uniformizer, since the spaces below will depend on the pseudo-uniformizer (or more precisely, on the norm defined by the pseudo-uniformizer). 
 Let $\bb{F}_p((\varpi^{1/p^{\infty}}))$ be the perfectoid field parametrizing pseudo-uniformizers in $\Perf$ and let $\Perf_{\varpi}$ be the slice category $\Perf_{/\Spa \bb{F}_p((\varpi^{1/p^{\infty}}))}$, equivalently, $\Perf_{\varpi}$ is the category of perfectoid spaces in characteristic $p$ equipped with a pseudo-uniformizer $\varpi$, and with maps preserving the pseudo-uniformizer. 
 
 For a compact interval $I=[s,r]\subset [0,\infty)$ and $S\in \Perf_{\varpi}$ consider the open subspace
  \[
\n{Y}_{S,I}=\{|p|^r\leq |[\varpi]|\leq |p|^s\}\subset \n{Y}_{S}.
\]
 If $S={\rm Spa}(R,R^+)$ then $\n{Y}_{S,I}$ is an affinoid space, more precisely 
  a rational subspace of $\Spa(W(R^+))$.
 
  For an arbitrary interval $I\subset [0,\infty)$ 
let $\n{Y}_{S,I}\subset \n{Y}_{S}$
be the open subspace 
 \[
 \n{Y}_{S,I}=\bigcup_{J\subset I} \n{Y}_{S,J}
 \]
 where $J$ runs over all the compact subintervals of $I$, in particular 
 $$\n{Y}_{S}=\n{Y}_{[0,\infty), S},\,\, Y_S=\n{Y}_{(0,\infty), S}.$$ 
  
  \begin{remark} There is a continuous radious map
\[
\rad_{\varpi}: |\n{Y}_{S}|^{\rm Berk}\to [0,\infty), \,\, \rad(x)=\frac{\log |[\varpi](x)|}{\log |p(x)|},
\]
such that $\rad_{\varpi}\circ \varphi_{S}=p\rad$, where $|\n{Y}_{S}|^{\rm Berk}$ is the maximal Hausdorff quotient of 
$|\n{Y}_S|$. We have $\n{Y}_{S, I}\subset {\rm rad}_{\varpi}^{-1}(I)$, with equality for an open interval $I$ (but not for a compact interval!).
\end{remark}

\begin{definition}\label{def:PeriodSheaves}
Let $I\subset [0,\infty)$ be a compact interval. We define presheaves 
$\bb{A}_{\inf}$, $\bb{B}_I$ and $\bb{B}_I^+$ on affinoid spaces 
$S\in \Perf_{\varpi}$ 
by
$$\bb{A}_{\inf}(S)=W(R^+),\,\,  \bb{B}_{I}^{+}(S)=\s{O}^{+}(\n{Y}_{S,I}),\,\,\bb{B}_{I}(S)=\s{O}(\n{Y}_{S,I}).$$
\end{definition}

\begin{proposition}\label{PropositionSheavesvSite}
 These presheaves are $v$-sheaves on affinoid perfectoid spaces in $ \Perf_{\varpi}$, thus extend to 
 $v$-sheaves on $\Perf_{\varpi}$.
\end{proposition}

\begin{proof}
That $\bb{A}_{\inf}$ is a sheaf follows immediately from the fact that $\s{O}^+$ is a $v$-sheaf on 
$\Perf$. It remains to see that if $T\to S$ is a $v$-cover of affinoids in $\Perf_{\varpi}$, then 
the natural sequences 
$$0\to \s{O}^{(+)}(\n{Y}_{S,I})\to \s{O}^{(+)}(\n{Y}_{T,I})\to \s{O}^{(+)}(\n{Y}_{T\times_S T,I})$$
are exact. The case of $\s{O}$ is proved in the proof of \cite[Proposition II.2.1]{FarguesScholze}, and the case of 
$\s{O}^+$ follows from this and the fact that $\n{Y}_{T,I}\to \n{Y}_{S,I}$ is surjective (which follows as well from loc.cit.). 
\end{proof}

\subsubsection{Reductions of integral period sheaves}
 
 For a positive rational number $r=d/h$ (with $d,h$ relatively prime positive integers), we will say that 
    $\varpi^{1/r}$ exists in $S={\rm Spa}(R,R^+)$ if $\varpi$ has a $d$th root in $R$, and then we let 
    $\varpi^{1/r}$ an $h$-power of this root. Note that $\varpi^{1/r}$ exists \'etale locally on $S$. If 
    $A$ is a $W(R^+)$-algebra, let $A\langle T\rangle$ be the $[\varpi]$-adic completion of the polynomial ring $A[T]$.
    
\begin{theorem}\label{bplus}
   Let $I=[s,r]\subset (0,\infty)$ and let $S=\Spa(R,R^+)\in \Perf_{\varpi}$ be affinoid perfectoid such that $\varpi^{1/r}$ exists in $S$. Then 
   
   a) There is a natural isomorphism of $W(R^+)$-algebras 
   $$\bb{B}_{[0,r]}^+(S)\simeq W(R^+)\langle T\rangle/( [\varpi^{1/r}]T-p),\quad 
\bb{B}_{[0,r]}^+(S)/[\varpi^{1/r}]\simeq (R^+/\varpi^{1/r}) [T].$$

b) If $\varpi^{1/s}$ also exists in $S$, then there is a natural exact sequence 
$$0\to \bb{B}_{[0,r]}^+(S)\langle U\rangle \xrightarrow{ p[\varpi^{-1/r}]U-[\varpi^{\frac{1}{s}-\frac{1}{r}}]} \bb{B}_{[0,r]}^+(S)\langle U\rangle\to 
\bb{B}_I^+(S)\to 0$$
and an isomorphism of $W(R^+)$-algebras 
$$\bb{B}_I^+(S)/[\varpi^{1/r}]\simeq (R^+/\varpi^{1/r})[T,U]/(TU-\varpi^{1/s-1/r}).$$
In particular, we have an isomorphism of $W(R^+)$-modules
$$\bb{B}_I^+(S)/[\varpi^{1/r}]\simeq \bigoplus_{n\geq 1} (R^+/\varpi^{1/r}).$$ 
\end{theorem}

\begin{proof} Let $$\alpha=\varpi^{1/r},\,\, \beta=\varpi^{1/s},\,\, t=\frac{p}{[\alpha]}.$$ 
 For $?\in \{+, 0\}$ let 
$$A^{?}=W(R^{?}), \,\, B^{?}=A^{?}[t]\subset W(R).$$

  By definition of rational localizations, 
  $\bb{B}^+_{[0,r]}(S)$ is the integral closure of the  $[\varpi]$-adic completion 
    of $B^+=A^+[t]$ in $B^+[1/[\varpi]]$ and 
  $\bb{B}^+_{[s,r]}(S)$ is the integral closure of the $p$-adic (or equivalently the $[\varpi]$-adic) completion of $B^+[[\beta]/p]=A^+[t, [\beta]/p]\subset W(R)[1/p]$ in $B^+[1/p]=B^+[\beta]/p][1/[\varpi]]$. We claim that $\bb{B}^+_{[0,r]}(S)$ is already the $[\varpi]$-completion of $B^+$ (resp. that  $\bb{B}^+_{[s,r]}(S)$ is the $[\varpi]$-adic completion of $B^+[\beta]/p]$), by \cite[Lemma 5.1.2.(i)]{BhattPerfectoid} it suffices to show that $B^+$ (resp. $B^+[[\beta]/p]$) is already integrally closed.   First, some preliminaries:

\begin{lemma}\label{intclosed4}
Let $C\to D$ be a map of rings.

a) If $P\in D[X]$ is integral over $C[X]$, then each coefficient of $P$ is integral over $C$.

b) If $x\in \{0,1\}$ and $P\in D[T,U]/(TU-x)$ is integral over $C[T,U]/(TU-x)$, then 
each coefficient of $P$ is integral over $C$.

\end{lemma}

\begin{proof} a) Let $Q_{n-1},...,Q_0\in C[X]$ such that $P^n+Q_{n-1}P^{n-1}+...+Q_0=0$. Take $N$ large enough and set 
$\tilde{P}=X^N+P$. Then $\tilde{P}^n+R_{n-1}\tilde{P}^{n-1}+...+R_0=0$, with $R_i\in C[X]$ and $R_0$ monic. We can thus find 
$\tilde{S}\in D[X]$ monic such that $\tilde{P}\tilde{S}=R_0$. There is a finite injective map $D\to D'$ such that in 
$D'[X]$ we have $\tilde{P}=\prod_{i=1}^n (X-a_i)$ and $\tilde{S}=\prod_{i=1}^m (X-b_i)$ for some $a_i,b_j\in D'$. Then 
$R_0(a_i)=R_0(b_j)=0$, thus $a_i, b_j$ are integral over $C$ and so the coefficients of $\tilde{P}$ are also integral over 
$C$, which yields the result.

b) This follows immediately from a), since integral closure commutes with localization and  
for any ring $R$ one has a natural embedding $R[T,U]/(TU)\subset R[T]\times R[U]$.
\end{proof}

   The pseudo-uniformizer $\varpi$ gives rise to a sub-multiplicative norm 
   $|-|_{\varpi}$ on $R$ such that $|\varpi|_{\varpi}\cdot |\varpi^{-1}|_{\varpi}=1$. Let 
$|-|_{\rm sp}$ be the corresponding spectral norm on $R$ (so $|f|_{\rm sp}=\lim_{n\to\infty} |f^n|_{\varpi}^{1/n}$), a sub-multiplicative and power-multiplicative norm on 
 $R$, with unit ball $R^0$ and defining the Banach topology of $R$ (since $R$ is uniform).

\begin{lemma} \label{intclosedzero}
 Let 
 $$S=\{\sum_{n\geq 0} [a_n]p^n\in W(R)|\,\, \lim_{n\to\infty} \alpha^n a_n=0\}.$$
 Then $S[1/p]$ is a subring of $W(R)[1/p]$ containing $B^0[1/p]$, and for 
 $\gamma\in \{\alpha, \beta\}$
 setting 
 $$|\sum_{n>-\infty} [a_n]p^n|_{\gamma}=\max_{n}|a_n\gamma^n|_{\rm sp}$$
 yields a sub-multiplicative and power-multiplicative norm on $S[1/p]$. 
\end{lemma}

\begin{proof}
 This follows from Proposition $5.1.2$ in \cite{KL}, observing that $|f \gamma|_{\rm sp}=|f|_{\rm sp}\cdot |\gamma|_{\rm sp}$ 
 for $f\in R$ (since $|\varpi|_{\varpi}\cdot |\varpi^{-1}|_{\varpi}=1$), thus the arguments in loc.cit. for $S$ apply to $S[1/p]$ as well.
\end{proof}

 \begin{lemma}\label{intclosed2}
  Let $\gamma\in R^+$ be divisible by $\alpha$ in $R^+$ and let $b\in B^+$. We have $[\gamma]b\in pB^+$ if and only if 
  $b\in tB^+$.
 \end{lemma}

\begin{proof}
 One direction is clear since $[\alpha]t=p$. Suppose that 
 $[\gamma]b\in pB^+$. Since $b\in B^+=A^+[t]$ and $[\gamma] (tB^+)\subset pB^+$, we may assume that 
 $b\in A^+$. But then $[\gamma] b\in pB^+\subset pW(R)$, so $b\in pA^+\subset tB^+$.
\end{proof}

\begin{lemma}\label{intclosed3}
a) The map $T\mapsto t$ yields an isomorphism of $A^+$-algebras
$$A^+[T]/([\alpha]T-p)\simeq B^+$$
and $[\alpha]T-p$ is a non-zero divisor in $A^+[T]$. In particular 
$$B^+/([\alpha])\simeq (R^+/\alpha)[T].$$

b) The maps $U\mapsto [\beta]/p$ and $T\mapsto t, U\mapsto [\beta]/p$ induce 
 isomorphisms 
 $$B^+[U]/(tU-[\beta/\alpha])\simeq B^+[[\beta]/p]\simeq A^+[T, U]/([\alpha]T-p, TU-[\beta/\alpha])$$
 and $tU-[\beta/\alpha]$ is a non-zero divisor in $B^+[U]$.  
 In particular
 $$B^+[[\beta]/p]/([\alpha])\simeq (R^+/\alpha)[T,U]/(TU-\beta/\alpha).$$

\end{lemma}

\begin{proof}
a) Since $\alpha$ is a unit in $R$ (in particular a non-zero divisor in $R^+$) and $A^{+}/(p)\simeq R^{+}$, the sequence $p, [\alpha]$ is regular in $A^+$, which yields the result.

b) It suffices to prove the isomorphism $B^+[U]/(tU-[\beta/\alpha])\simeq B^+[[\beta]/p]$, the remaining parts being consequences of part a).  
The map is obviously well-defined and surjective, and 
 $f:=tU-[\beta/\alpha]$ is a non-zero divisor in $B^+[U]$ since $t$ is a non-zero divisor in 
 $B^+=A^+[t]\subset W(R)$. 
   Let $P=b_0+b_1U+...+b_dU^d\in B^+[U]$ of degree $d>0$ such that $P([\beta]/p)=0$. Then $b_0p^d+b_1p^{d-1}[\beta]+...+b_d[\beta]^d=0$, thus $b_d[\beta^d]\in pB^+$ and by  \Cref{intclosed2} we have $b_d\in tB^+$. Writing 
   $b_d=tc_d$ with $c_d\in B^+$, $P-c_{d}fU^{d-1}$ still vanishes at $[\beta]/p$ and has degree smaller than $d$, thus by induction we obtain $P\in fB^+[U]$, finishing the proof.   
\end{proof}

\begin{lemma}\label{intclosed1}
a) $B^+=A^+[t]$ is integrally closed in $B^+[1/[\varpi]]$. 

b) $B^+[[\beta]/p]$ is integrally closed in $B^+[1/p]$.
\end{lemma}

\begin{proof} 
a) Recall the ring $S$ from  
 \Cref{intclosedzero}. We claim that $B^0=A^0[t]$ is the unit ball of the norm 
 $|-|_{\alpha}$ on $B^0[1/[\varpi]]$. 
  Clearly $B^0$ is contained in the unit ball, so take 
 $f=\sum_{n\geq 0} [a_n] p^n\in B^0[1/[\varpi]]$ with $|f|_{\alpha}\leq 1$. 
 Then all $\alpha^n a_n$ belong to $R^0$, and they tend to $0$ as $n\to \infty$. Writing 
 $$f=\sum_{n\geq 0} [a_n \alpha^n] t^n,$$
 we see that
 $f$ is in the $[\varpi]$-adic completion of $B^0=A^0[t]$. But $B^0$ is $[\varpi]$-torsion free, thus 
 the intersection of its $[\varpi]$-adic completion with $B^0[1/[\varpi]]$ is $B^0$, proving the claim.
 
 Now let $f\in B^+[1/[\varpi]]$ be integral over $B^+$, in particular $f\in B^0[1/[\varpi]]$ is integral over $B^0$. Since 
 $|-|_{\alpha}$ is sub-multiplicative and power-multiplicative, with unit ball $B^0$, we must have $f\in B^0$ and 
 $f$ is integral over $B^+$. Since 
 $[\alpha] B^0\subset B^+$ (because $[\alpha] A^0\subset A^+)$, in order to show that $f\in B^+$ it suffices to check that the 
 image of $f$ in $B^0/([\alpha])$ belongs to the image of the natural map $B^+/([\alpha])\to B^0/([\alpha])$. 
  By  \Cref{intclosed3} the image of 
 $f$ in $B^0/([\alpha])$ corresponds to a polynomial $P\in (R^0/\alpha)[T]$, and we want to show that 
 $f$ belongs to the image of $(R^+/\alpha)[T]\to (R^0/\alpha)[T]$. 
 Now $P$ is integral over 
 $B^+/([\alpha])\simeq (R^+/\alpha)[T]$, so by  \Cref{intclosed4} the coefficients of $P$ are integral over 
 $R^+/\alpha$. Since $R^+$ is integrally closed in $R$ and contains $\alpha R^0$, any element of 
 $R^0/\alpha$ that is integral over $R^+/\alpha$ must be in the image of $R^+/\alpha\to R^0/\alpha$, hence we are done. 

b) As in the proof of part a), 
$B^0[[\beta]/p]$ is the unit ball for the power-multiplicative and sub-multiplicative norm $\max(|-|_{\alpha}, |-|_{\beta})$ on $B^0[1/p]=B^0[[\beta]/p][1/[\varpi]]$, hence any element of $B^+[1/p]$ that is integral over $B^+[[\beta]/p]$ must belong to 
$B^0[[\beta]/p]$. By 
 \Cref{intclosed3}
 $$B^+[[\beta]/p]/([\alpha])\simeq (R^+/\alpha)[T,U]/(TU-\beta/\alpha)$$
 and $\beta/\alpha=\varpi^{\frac{1}{s}-\frac{1}{r}}$ is either a unit (if $s=r$) or is topologically nilpotent (is $s<r$). In the first case, the same argument as in a) allows to conclude by working modulo $([\alpha])$. In the second case, it suffices to work modulo $([\gamma])$, where 
 $\gamma$ is a positive power of $\varpi$ dividing both $\alpha$ and $\beta/\alpha$, so that $B^+[[\beta]/p]/([\gamma])\simeq (R^+/\gamma)[T,U]/(TU)$ and we can still apply  \Cref{intclosed4}.
   \end{proof}

  At this step we know that 
  $$\bb{B}^+_{[0,r]}(S)=\widehat{B^+},\,\, \bb{B}^+_{[s,r]}(S)=\widehat{B^+[[\beta]/p]},$$
  the completion being $[\varpi]$-adic, or equivalently $[\alpha]$-adic. 
  
  Let $$I=([\alpha]T-p) A^+[T],\,\, J=(tU-[\beta/\alpha])B^+[U].$$
 By  \Cref{intclosed3} we have   
  $$B^+\simeq A^+[T]/I, \,\, B^+[[\beta]/p]\simeq B^+[U]/J.$$
  
  \begin{lemma}\label{intclosed6}
   We have $I\cap [\alpha]A^+[T]=[\alpha] I$ and $J\cap [\alpha] B^+[U]=[\alpha] J$.
  \end{lemma}
  
  \begin{proof}
  This follows from the isomorphisms preceding the lemma and the fact that $B^+$ and 
  $B^+[[\beta]/p]$ are $[\alpha]$-torsion free.
  \end{proof}
  
   By induction it follows that $I\cap [\alpha]^n A^+[T]=[\alpha]^n I$ for all $n$, so letting $f=[\alpha]T-p$ we obtain an exact sequence 
  $$0\to A^+[T]/[\alpha]^n\xrightarrow{f} A^+[T]/[\alpha]^n\to B^+/[\alpha]^n\to 0$$
  for all $n$. Passing to the limit yields an exact sequence 
  $$0\to A^+\langle T\rangle \xrightarrow{f} A^+\langle T\rangle \to \widehat{A^+[t]}=\bb{B}^+_{[0,r]}(S)\to 0$$
  where 
$A^+\langle T\rangle$ is the $[\varpi]$-adic completion of $A^+[T]$. 
 Thus
  $$\bb{B}^+_{[0,r]}(S)/[\alpha]\simeq A^+\langle T\rangle/([\alpha], [\alpha]T-p)\simeq (R^+/\varpi^{1/r})[T].$$
  This proves part a) of  \Cref{bplus}.

  Similar arguments yield an exact sequence 
 $$0\to B^+\langle U\rangle \xrightarrow{tU-[\beta/\alpha]} B^+\langle U\rangle\to \bb{B}^+_{[s,r]}(S)\to 0,$$
 where $B^+\langle U\rangle$ is the $[\varpi]$-adic completion of $B^+[U]$, which is also $\bb{B}^+_{[0,r]}(S)\langle U\rangle$.
 In particular (recall that $[\alpha]\mid p$ in $B^+$)
 $$\bb{B}^+_{[s,r]}(S)/([\alpha])\simeq B^+[U]/([\alpha], tU-[\beta/\alpha])\simeq $$
 $$A^+[T, U]/([\alpha], TU-[\beta/\alpha], [\alpha]T-p)=
 (R^+/\varpi^{1/r})[T, U]/(TU-\varpi^{1/s-1/r}).$$
  Finally, observe that if $C$ is any ring that is $\varpi$-adically complete, then 
  $C[T,U]/(TU-\varpi^{1/s-1/r})$ is a $\varpi$-completely free $C$-module of countable rank (for example, by lifting a basis of its $\varpi^{\frac{1}{s}-\frac{1}{r}}$-reduction and using topological Nakayama's lemma). 
\end{proof}

\begin{corollary}\label{bpluscor1}
 Consider the category $\mathrm{Perf}_{\varpi}$ of perfectoid spaces in characteristic $p$ with pseudo-uniformizer $\varpi$. Let  $I=[s,r]\subset (0,\infty)$ be an interval such that $\varpi^{1/s}$ and $\varpi^{1/r}$ exist in $\mathbb{F}_p((\varpi^{1/p^{\infty}}))$. Then  there are natural isomorphisms of $v$-sheaves of rings over $\mathrm{Perf}_{\varpi}$
$$
\bb{B}_{[0,r]}^+/([\varpi^{1/r}])\simeq (\widehat{\s{O}}^+/(\varpi^{1/r})) [T]$$
and 
$$\bb{B}_I^+/([\varpi^{1/r}])\simeq (\widehat{\s{O}}^+/(\varpi^{1/r})) [T, U]/(TU-\varpi^{\frac{1}{s}-\frac{1}{r}}).$$
In particular, there is a natural isomorphism of $\mathbb{A}_{\inf}$-modules 
$$\bb{B}_I^+/([\varpi^{1/r}])\simeq \bigoplus_{n\geq 1} (\widehat{\s{O}}^+/(\varpi^{1/r})).$$

\end{corollary}
      \begin{proof}
      The follows directly from \cref{PropositionSheavesvSite} and  \cref{bplus}.
      \end{proof}

    In the following corollary we consider almost mathematics with respect to the ideal of $\bb{A}_{\inf}(\bb{F}_p((\varpi^{1/p^{\infty}})))$ generated by $([\varpi]^{1/p^n})_{n\in \bb{N}}$. 

\begin{corollary}\label{LemmaVanishingCohoBI}
Let $S=\Spa(R,R^+)\in \Perf_{\varpi}$ be affinoid perfectoid.

\begin{enumerate}

\item  We have an almost equivalence of derived $(p,[\varpi])$-adically complete complexes
\[
\bb{A}_{\inf}(R^+) =^a R\Gamma_{v}(S,\bb{A}_{\inf}). 
\]

\item  For $I\subset[0,\infty)$ a compact interval we have a natural 
almost equivalence of derived  $[\varpi]$-adically complete complexes
\[
\bb{B}_I^+(S)=^a R\Gamma_v(S,\bb{B}_I^+).
\]
and a natural equivalence
\[
\bb{B}_I(S)=R\Gamma_{v}(S, \bb{B}_I)
\]

\end{enumerate}
\end{corollary}

\begin{proof}
By \cite[Proposition 8.8]{scholze2022etale} we have a natural almost equivalence $R^+=^a R\Gamma_{v}(S,\s{O}^+)$ and so an almost equivalence modulo any pseudo-uniformizer. Since $\bb{A}_{\inf}/(p,[\varpi])=^a\widehat{\s{O}}^+/\varpi$ as almost $v$-sheaves, by derived Nakayama's lemma we have an almost equivalence of derived  $(p,[\varpi])$-complexes
\[
\bb{A}_{\inf}(R^+)=^a R\Gamma_{v}(S,\bb{A}_{\inf})
\]
proving (1). The second assertion in (2) follows from the first one by inverting pseudo-uniformizers. For the first, it suffices to use again 
the almost acyclicity of $\widehat{\s{O}}^+$ and derived Nakayama's lemma, together with the isomorphism of 
$\widehat{\s{O}}^+$-modules of   \Cref{bpluscor1}.
\end{proof}

\subsection{Review of solid locally analytic representations}
\label{SubsecLocAnRep}

 In this section we recall the main properties of solid locally analytic representations (introduced in \cite{RRLocallyAnalytic,RJRCSolidLocAn2}) that will be needed later on. 
 
 We fix a compact $p$-adic Lie group $G$, with Lie algebra 
 $\mathfrak{g}$, and a finite extension $K$ of $\bb{Q}_p$, with ring of integers $\mathcal{O}_K$. If $X$ is a T1 topological abelian group let $\underline{X}$ be the corresponding condensed abelian group (so $\underline{X}(S)=C^{\rm cont}(S,X)$ for a profinite set $S$).
 
 \subsubsection{Solid algebra}  Let $\Solid$ be the abelian closed symmetric monoidal category of solid abelian groups (see \cite[Definition 5.1, Theorems 5.8 and 6.2]{ClausenScholzeCondensed2019}). For a ring object $R$ of $\Solid$ let $\Solid_{R}$ be the abelian closed symmetric monoidal category of solid $R$-modules (i.e. $R$-module objects in $\Solid$) and let $\ob{D}(R)$ be the derived $\infty$-category of $\Solid_R$. 
 
 The condensed abelian groups 
$\underline{\mathcal{O}_K}$ and $\underline{K}$ are ring objects in $\Solid$. We will write 
$\Sol_K$ instead of $\Sol_{\underline{K}}$ (and call its objects solid $K$-vector spaces or simply $K_{\sol}$-vector spaces) and $\ob{D}(K_{\sol})$ instead of $\ob{D}(\underline{K})$. The category $\Sol_K$ has the following properties:

\begin{itemize}

\item  Let ${\rm LC}_K$ be the category of complete locally convex $K$-vector spaces and let ${\rm LC}_K^{\rm cp. gen}$ be its full subcategory of compactly generated topological $K$-vector spaces. For any $V\in {\rm LC}_K$ we have $\underline{V}\in \Sol_K$ and the resulting functor 
${\rm LC}_K^{\rm cp. gen}\to \Sol_K$ is fully faithful. For $V\in {\rm LC}_K^{\rm cp. gen}$ (e.g. any metrizable complete locally convex $K$-vector space)
we will simply write $V$ instead of $\underline{V}$ from now on. 

\item  The Smith spaces (for various sets $I$)  $${\rm Sm}_I:=(\prod_{I} \mathcal{O}_K)[1/p]$$ form a collection of compact projective generators of $\Sol_K$ and $\Sol_K$ is a closed symmetric monoidal category for a tensor product $\otimes_{K_{\sol}}$ which commutes with colimits separately in each variable and satisfies $${\rm Sm}_I\otimes_{K_{\sol}} {\rm Sm}_J={\rm Sm}_{I\times J}.$$

\item   The duality functor $$V\mapsto V^{\vee}:=\underline{\rm Hom}_K(V, K)$$ in $\Sol_K$ induces an exact anti-equivalence between 
 Smith spaces and (solid $K$-vector spaces attached to) $K$-Banach spaces, and more generally between LS-spaces\footnote{i.e. countable filtered colimits of Smith spaces, with injective transition maps.} and $K$-Fr\'echet spaces. If $V,W\in \Sol_K$ and one is a $K$-Fr\'echet and the other is an LS-space, then 
 $$\underline{\rm Hom}_K(V, W)=V^{\vee}\otimes_{K_{\sol}} W,$$
 in particular $(-)^{\vee}$ is symmetric monoidal on $K$-Fr\'echet spaces and on LS-spaces. If 
 $V,W$ are classical $K$-Fr\'echet spaces then $V\otimes_{K_{\sol}} W$ is simply the usual completed tensor product of $V$ and $W$.  Finally, 
 $K$-Fr\'echet spaces and LS-spaces\footnote{More generally, any quasi-separated object of $\Sol_K$.} are flat, i.e. if $V$ is such a space then 
 $V\otimes^L_{K_{\sol}} W$ is concentrated in degree $0$ for $W\in \Sol_K$.
 
 \end{itemize}

 \begin{remark} It is known that  for a separable Banach $K$-vector space  $B$, the concentration of $R\underline{\rm Hom}_K(B, K)$ in degree zero implies the negation of the continuum hypothesis, see \cite[Theorem B]{MR4984500}. However, for a $K$-Fr\'echet space $V$ one has  $R\underline{\rm Hom}_K(V^{\vee}, K)=V$ .
 \end{remark}

  \subsubsection{Rings of functions and distributions}
    Let 
      $$C^{\rm sm}(G, K)\subset C^{\rm la}(G, K)\subset C^{\rm cont}(G,K)$$
      be the spaces of $K$-valued smooth, locally analytic and continuous functions on $G$, endowed with the left regular action of $G$ (i.e. $(g.f)(x)=f(g^{-1}x)$). The first two are LS spaces and the third is a $K$-Banach space. Let $?\in \{{\rm sm}, {\rm la}, {\rm cont}\}$ and define 
      $$D^{?}(G, K)=C^{?}(G, K)^{\vee},$$
      thus $D^{\rm sm}(G,K)$ and $D^{\rm la}(G,K)$ are $K$-Fr\'echet spaces and $D^{\rm cont}(G,K)$ is a Smith space. 
       Note that $C^?(G,K)$ is naturally a commutative $K$-algebra in $\Sol_K$, and the multiplication on 
      $G$ induces a co-multiplication\footnote{We use here the isomorphism
      $C^?(G\times G,K)\simeq C^?(G, K)\otimes_{K_{\sol}} C^?(G,K)$.}
       $C^?(G,K)\to C^?(G,K)\otimes_{K_{\sol}} C^?(G,K)$, which by duality gives a (non-commutative in general) $K$-algebra structure on $D^?(G,K)$. 
      
\begin{convention} If $G$ is clear from the context we also write 
      $C^?$ and $D^?$ instead of $C^?(G,K)$ and $D^?(G,K)$. In general, we write     $K_{\sol}[G]=D^{\rm cont}(G,K)$   as  it agrees with the free solid $K$-vector space generated by $G$.     
 \end{convention}   
                
                 There are natural isomorphisms of solid $K$-algebras 
       $$D^{\rm sm}\simeq \varprojlim_{H} K[G/H],\,\, D^{\rm cont}\simeq (\varprojlim_{H} \mathcal{O}_K[G/H])[1/p],$$
      $H$ running through open normal subgroups of $G$.
       
        The story is more complicated for $D^{\rm la}$. Fix a sufficiently small $\bb{Z}_p$-lattice $\mathcal{L}$ in $\mathfrak{g}$. By exponentiation and using the Baker-Campbell-Hausdorff formula, it induces a normal open uniform pro-$p$ subgroup $G^u$ of $G$, which is the set of 
        $\bb{Q}_p$-points of an affinoid group $\bb{G}^u$, whose underlying adic space is an affinoid polydisc of radius $1$. For a positive rational number 
        $h$ let $\bb{G}^{u}_h$ be the  affinoid polydisc of radius $p^{-h}$ in $\bb{G}^u$ and let $\bb{G}_h=G\bb{G}^{u}_h$, a finite disjoint union of  affinoid polydiscs. Define $C^h:=C^{h}(G, K):=\mathcal{O}(\bb{G}_h)\otimes_{\bb{Q}_p} K$, the ring of $K$-rational analytic functions on $\bb{G}_h$, a $K$-Fr\'echet space, and define $D^h:=D^h(G,K):=(C^h)^{\vee}$, an LS-space with an associative solid $K$-algebra structure. The rings 
        $C^h$ and $D^h$ depend on the choice of $\mathcal{L}$, but the ind-system $\{C^h\}_h$ does not, has colimit $C^{\rm la}$, and the pro-system 
        $\{D^h\}_h$ is independent of choices and has limit $D^{\rm la}$. It is useful to also consider the Stein group $\bb{G}_{<h}=\bigcup_{h'<h} \bb{G}_{h'}$, its space of functions $C^{<h}:=C^{<h}(G,K):=\mathcal{O}(\bb{G}_{<h})\otimes_{\bb{Q}_p} K$, and the distribution algebra $D^{<h}:=D^{<h}(G,K):=(C^{<h})^{\vee}$. 
        \begin{remark}
        
        \begin{enumerate}
        
        \item A map $A\to B$ of solid $K$-algebras is called idempotent if the natural map $B\otimes^L_{A, \sol} B\simeq B$
       is an equivalence. In this case, the derived $\infty$-category $\ob{D}(B)$ of solid $B$-modules embeds fully faithfully into 
       $\ob{D}(A)$.       
       The maps $D^{\rm la}\to D^{<h}$ and $D^{\rm cont}\to D^{\rm la}$ are idempotent.  
       
       \item The left (resp. right) regular action of $G$ on $C^?(G,K)$ extends to a structure of $D^?(G,K)$-module. 
        \end{enumerate}
        \end{remark}

   \subsubsection{Derived solid, locally analytic and smooth representations}     The category of \emph{derived solid $G$-representations over $K$} is by definition the derived $\infty$-category $$\Rep^{\ob{sol}}_{K_{\sol}}(G):=\ob{D}(K_{\sol}[G])=\ob{D}(D^{\rm cont}(G,K)).$$

       Let $?\in \{\text{\rm sm, la}\}$. Since $C^?$ is a $D^?$-bimodule, we can define an endo-functor  
     $$(-)^{R?}: \ob{D}(D^?)\to \ob{D}(D^?), \,\, V^{R?}:=R\underline{{\rm Hom}}_{D^?}(K, 
     V\otimes^L_{K_{\sol}} C^?),$$
     where we use the left-regular action of $G$ on $C^?$ to form $V\otimes^L_{K_{\sol}} C^?\in \ob{D}(D^?)$ and the right-regular action of 
     $G$ on $C^?$ to turn $V^{R?}$ into an object of $\ob{D}(D^?)$. The endo-functor $(-)^{R?}$ is colimit-preserving, comes with a natural transformation to the identity endo-functor, commutes with restriction to an open subgroup of $G$ and is idempotent (see \cite[Propositions 3.2.6, 3.2.10 and Lemmas 5.1.6, 5.1.8]{RJRCSolidLocAn2}). Let ${\rm Rep}_{K_{\sol}}^?(G)$ be the full sub-category of $\ob{D}(D^?)$ consisting of those 
     $V\in \ob{D}(D^?)$ for which the natural map $V^?\to V$ is an equivalence. 
        
    \begin{prop}\label{PropStabilityLocAn}
 Let $?\in\{la,sm\}$
    \begin{enumerate}

    \item An object $V$ of $\ob{D}(D^{?})$ belongs to  
    ${\rm Rep}_{K_{\sol}}^{?}(G)$ if and only if each $H^i(V)$ belongs to ${\rm Rep}_{K_{\sol}}^{?}(G)$, 
    and ${\rm Rep}_{K_{\sol}}^{?}(G)$ is the derived $\infty$-category of its heart, for the resulting $t$-structure induced by $\ob{D}(D^{?})$. The heart of ${\rm Rep}_{K_{\sol}}^{?}(G)$ is the Grothendieck abelian category of  comodules for the exact functor $C^{?}\otimes_{K_{\sol}} (-)$ on $\Sol_K$.
    
    \item The functor $(-)^{R?}: \ob{D}(D^{?})\to \ob{D}(D^{?})$ factors through ${\rm Rep}_{K_{\sol}}^{?}(G)$
    and yields a right adjoint to the inclusion ${\rm Rep}_{K_{\sol}}^{?}(G)\to \ob{D}(D^{?})$. It satisfies the projection formula (for 
    $V,W\in  \ob{D}(D^{?})$)
    $$(V\otimes^L_{K_{\sol}} W^{R?})^{R?}\simeq V^{R?}\otimes^L_{K_{\sol}} W^{R?}.$$
    In particular ${\rm Rep}_{K_{\sol}}^{?}(G)$ is stable under colimits and solid derived tensor products over $K$. 
    \end{enumerate}
    \end{prop}     
    
    \begin{proof} All references are to \cite{RJRCSolidLocAn2}. Part (1) for the locally analytic case follows from Propositions 3.3.5,   3.3.6  and Lemma 4.3.4.. For the smooth case it follows from  Proposition 5.1.10, and 5.4.2 (rather its proof). 
     
        Part (2) follows from Corollaries 3.2.14 and 3.3.7. in the locally analytic setting.  In the smooth setting, the right adjoint  claim is Proposition 5.1.11, the projection formula is easier since taking smooth vectors is a $t$-exact functor and this can be directly shown in the abelian category. 
  \end{proof}

     Let $?\in \{\ob{sol}, \rm la, \rm sm\}$. 
    The trivial $G$-representation functor ${\rm triv}: \ob{D}(K_{\sol})\to \Rep^?_{K_{\sol}}(G)$
 is colimit-preserving, with right adjoint\footnote{We are abusing notation here, by writing $K$ instead of $K_{\sol}$.} $$R\Gamma^?(G,-)=R\underline{\rm Hom}_{D^?}(K, -): \Rep^?_{K_{\sol}}(G)\to \ob{D}(K_{\sol})
    $$ 
    is colimit-preserving (since $K$ is a compact $D^?$-module). This right adjoints is also the right derived functor of the functor of $G$-invariants at the level of hearts. When $?=\rm sol$, we will usually write simply $R\Gamma(G,-)$ instead of $R\Gamma^{\rm sol}(G,-)$.

     The natural maps 
 $D^{\rm cont}\to D^{\rm la}\to D^{\rm sm}$
 induce forgetful functors 
 $$\Rep^{\rm sm}_{K_{\sol}}(G)\to \Rep^{\rm la}_{K_{\sol}}(G)\to \Rep^{\rm sol}_{K_{\sol}}(G),$$
 the second one being fully faitfhful since $D^{\rm la}$ is idempotent over $D^{\rm cont}=K_{\sol}[G]$. The first functor is \emph{not} fully faithful, but it has 
 a right adjoint given by 
 $$R\Gamma(\mathfrak{g}, -):=R\underline{{\rm Hom}}_{U(\mathfrak{g})}(K, -): {\rm Rep}^{\rm la}(G)\to {\rm Rep}^{\rm sm}(G),$$
 see \cite[Proposition 6.2.1]{RJRCSolidLocAn2}. Using this, one easily obtains the following: 

  \begin{proposition}\label{PropositionComparisonLocAnCohomology}
   For any $V\in\Rep_{K_{\sol}}(G)$ there are natural isomorphisms
   $$R\Gamma(G,V)\simeq R\Gamma(G,V^{R\rm la})\simeq R\Gamma^{\rm la}(G, V^{R\rm la})\simeq R\Gamma^{\rm sm}(G, 
   R\Gamma(\mathfrak{g}, V)).$$
  \end{proposition}

 The functor ${\rm LC}_K\to \Sol_K, V\mapsto \underline{V}$ induces a functor from the category of continuous representations of $G$ on complete locally convex $K$-vector spaces to the heart of $\Rep_{K_{\sol}}^{\ob{sol}}(G)$, and $H^i(R\Gamma(G,\underline{V}))(*)$ is isomorphic to the continuous cohomology $H^i(G,V)$ computed using continuous $V$-valued cochains (see \cite[Lemma 5.2]{RRLocallyAnalytic}). 
 
 \subsubsection{Detecting local analyticity}

     The main technical tool for proving local analyticity of representations is the following criterion.

\begin{proposition}\label{PropLocAnCriterion}
Let $A$ be a sous-perfectoid Banach $\bb{Q}_p$-algebra with pseudo-uniformizer $\varpi$ and let $G$ be a compact $p$-adic Lie group. Let $N$ be a connective complex of derived $\varpi$-complete $A^{\circ}$-linear solid $G$-representations. Let $a>b>0$ be rational numbers such that $\varpi^{a}$ and $\varpi^b$ are defined in $A$. Suppose that there is a trivial $A^{\circ}/\varpi^a$-linear  representation $M$  of $G$ and a $\varpi^b$-retract  of $(A^{\circ}/\varpi^a,\mathbb{Z})_{\sol}[G]$-modules $N/\varpi^a\xrightarrow{f} M\xrightarrow{g} N/\varpi^a$, i.e. $g\circ f$ is given by multiplication by $\varpi^b$ on $N/\varpi^a$. Then $N[\frac{1}{\varpi}]$ is a locally analytic representation of $G$. 
\end{proposition}

\begin{proof}   Since locally analytic representations  are stable under retracts of $\bb{Q}_p$-vector spaces by \Cref{PropStabilityLocAn}, after replacing $A$ by a finite \'etale extension $A'$ we can assume that given $d\in\N$,  $\varpi$ and $p$ admit a $d$-th power root in $A$. In particular, we can choose $d>1/b$ such that  that $p$ is divisible by $\varpi^{1/d}$ in $A^{',\circ}$. Thus, after replacing $\varpi$ by $\varpi^{1/d}$, and $N$ by $N\otimes_{A^{\circ}}^L  A^{',\circ}$, we can assume that $p$ is divisible by $\varpi^b$ in $A^{\circ}$. Let $|-|_{\varpi}$ be the spectral norm on $A$ (i.e. the power-multiplicative submultiplicative norm) making $|\varpi|_{\varpi}\in (0,1)$ a multiplicative unit (eg. we can fix $|\varpi|_{\varpi}=1/2$), in particular we have that $|p|_{\varpi}\leq |\varpi|^b_{\varpi}$.

Consider the space of $h$-analytic functions of $G$ given by $C^{h}(G,K)=\mathcal{O}(\bb{G}_h)$, it can be equally described as the induction $C^{h}(G,K)=\mathrm{Ind}_{G_h}^{G} C^{{\rm an}}(G_h, K)$ where $C^{{\rm an}}(G_h, K)=\mathcal{O}(\mathbb{G}_h^u)$ is a rigid analytic group whose adic space is an affinoid polydis of radius $p^{-h}$, and such that $G_h=\mathbb{G}_h^u(\mathbb{Q}_p)$.  

 In the following all tensor products are with respect to the induced solid structure from $\Z$ and are considered to be derived.  Consider the $\varpi$-complete $A^{\circ}$-linear representation $$C^{h}(G, A^{\circ}):= \mathcal{O}(\mathbb{G}_h)^{\circ}\otimes_{\Z_p} A^{\circ}.$$  For $P$ a derived $\varpi$-complete $A^{\circ}$-linear representation we set $$P_{\varpi^b}=P/^{\bb{L}}\varpi^b, \quad P^{h}:=R\Gamma(G, P\otimes_{A^{\circ}}  C^{h}(G,A^{\circ})).$$  The representation $A^{\circ, h}$ has a natural commutative algebra structure on the $\infty$-derived category of solid abelian groups, and $P^{h}$ has a natural structure of $A^{h}$-module.

  We claim that there is some $h>0$ such that the  fiber of the map 
 \begin{equation}\label{eqoj9oqjoqwmsda}
 N^{h}\otimes_{A^{h}} A\to N
 \end{equation}
 has $\varpi^b$-torsion cohomology groups. In particular, taking colimits as $h\to \infty$ and inverting $\varpi$, we have an isomorphism of $A$-linear $G$-representations
 \[
 (N[\frac{1}{\varpi}])^{{\rm la}}\otimes_{A^{{\rm la}}} A\xrightarrow{\sim} N[\frac{1}{\varpi}.]
 \]
 Since the action of $G$ on $A$ is trivial we have $A^{{\rm la}}=A$ and so $N[\frac{1}{\varpi}]$ is locally analytic, as desired.

To prove the claim we need the following  lemma, see \cite[Lemma 7.1.7]{BhattHodgeTheory}. 

\begin{lemma}\label{LemmaTorsionBhatt}
Let $A$ be a commutative ring and $f,g\in A$.  Let $M$ be a connective derived $f$-complete complex of $A$-modules. If for all $i\in \Z$ the cohomology groups $H^i(M/^{\bb{L}}fg)$ are annihilated by $g$, then so are the cohomology groups of $M$.
\end{lemma}
\begin{proof}
We can  assume without loss of generality that $A=\Z[X,Y]$, and that $f=X$ and $Y=g$. In particular, $g$ is a non-zero divisor on $A$.  We will use the decalage functor of \cite[Section 6]{IntegralPHodgeTheory}.  Since the cohomology groups of  $M/^{\bb{L}}fg$ are $g$-torsion, we have that $L\eta_{g} (M/^{\bb{L}}fg) = 0$ (see \cite[Lemma 6.4]{IntegralPHodgeTheory}). On the other hand, \cite[Lemma 7.9]{IntegralPHodgeTheory} yields 
\[
0=L\eta_{g} (M/^{\bb{L}}fg)= (L\eta_g M) /^{\bb{L}}f.
\]
By \cite[Lemma 6.19]{IntegralPHodgeTheory} the object $L\eta_g M$ is derived $f$-complete, thus derived Nakayama's lemma implies that $L\eta_g M=0$ proving that the cohomology groups of $M$ are annihilated by $g$ as wanted.  
\end{proof}
 
Thus, to show that the fiber of \eqref{eqoj9oqjoqwmsda} has cohomology groups killed by $\varpi^b$, it suffices to prove the analogue statement after reduction module $\varpi^a$.  By hypothesis $N_{\varpi^a}$  is a $\varpi^b$-retract of a trivial $A^{\circ}_{\varpi^a}$-linear $G$-representation $M$. Hence, the $\varpi^a$-reduction of  \eqref{eqoj9oqjoqwmsda} is a $\varpi^b$-retract of the fiber of  $M^{h}\otimes_{A^{\circ, h}}A^{\circ}\to M$ and therefore it suffices to see that this last map is an equivalence.  Consider 
\[
M^h=R\Gamma(G, M\otimes_{A^{\circ}} C^{h}(G,A^{\circ})),
\]
by taking $h$ large enough,  since $h$-analytic functions are induced from a subgroup $G_h\subset G$, we can assume without loss of generality that $G$ is an uniform pro-$p$-group. In that case, the trivial  $G$-representation $\mathbb{Z}_p$ admits a Lazard resolution by finite free $\Z_{p,\sol}[G]$-modules, and therefore group cohomology for $G$ satisfies the projection formula with respect to trivial representations. We deduce that $$M^h=M\otimes_{A^{\circ}} R\Gamma(G, C^{h}(G,A^{\circ}))= M\otimes_{A^{\circ}} A^{\circ,h},$$ which yields an isomorphism 
\[
M^h\otimes_{A^{\circ, h}} A^{\circ}= (M\otimes_{A^{\circ}} A^{\circ,h})\otimes_{A^{\circ, h}} A^{\circ} = M
\]
proving the claim and finishing the proof of the proposition.  
\end{proof}

\subsection{Geometric Sen theory}

\subsubsection{The geometric Sen operator}

 Let $C$ be a perfectoid field extension of $\qp(\mu_{p^{\infty}})$ and let $X={\rm Spa}(A, A^+)$ be a smooth affinoid rigid analytic variety over $C$ with toric coordinates, i.e. for which there is an \'etale map $$f: X\to {\rm \mathbb{T}}^d_{C}:={\rm Spa}(C\langle T_1^{\pm 1},\ldots, T_d^{\pm 1}\rangle)$$ for some 
  $d\geq 0$, which factors as a finite composition of rational localizations and finite \'etale maps. Let 
  $\mathcal{F}$ be a pro\'etale $\widehat{\s{O}}_X$-module of the form $\mathcal{F}=\mathcal{F}^0[1/p]$ for some $p$-adically complete $\widehat{\s{O}}_X^+$-module $\mathcal{F}^0$, such that $\mathcal{F}^0/p^{\varepsilon}\mathcal{F}^0$ is almost isomorphic to a free (possibly of infinite rank) 
$\s{O}_X^+/p^{\varepsilon}$-module, for some $\varepsilon>0$. 
 
 Let 
  $${\rm \mathbb{T}}^d_{\infty, C}={\rm Spa}(C\langle T_1^{\pm 1/p^{\infty}},\ldots, T_d^{\pm 1/p^{\infty}}\rangle)$$ be the natural perfectoid 
  $\Gamma:=\zp(1)^d$-torsor of ${\rm \mathbb{T}}^d_C$, and let $X_{\infty}\to {\rm \mathbb{T}}^d_{\infty, C}$ be the pullback of $X$. By the almost purity theorem 
  $X_{\infty}={\rm Spa}(A_{\infty}, A_{\infty}^+)$ is affinoid perfectoid. Almost purity also implies (see   \cite[Remark 3.2.2]{RCGeoSen}) that $\mathcal{F}$ is entirely described by the $A_{\infty}$-semi-linear representation  $V:=\mathcal{F}(X_{\infty})$ of $\Gamma$ (which is a relatively locally analytic ON $A_{\infty}$-Banach in the sense of   \cite[Definition 2.3.5]{RCGeoSen}).  More precisely, the obvious map induces a $\Gamma$-equivariant isomorphism  of pro\'etale sheaves 
  $$\mathcal{F}|_{X_{\infty}}\simeq V\widehat{\otimes}_{A_{\infty}} \widehat{\s{O}}_X|_{X_{\infty}}.$$ Classical Sen theory (in the form of \cite[Theorem 2.4.4]{RCGeoSen}) applied to $V$ yields the existence of some $h$ such that $V^{h-{\rm la}}$ is an ON Banach $A_{\infty}^{h-{\rm la}}$-module and for which the natural map 
  $$A_{\infty}\widehat{\otimes}_{A_{\infty}^{h-{\rm la}}} V^{h-{\rm la}}\to V$$
  is an isomorphism. 
  
   Using a compatible sequence of $p$-power roots of unity we identify $\Gamma$ with $\mathbb{Z}_p^d=\Gamma_1\times\ldots\times\Gamma_d$ (each $\Gamma_i$ being a copy of $\zp$). There is a unique $A_{\infty}$-linear continuous map $$\tilde{\theta}_V: V\to \bigoplus_{i=1}^d V \frac{dT_i}{T_i} (-1)$$
 with the property that for $v\in V^{h-{\rm la}}$
 $$\tilde{\theta}_V(1\widehat{\otimes} v)=\sum_{i=1}^d (\partial_i v) \frac{dT_i}{T_i} (-1),$$
 where $\partial_i$ is the infinitesimal action of $\Gamma_i$ on $V^{h-{\rm la}}$. The above map induces a 
     a geometric Sen operator 
  $$\theta_{\mathcal{F}}: \mathcal{F}\to \mathcal{F}\otimes_{\s{O}_X}\Omega^1_X(-1)=\bigoplus_{i=1}^d \mathcal{F} \frac{dT_i}{T_i} (-1),$$
which is the same as the one constructed in  \cite[Proposition 3.2.3]{RCGeoSen}, once one unravels its proof.

It turns that this construction is independent of the toric chart $f$ (\cite[Proposition 3.3.1 and Theorem 3.3.2]{RCGeoSen}), and so it globalizes. More precisely, suppose now that $X$ is a general smooth adic space over $C$. A pro\'etale $\widehat{\s{O}}_X$-module 
$\mathcal{F}$ is called \emph{relative locally analytic ON Banach} if \'etale locally on $X$ there is a $p$-adically complete $\widehat{\s{O}}_X^+$-lattice 
$\mathcal{F}^0$ in $\mathcal{F}$ and some $\varepsilon>0$ for which $\mathcal{F}^0/p^{\varepsilon}\mathcal{F}^0$ is almost isomorphic to a free (possibly of infinite rank) 
$\s{O}_X^+/p^{\varepsilon}$-module. The next result is proved in  \cite[Theorem 3.3.2]{RCGeoSen}. 

\begin{prop}\label{gsen}
For any relative locally analytic ON Banach module $\mathcal{F}$ on $X$ there is a canonical 
$\widehat{\s{O}}_X$-linear Higgs field $$\theta_{\mathcal{F}}: \mathcal{F}\to \mathcal{F}\otimes_{\s{O}_X}\Omega^1_X(-1)$$
functorial in $\mathcal{F}$ and endowed with canonical isomorphisms $$R^n\nu_*\mathcal{F}\simeq
\nu_* H^n[0\to \mathcal{F}\to \mathcal{F}\otimes_{\s{O}_X}\Omega^1_X\to...\to \mathcal{F}\to \mathcal{F}\otimes_{\s{O}_X}\Omega^d_X(-d)\to 0],$$
where $d=\dim X$ and $\nu: X_{\rm proet}\to X_{\rm et}$ is the natural projection of sites. We have $\theta_{\mathcal{F}}=0$ if and only if $\mathcal{F}=\widehat{\s{O}}_X\widehat{\otimes}_{\s{O}_X} \mathcal{G}$ for a (locally  in the \'etale site) orthonormalisable Banach $\s{O}_X$-module $\mathcal{G}$. 
  \end{prop}
  
  \subsubsection{The case of pro\'etale torsors} The following construction will occur in various forms throughout the paper:

\begin{definition}\label{DefinitionProetaleSystemsShimura}

 Let $X_{\infty}\to X$ be a $\underline{G}$-torsor over a locally spatial diamond $X$, with $G$ a locally profinite group. Let $V$ be a $p$-adically complete or $\bb{Q}_p$-Banach continuous representation of a closed subgroup $K$ of $G$. The sheaf 
 $\n{F}_{V}$ attached to $V$ is the $v$-sheaf on $X/\underline{K}$ obtained via descent from the $\underline{K}$-equivariant sheaf $\underline{V}$ on $X$
 (recall that $\underline{V}(S)=\Cont(|S|,V)$ for $S\in X_v$ affinoid perfectoid). 
 For an ind-system $V=``\varinjlim_{i}" V_i$ of $p$-adically complete or Banach representations we define $\n{F}_{V}:=\varinjlim_i \n{F}_{V_{i}}$.  Finally, if $\s{A}$ is a period sheaf on the pro\'etale site of $X$, such as $\widehat{\s{O}}$ or $\bb{B}_I$, we let $\n{F}_{V,\s{A}}:= \n{F}_{V}\otimes_{\bb{Z}_p} \s{A}$ denote its solid base change to $\s{A}$.
\end{definition}

Now let $G$ be a compact $p$-adic Lie group, with Lie algebra $\mathfrak{g}$ endowed with the adjoint representation. Let $X_{\infty}\to X$ a pro\'etale $G$-torsor. For any locally analytic 
representation $V$  of $G$ the sheaf $\n{F}_{V,\widehat{\s{O}}}$   is relative locally analytic ON Banach and we let $$\theta_V:=\theta_{\mathcal{F}_{V,\widehat{\s{O}}}}$$ be its geometric Sen operator. By \cite[Theorem 3.3.4]{RCGeoSen} there is a universal 
Sen operator $$\theta_{X_{\infty}}: \mathcal{F}_{\mathfrak{g}^{\vee},\widehat{\s{O}}}\to  \widehat{\s{O}}_X(-1)\otimes_{\s{O}_X}\Omega^1_X,$$
associated to the torsor $X_{\infty}\to X$, encoding all 
$\theta_V$ (for representations $V$ as above), in the sense that $\theta_V$ is the composite 
$$\mathcal{F}_{V,\widehat{\s{O}}}\to \mathcal{F}_{V\otimes_{\qp} \mathfrak{g}^{\vee},\widehat{\s{O}}}\simeq \mathcal{F}_{V,\widehat{\s{O}}} \otimes_{\qp} \n{F}_{\f{g}^{\vee},\widehat{\s{O}}}\to  \mathcal{F}_{V,\widehat{\s{O}}}\otimes_{\s{O}_X} \Omega^1_X(-1),$$
where:

\begin{itemize}
\item The first map is 
 induced by 
the derivation map $d_V: V\to V\otimes_{\bb{Q}_p} \mathfrak{g}^{\vee}$, encoding the infinitesimal action of $\mathfrak{g}$ on $V$.

\item The second,  map is the obvious isomorphism.

\item The third map is the base change of 
$\theta_{X_{\infty}}$.
\end{itemize}

  The following result summarises the main properties of this construction:
  
  \begin{prop} \label{geoSentorsor} Suppose that $X={\rm Spa}(A, A^+)$ a smooth affinoid over $C$, with toric coordinates and that
   $X_{\infty}\to X$ is a pro\'etale $G$-torsor, for a compact $p$-adic Lie group $G$. Let $G$ act by the left regular action on $V:=\mathcal{C}^{\rm la}(G, \qp)$ and let
     $A_{\infty}^{\rm la}$ be the space of locally analytic vectors for the action of $G$ on $A_{\infty}$. Then the following holds.

\begin{enumerate}     
     \item  We have $$A_{\infty}^{\rm la}=H^0_{\proet}(X, \mathcal{F}_{V,\widehat{\s{O}}})$$
     and any element in the image of the map $A_{\infty}^{\rm la}\otimes_A (\Omega^1_A)^{\vee}(1)\to \mathfrak{g}\otimes_{\qp} A_{\infty}^{\rm la}$ (induced by the adjoint of $\theta_{X_{\infty}}$) acts trivially (by derivations) on $A_{\infty}^{\rm la}$. 
     
     \item  If $\theta_{X_{\infty}}$ is surjective, then we have a quasi-isomorphism 
     $${\rm R\Gamma}_{\proet}(X, \mathcal{F}_{V,\widehat{\s{O}}})\simeq A_{\infty}^{\rm la}[0].$$
     
     \end{enumerate}
     
  \end{prop}
    
   \begin{proof}
   This is proved in Corollary 3.2.6 and Proposition 3.2.7 of \cite{RCGeoSen}.
   \end{proof}

\section{Locally analytic vectors of towers of rigid spaces}
\label{SubsecLocAnTower}

     The goal of this section is to prove the following theorem, which is the main innovation of this article. We refer the reader to  \Cref{SubsecLocAnRep}  for 
     recollections on locally analytic vectors
     and the categories of solid $\bb{Q}_p$-linear representations of compact $p$-adic Lie groups.

\begin{theorem}\label{theo:LocAnTwoTowers} Let $G,H$ be compact $p$-adic Lie groups and let 
$X$ be a smooth quasi-compact and separated rigid space over $\bb{C}_p$, endowed with a continuous action of $H$. Let
$\widetilde{X}\to X^{\lozenge}$ be an $H$-equivariant pro\'etale $G$-torsor.

\begin{enumerate}

\item The complex
$R\Gamma_v(\widetilde{X}, \widehat{\s{O}}_X)$ has a natural structure of solid $G\times H$-representation over $\bb{Q}_{p}$. 

\item For any locally analytic representation $V$ of $G$ on a $\bb{Q}_p$-Banach space the solid $H$-representation 
$R\Gamma(G, V\otimes^L_{\bb{Q}_p, \square} R\Gamma_v(\widetilde{X}, \widehat{\s{O}}_X))$ is locally analytic. 

\item The solid $H$-representation $R\Gamma_v(\widetilde{X}, \widehat{\s{O}}_X)^{RG-\rm la}$ is locally analytic, in particular 
\[
R\Gamma_v(\widetilde{X}, \widehat{\s{O}}_X)^{R(G\times H)-\rm la}\xrightarrow{\sim } R\Gamma_v(\widetilde{X}, \widehat{\s{O}}_X)^{RG-\rm la}
\]
 is an equivalence.

 \item The above results in (1), (2), (3) hold for $R\Gamma_{v}(\widetilde{X}, \bb{B}_I)$ instead of $R\Gamma_v(\widetilde{X}, \widehat{\s{O}}_X)$, for any compact interval with rational endpoints $I=[s,r]\subset (0,\infty)$. 

\end{enumerate}

\end{theorem}

An immediate corollary of the theorem above is that the $\bb{B}_I$-cohomology groups of qcqs smooth rigid spaces endowed with an action of a $p$-adic Lie group are locally analytic (simply take $G=1$ in the above theorem):

\begin{corollary}\label{CohoEquivRigidSpaces}
Let $X$ be a  smooth quasi-compact and separated rigid space over $\mathbb{C}_p$, endowed with a continuous action of a $p$-adic Lie group $H$. Then for any compact interval $I\subset (0,\infty)$ with rational endpoints the $v$-cohomologies $R\Gamma_v(X, \widehat{\s{O}}_X)$ and 
$R\Gamma_v(X, \bb{B}_I)$ are solid locally analytic $H$-representations.
\end{corollary}

\begin{corollary}
Let $X$ be a smooth quasi-compact and separated rigid space over $\mathbb{C}_p$, endowed with a continuous action of a $p$-adic Lie group $H$. Then 
$R\Gamma_{\proet}(X, \bb{Q}_p)$ has a natural structure of locally analytic $H$-representation. 
\end{corollary}

\begin{proof}
By \cite[Proposition II.2.5]{FarguesScholze} we have a short exact sequence of pro\'etale sheaves
\[
0\to \bb{Q}_p\to \bb{B}_{[1,p]} \xrightarrow{\varphi-1} \bb{B}_{[1,1]}\to 0.
\]
This yields a fiber sequence\footnote{We implicitly use here that since $X$ is qcqs, \'etale, pro\'etale and $v$-cohomology 
with $\bb{Q}_p$-coefficients are the same.}
$$R\Gamma_{\proet}(X, \bb{Q}_p)\to R\Gamma_v(X, \bb{B}_{[1,p]})\to R\Gamma_v(X, \bb{B}_{[1,1]})$$
and the result follows from the previous corollary. 
\end{proof}

  In the remaining sections we fix the setup in  \Cref{theo:LocAnTwoTowers}. Let us remark that part (3) is a direct consequence of parts (1) and (2): writing 
  $C^{\rm la}(G, \bb{Q}_p)=\varinjlim_{n} V_n$ as a countable colimit of locally analytic representations of 
  $G$ on $\bb{Q}_p$-Banach spaces and using that solid tensor product and $R\Gamma(G,-)$ commute with colimits, we have
  $$R\Gamma_{v}(\widetilde{X}, \widehat{\s{O}}_X)^{RG-\rm la}=R\Gamma(G, R\Gamma_v(\widetilde{X}, \widehat{\s{O}}_X)\otimes^L_{\bb{Q}_p, \square} C^{\rm la}(G, 
  \bb{Q}_p))\simeq \varinjlim_{n} R\Gamma(G, R\Gamma_v(\widetilde{X}, \widehat{\s{O}}_X)\otimes^L_{\bb{Q}_p, \square} V_n)$$
  and each object in the last colimit is a locally analytic representation of $H$ by part (2). Furthermore, the case of the sheaves $\bb{B}_{I}$ will imply the case for $\widehat{\s{O}}_X$ by taking a short exact sequence 
  \begin{equation}\label{eqojoajepaeg}
  0\to \bb{B}_I\xrightarrow{\xi} \bb{B}_{I}\to \widehat{\s{O}}_X\to 0 
  \end{equation}
  where $\xi$ is a generator of the map $\theta\colon \bb{A}_{\inf}(\bb{C}_p)\to \n{O}_{\bb{C}_p}$, and $I$ is any interval such that $\Spa \bb{B}_I(\bb{C}_p)$ contains $\Spa \bb{C}_p$.
  
\subsection{$v$-cohomology as solid quasi-coherent cohomology}\label{ss:SolidQuasiCohCohomology}  
  
  In order to see the $v$-cohomology groups of \Cref{theo:LocAnTwoTowers} as solid abelian groups we use the six functor formalism of  quasi-coherent sheaves on Fargues--Fotaine curves following \cite{Anschutz-LeBras-Mann}. We will not recall all the details, only the minimum that is needed for us. Given a $v$-stack $S$ there is a category of solid quasi-coherent sheaves on the different curves $\n{Y}_S$, $Y_S$  and $X_S$ associated to $S$, in the notation of \cite[Section 4]{Anschutz-LeBras-Mann} this corresponds to the six functor formalisms sending $S$ to   $\n{D}_{[0,\infty]}(S)$, $\n{D}_{(0,\infty)}(S)$ and $\n{D}_{\ob{FF}}(S)$ respectively.

 \begin{lemma}\label{LemmaSolidCohoisvCoho}
Let $G$ be a compact $p$-adic Lie group. There is a natural symmetric monoidal equivalence 
\[
\n{D}_{(0,\infty)}(\Spd (\bb{F}_p)/\underline{G})=\ob{D}(\AnSpec \bb{Q}_p/ G^{\ob{prof}})
\]
where the right term is as in \cite[Definition 6.1.1]{RJRCSolidLocAn2}, equivalently, the right term is given by $C(G,\bb{Q}_{p})$-comodules on solid $\bb{Q}_p$-vector spaces. In particular, after applying the forgetful functor $\ob{D}(\AnSpec \bb{Q}_p/ G^{\ob{prof}})\to \ob{D}(D^{\ob{cont}}(G,\bb{Q}_p))$ obtained by sending a $C(G,\bb{Q}_p)$-comodule to its underlying $D^{\ob{cont}}(G,\bb{Q}_p)$-module structure, we have a natural symmetric monoidal functor 
\[
\n{D}_{(0,\infty)}(\Spd (\bb{F}_p)/\underline{G})\to \ob{D}(D^{\ob{cont}}(G,\bb{Q}_p)).
\]
\end{lemma}  
  \begin{proof}
  The first statement is \cite[Corollary 6.4.3]{Anschutz-LeBras-Mann}. The forgetful functor $\ob{D}(\AnSpec \bb{Q}_p/ G^{\ob{prof}})\to \ob{D}(D^{\ob{cont}}(G,\bb{Q}_p))$ is first defined at the level of abelian categories by \cite[Lemma 6.1.3]{RJRCSolidLocAn2} and then by passing to derived $\infty$-categories. 
  \end{proof}
  
\begin{corollary}\label{CorollaryFFCurve}
Let $G$ be a compact $p$-adic Lie group. There is a natural equivalence of symmetric monoidal categories 
\[
\n{D}_{(0,\infty)}(\Spd \bb{C}_p)=\ob{D}_{\sol}(Y_{\bb{C}_p}/G^{\ob{prof}})
\]
where the right term is the category of solid  quasi-coherent sheaves on the relative classifying stack of $G^{\ob{prof}}$ over the adic space $Y_{\bb{C}_p}$. More precisely, 
\[
\ob{D}_{\sol}(Y_{\bb{C}_p}/G^{\ob{prof}})=\ob{coMod}_{C(G,\bb{Q}_p)}(\ob{D}_{\sol}(Y_{\bb{C}_p}))
\]
is the category of $C(G,\bb{Q}_p)$-comodules on solid quasi-coherent sheaves on the adic space $Y_{\bb{C}_p})$. In particular, we have a forgetful functor 
\[
\ob{D}_{\sol}(Y_{\bb{C}_p}/G^{\ob{prof}})=\ob{coMod}_{C(G,\bb{Q}_p)}(\ob{D}_{\sol}(Y_{\bb{C}_p}))\to \ob{Mod}_{D^{\ob{cont}}(G,\bb{Q}_p)}(\ob{D}_{\sol}(Y_{\bb{C}_p}))
\]
towards $D^{\ob{cont}}(G,\bb{Q}_p)$-modules on $\ob{D}_{\sol}(Y_{\bb{C}_p})$.
\end{corollary}
\begin{proof}
The first statement exactly by the same argument as  \cite[Corollary 6.4.3]{Anschutz-LeBras-Mann} noticing that $\ob{D}_{(0,\infty)}(\Spd \bb{C}_p)$ is precisely $\ob{D}_{\widehat{\sol}}(Y_{\bb{C}_p})=\ob{D}_{\sol}(Y_{\bb{C}_p})$ by  Proposition 2.1.2 (i) of  \textit{loc. cit.}. For the second statement one can argue in different ways.  Explicitly, one writes $\ob{D}_{\sol}(Y_{\bb{C}_p})$ as limit of categories of quasi-coherent sheaves of affinoids $Y_{I,\bb{C}_p}$ with $I\subset (0,\infty)$ a closed interval with rational ends reducing the question to $Y_{I,\bb{C}_p}$. In that case, one first  do the passage from comodules to modules at the level of abelian categories first and then to derived categories as in \Cref{LemmaSolidCohoisvCoho}. A more sophisticated and soft argument is by taking base change of the $\ob{D}(\bb{Q}_{p,\sol})$-linear forgetful functor $\ob{coMod}_{C(G,\bb{Q}_p)}(\ob{D}(\bb{Q}_{p,\sol}))\to \ob{D}(D^{\ob{cond}}(G,\bb{Q}_p))$ along the symmetric monoidal functor $\ob{D}(\bb{Q}_{p,\sol})\to \ob{D}_{\sol}(Y_{\bb{C}_p})$ with respect to Lurie's tensor product in $\mathrm{Pr}^L$, see \cite[Section 4.8]{HigherAlgebra}.
\end{proof}

\begin{lemma}\label{PropositionClassifyingStack}
Let $X$ be a  quasi-compact and separated compactifiable bdcs diamond over $\Spd \bb{C}_p$, let $\overline{X}$ denote Huber's compactification over $\bb{C}_p$ and  consider the natural map $f\colon \overline{X}\to \Spd \bb{C}_p$. Then the map $f$ is proper and bdcs, hence  cohomologically proper for the six functor formalism $\n{D}_{(0,\infty)}$. Moreover, let $I\subset (0,\infty)$ be a closed interval with rational ends and let $\bb{B}_{I}\in \n{D}_{(0,\infty)}(\overline{X})$ be the sheaf that sends an affinoid perfectoid $S\to \overline{X}$ to the nuclear sheaf on $Y_{S}$ given by $\bb{B}_I(S)$.  Then $f_* \bb{B}_I\in \ob{D}_{\sol}(Y_{\bb{C}_p})$ is supported in any strict affinoid subspace $Y_{J,\bb{C}_p}$ for $I\subsetneq J$. Moreover, under the forgetful functor $\ob{D}_{\sol}(Y_{J,\bb{C}_p})\to \ob{D}(\bb{Q}_{p,\sol})$ given by global sections, $f_* \bb{B}_{I}$ is naturally identified with $R\Gamma_v(X,\bb{B}_I)$.

\end{lemma}
\begin{remark}\label{RemarkNuclearSheaves}
 Let $S$ be an affinoid perfectoid space over $\bb{F}_p$. By  \cite[Proposition 2.1.2]{Anschutz-LeBras-Mann} (ii) there is not distinction between nuclear sheaves on $\ob{D}_{\sol}(Y_{S})$ and $\ob{D}_{\widehat{\sol}}(Y_{S})$, hence the functor $S\mapsto \bb{B}_{I}(S)$ also defines an object in $\ob{D}_{\widehat{\sol}}(Y_S)$ whose formation satisfies $v$-descent, eg. by \Cref{LemmaVanishingCohoBI}. 
\end{remark}
\begin{proof}
The fact that $f$ is cohomologically proper follows from \cite[Corollary 3.3.5]{Anschutz-LeBras-Mann} and \cite[Theorem 1.2.4 (ii)]{MannSix} since the map $f\colon \overline{X}\to \Spd \bb{C}_p$ is a bdcs proper map of $v$-stacks (see Definition 3.6.9 of \cite{MannSix} for the definition of bdcs map of stacks). Indeed, the fact that $f$ is bdcs follows from \cite[Lemma 3.5.10]{MannSix} and the first step in the proof of Theorem 3.6.12 of \textit{loc. cit.}   It  is clear that $f_*\bb{B}_I$ is supported in $Y_{J,\bb{C}_p}$ since it contains the closure of $Y_{I,\bb{C}_p}$ in $Y_{\bb{C}_p}$. For the identification of $f_*\bb{B}_{I}$ with $v$-cohomology,  by descent it suffices to notice that if $S$ is an affinoid perfectoid space over $\Spd(\bb{C}_p)$ and $g\colon S\to \Spd(\bb{C}_p)$ is the structural map, then $g_* \bb{B}_{I}= \bb{B}_I(S)$ which is obvious by the construction of the six functor formalism. 
\end{proof}

  \begin{proposition}\label{PropositionProofTheorem1}
  Let $X$ be a quasi-compact and separated bdcs diamond  over $\Spd \bb{C}_p$ endowed with an action of a compact $p$-adic Lie group $H$, and let $\overline{X}$ denote Huber's compactification of $X$ over $\bb{C}_p$. Consider the map of diamonds $f\colon \overline{X}/\underline{H}\to \Spd(\bb{C}_p)/\underline{H}$.  Then, under the equivalences of categories \Cref{CorollaryFFCurve}, the pullback to $Y_{\bb{C}_p}$ of  $f_* \bb{B}_I$ is given by $R\Gamma_v(X,\bb{B}_I)$. Thus, $f_*\bb{B}_I$ induces a structure of solid $H$-representation to the $v$-cohomology of $\bb{B}_I$ over $X$.   Applying this to the diamond $\widetilde{X}$ endowed with the $G\times H$ action of \Cref{theo:LocAnTwoTowers} we get part (1) of the theorem.   
  \end{proposition}
  \begin{proof}
The map $f\colon \overline{X}/\underline{H}\to \Spd \bb{C}_p /\underline{H}$ is proper and bdcs after \Cref{PropositionClassifyingStack}, hence it is cohomologically proper after \cite[Theorem 3.6.12]{MannSix}. It follows that $f_*$ satisfies base change. The fact that $f_*\bb{B}_I$ is just $R\Gamma_v(X,\bb{B}_I)$ endowed with it natural $H$-module structure follows from \Cref{PropositionClassifyingStack} by base change along the map $\Spd \bb{C}_p\to \Spd \bb{C}_p/\underline{H}$.  
  \end{proof}

\begin{remark}\label{RemarkExplicitDescriptionAction}
Keep the notation of \Cref{PropositionProofTheorem1}. We give a more explicit description of the solid structure of the representation $R\Gamma_v(X,\bb{B}_I)$. We can write $\bb{B}_I=\bb{B}_{I}^+[\frac{1}{[\varpi]}]$ where $\varpi$ is a pseudo-uniformizer of $\bb{C}_p^{\flat}$. Since $X$ is quasi-compact we have that $R\Gamma_v(X,\bb{B}_I)=R\Gamma_v(X,\bb{B}_I^+)[\frac{1}{[\varpi]}]$. Then, $R\Gamma_v(X,\bb{B}_I^+)$ is naturally a derived $[\varpi]$-complete module, that is the natural map  $R\Gamma_v(X,\bb{B}_I^+)=\varprojlim_{n}R\Gamma_v(X,\bb{B}_I^+/[\varpi]^n)$ is an equivalence. Now, since $X$ is a  locally special diamond, we have that $R\Gamma_v(X,\bb{B}_I^+/\varpi^n)=R\Gamma_{\et}(X,\bb{B}^+_I/[\varpi]^n)$. In total, we have that 
\[
R\Gamma_v(X,\bb{B}_I)=\bigg(\varprojlim_n \big( R\Gamma_{\et}(X,\bb{B}^+_I/[\varpi]^n) \big) \bigg) [\frac{1}{[\varpi]}].
\]
The terms  $R\Gamma_{\et}(X,\bb{B}^+_I/[\varpi]^n)$  are naturally seen as discrete $\bb{Z}_p$-modules  endowed with a smooth $H$-action, hence $R\Gamma_v(X,\bb{B}_I)$ is essentially a  Banach representation in a suitable derived sense, i.e. it is the generic fiber of a \textit{derived $p$-complete} $H$-representation whose special fiber is discrete. \Cref{PropositionProofTheorem1} says that the previous ad hoc construction of the $H$-action is formally obtained from a suitable six functor formalism. 
\end{remark}

\subsection{Independence of locally analytic vectors}\label{ss:IndependencelocAn}  
  
  In this section we shall prove parts (2) and (3) of \Cref{theo:LocAnTwoTowers}.  We recall that it suffices to prove it for the sheaves $\bb{B}_I$ thanks to the fiber sequence \eqref{eqojoajepaeg}.

\begin{lemma}\label{LemmaIndependenceVectors}
Keep the notation of \Cref{theo:LocAnTwoTowers}. That is, $X$ is a quasi-compact and separated smoth adic space over $\bb{C}_p$ endowed with an action of a $p$-adic Lie group $H$, and $\widetilde{X}\to X$ is a $H$-equivariant $G$-torsor with $G$ a $p$-adic Lie group. Then part (3) of the theorem implies part (2).
\end{lemma}  
\begin{proof}
Suppose that  (3) of the theorem holds, that is, that the natural map 
\[
R\Gamma_v(\widetilde{X}, \bb{B}_I)^{R(G\times H)-\rm la}\xrightarrow{\sim} R\Gamma_v(\widetilde{X}, \bb{B}_I)^{RG- \rm la}
\]
is an equivalence. Let $V$ be a locally analytic $\bb{Q}_p$-Banach representation of $G$. We want to show that $R\Gamma(G, V\otimes^L_{\bb{Q}_p,\sol} R\Gamma_v(\widetilde{X},\bb{B}_I))$ is a locally analytic representation of $H$. We have natural equivalences of  $H$-representations 
\[
\begin{aligned}
R\Gamma(G, V\otimes^L_{\bb{Q}_p,\sol} R\Gamma_v(\widetilde{X},\bb{B}_I)) & = R\Gamma(G, V\otimes^L_{\bb{Q}_p,\sol} R\Gamma_v(\widetilde{X},\bb{B}_I)^{RG-\rm la}) \\
& =   R\Gamma(G, V\otimes^L_{\bb{Q}_p,\sol} R\Gamma_v(\widetilde{X},\bb{B}_I)^{R(H\times G)-\rm la})
\end{aligned}
\]
where the first equivalence follows from  the comparison theorems for group cohomology \Cref{PropositionComparisonLocAnCohomology} and the projection formula of locally analytic vectors \Cref{PropStabilityLocAn} (2), and the second equivalence holds by assumption. It is clear the the last representation of $H$ is locally analytic, proving what we wanted. 
\end{proof}

\begin{proposition}\label{PropositionCompatibilityLocAn}
Keep the notation of \Cref{theo:LocAnTwoTowers}. Then the natural map 
\[
R\Gamma_v(\widetilde{X}, \bb{B}_I)^{R(G\times H)-\rm la}\xrightarrow{\sim} R\Gamma_v(\widetilde{X}, \bb{B}_I)^{RG- \rm la}
\]
is an equivalence, 
\end{proposition}
\begin{proof}
We want to show that the action of $H$ on $R\Gamma_v(\widetilde{X}, \bb{B}_I)^{RG- \rm la}$ is locally analytic. By definition we have that 
\[
R\Gamma_v(\widetilde{X}, \bb{B}_I)^{RG- \rm la}=R\Gamma(G, C^{la}(G,\bb{Q})_p\otimes^L_{\bb{Q}_{p},\sol} R\Gamma_v(\widetilde{X}, \bb{B}_I)).
\] 
Writing $X=\bigcup_{i\in J} U_i$ as a finite union of open affinoid subspaces, there is an open subgroup $H_0\subset H$ acting on  any finite intersection $U_{J_0}$ of the  $U_i$'s with $J_0\subset J$. Hence, we have an isomorphism of $G\times  H_0$-representations
\[
R\Gamma_v(\widetilde{X},\bb{B}_I)=\varprojlim_{J_0\subset J} R\Gamma_v(\widetilde{X}_{U_{J_0}}, \bb{B}_I).
\]
Since locally analytic representations are stable under finite limits, we can shrink $X$ and assume without loss of generality that it is affinoid and even that it admits an \'etale map $X\to \bb{T}^d_{\bb{C}_p}$ to a torus that factors as a finite composite of  finite \'etale maps and rational localizations.

Writing $C^{la}(G,\bb{Q})_p=\varinjlim_{h} C^h(G,\bb{Q}_p)$ as colimit of spaces of analytic functions of $G$ (in particular of analytic Banach representations of $G$), we have that 
\[
R\Gamma_v(\widetilde{X}, \bb{B}_I)^{RG- \rm la}=\varinjlim_{h} R\Gamma(G, C^h(G,\bb{Q}_p)\otimes^L_{\bb{Q}_{p},\sol} R\Gamma_v(\widetilde{X}, \bb{B}_I)),
\]
hence it suffices to prove that each term $C^h(G,\bb{Q}_p)\otimes^L_{\bb{Q}_{p},\sol} R\Gamma_v(\widetilde{X}, \bb{B}_I)$ as a locally analytic action of $H$. Write $C^h(G,\bb{Q}_p)= V_h[\frac{1}{p}]$ where $V_h=C^h(G,\bb{Z}_p)$ is the space of bounded  power functions of  $C^h(G,\bb{Q}_p)$, in particular $V_h$ is a $p$-complete and torsion free $G$-representation for which the action on $V_h/p^s$ factors through a finite quotient of $G$ for all $s> 0$. Now, write $I=[s,r]\subset (0,\infty)$ with $s,r$ rational numbers. 

By \cite[Proposition 2.12.10 (i)]{MannSix} the tensor product $V_{h}\otimes^L_{\bb{Z}_p,\sol} \bb{B}^+_I$ is derived $p$-adically complete, since $[\varpi]$ is another pseudo-uniformizer of $\bb{B}^+_I$ comparable with $p$ (i.e. $|[\varpi]^{1/s}|\leq |p|\leq |[\varpi]^{1/r}|$), the tensor is also $[\varpi]$-adically complete. Furthermore, its mod $[\varpi]^{1/r}$-fiber is equal to the tensor of discrete sheaves
\[
(V_{h}\otimes^L_{\bb{Z}_p,\sol} \bb{B}^+_I)/[\varpi^{1/r}]= (V_{h}/p)\otimes^L_{\bb{Z}/p} \bb{B}^+_I/[\varpi]^{1/r}. 
\]
Let $G_0\subset G$ be an open subgroup such that the action of $G_0$ on $V_{h}/p$ is trivial and isomorphic to a direct sum of copies of $\bb{Z}/p$. Then, by \Cref{bpluscor1} we have a natural $G_0\times H$-equivariant almost isomorphism of sheaves on $\widetilde{X}$
\[
(V_{h}\otimes^L_{\bb{Z}_p,\sol} \bb{B}^+_I)/[\varpi^{1/r}] =^a \bigoplus_i \widehat{\s{O}}^+_{\widetilde{X}}/\pi^{1/r}
\]
  where $\pi^{1/r}\in \n{O}_{\bb{C}_p}$ is any element whose fiber mod $p$ is equivalent to $\varpi^{1/r}\in \n{O}_{\bb{C}_p^{\flat}}$. Hence, by  the locally analytic criterion \Cref{PropLocAnCriterion} it suffices to show that there is an open compact subgroup $H_0\subset H$ such that the $\bb{B}_{I}^+(\bb{C}_p)$-linear $H_0$-representation 
  \[
  R\Gamma(G_0, R\Gamma_v(\widetilde{X}, \widehat{\s{O}}^+_{\widetilde{X}}/\pi^{1/r}))
  \]
  is a $\pi^{1/2r}$-retract of a trivial $H_0$-representation.  After modifying $X$ by a finite \'etale cover determined by the quotient $G\to G/G_0$, we can assume that $G=G_0$ so that  
  \[
    R\Gamma(G_0, R\Gamma_v(\widetilde{X}, \widehat{\s{O}}^+_{\widetilde{X}}/\pi^{1/r}))=^{a} R\Gamma_v(X, \widehat{\s{O}}^+_{\widetilde{X}}/\pi^{1/r})). 
  \]

   Let $T_{1},\ldots, T_d$  be the coordinates of the toric chart  $X\to \bb{T}^d_{\bb{C}_p}$. Let $X_n$ be the $\bb{Z}(1)^d/p^n$-torsor obtained by taking $p^n$-th power roots of $T$ and let $X_{\infty}=\varprojlim_n X_n$ be the perfectoid $\bb{Z}_p(1)^{d}$-torsor over $X$. Then 
   \[
   R\Gamma_v(X, \widehat{\s{O}}^+_{\widetilde{X}}/\pi^{1/r})) = R\Gamma(\bb{Z}_p(1)^d,  \widehat{\s{O}}^+(X_{\infty})/\pi^{1/r})
   \]is computed via a Kozsul complex. For $n\geq 1$ consider the cofiber  $Q(n)=\ob{cofib}(\s{O}^+(X_n)\to \widehat{\s{O}}^+(X_{\infty}))$ where $\s{O}^+(X_n)\subset \s{O}(X_n)$ is the ring of integral elements of $X_n$. Then, the bounds arising from geometric Sen theory \cite[Definition 2.2.1 and Proposition 2.2.14]{RCGeoSen} imply that there is some $n\gg 0$ such that  the multiplication by $\pi^{1/2r}$ on $R\Gamma(\Z_p(1)^d, Q(n))$ is homotopic to zero and therefore multiplication by $\pi^{1/2r}$ on  $R\Gamma_v(X, \widehat{\s{O}}_{X}^+/\pi^{1/r})$ factors through a map 
   \[
   R\Gamma_v(X, \widehat{\s{O}}_{X}^+/\pi^{1/r})\to R\Gamma(\bb{Z}_p(1)^d, \s{O}^+(X_n)/\pi^{1/r})\to  R\Gamma_v(X, \widehat{\s{O}}_{X}^+/\pi^{1/r}). 
   \]
   Now, the action of $\bb{Z}_p(1)^d$ on $\s{O}^+(X_n)$ factors through the finite quotient $\bb{Z}_p(1)^d/p^n$, and by \cite[Lemma 2.3]{ScholzeLubinTate} there is an open subgroup $H_0(n)\subset H$ that lifts to an action on $X_n$ that commutes with that of $\bb{Z}_p(1)^d/p^n$. This implies that the map 
   \[
   R\Gamma(\bb{Z}_p(1)^d, \s{O}^+(X_n)/\pi^{1/r})\to  R\Gamma_v(X, \widehat{\s{O}}_{X}^+/\pi^{1/r})
   \]
   is $H^0(n)$-equivariant. It is left to see that we can pick $H_0(n)$ so that $ R\Gamma_v(X, \widehat{\s{O}}_{X}^+/\pi^{1/r})\to R\Gamma(\bb{Z}_p(1)^d, \s{O}^+(X_n)/\pi^{1/r})$ is also $H^0(n)$-equivariant, for that, notice that by the Sen axiom \cite[Definition 2.2.1 (CST2)]{RCGeoSen}, given $l\gg 0$ there is $n\gg 0$ such that the cohomology groups of the fiber of the map $R\Gamma(\bb{Z}_p(1)^d, \s{O}^+(X_n)/\pi^{1/r})\to  R\Gamma_v(X, \widehat{\s{O}}_{X}^+/\pi^{1/r})$ are killed by $\pi^{1/l}$, since this fiber has  cohomology groups supported in cohomological degrees $[0,d]$ in the almost category, we observe that $\pi^{d/l}$ acts homotopically equivalently to zero as $H^0(n)$-representation (this follows from an induction step as multiplication by $\pi^{1/l}$ kills the top cohomology group of the given complex $H^0(n)$-equivariantly). Taking $l$ such that $l/d>2r$, we deduce that multiplication by $\pi^{1/2r}$ on  $R\Gamma_v(X, \widehat{\s{O}}_{X}^+/\pi^{1/r})$ will factor through a $H_0(n)$-equivariant map 
   \[
   R\Gamma_v(X, \widehat{\s{O}}_{X}^+/\pi^{1/r})\to R\Gamma(\bb{Z}_p(1)^d, \s{O}^+(X_n)/\pi^{1/r})
   \]
proving that $ R\Gamma_v(X, \widehat{\s{O}}_{X}^+/\pi^{1/r})$  is an $H_0(n)$-equivariant  $\pi^{1/r}$-retract of $R\Gamma(\bb{Z}_p(1)^d, \s{O}^+(X_n)/\pi^{1/r})$ as wanted. To finish the proof by applying \Cref{PropLocAnCriterion}, we need to show that there is some $H_0$ such that the restriction to $R\Gamma(\bb{Z}_p(1)^d, \s{O}^+(X_n)/\pi^{1/r})$ is trivial, but this follows from the fact that $X_n$ is of finite type over $\bb{C}_p$, so that  $H_0(n)$ acting on $\s{O}^+(X_n)/\pi^{1/r}$ factors through a finite quotient. Replacing $H_0(n)$ by a smaller compact open subgroup, we can assume that $H_0(n)$ acts trivially on $\s{O}^+(X_n)/\pi^{1/r}$ and thus on the complex $R\Gamma(\bb{Z}_p(1)^d, \s{O}^+(X_n)/\pi^{1/r})$ proving what we wanted. 
\end{proof}

\section{Local Shimura varieties, preliminaries}
\label{SectionLocalShimura}

In this section we recall the construction and main properties of local Shimura varieties, departing from \cite{ScholzeWeinspadicgeometry}
and using instead the geometry of the stack of $\bbf{G}$-bundles on the Fargues-Fontaine curve \cite{FarguesScholze} (the resulting objects are of course the same). 
This chapter contains no original result and its only purpose is to fix notations and recall the key results that we will need in the next chapter. 
 
\subsection{Group theory}\label{groupth}
    Fix once and for all an algebraic closure $\overline{\bb{Q}}_p$ of $\bb{Q}_p$ and let $$\Gamma:={\rm Gal}(\overline{\bb{Q}}_p/\bb{Q}_p).$$ Let 
   $k=\overline{\bb{F}}_p$ be the residue field of $\overline{\bb{Q}}_p$ and let $\breve{\bb{Q}}_p=W(k)[1/p]\subset \widehat{\overline{\bb{Q}}}_p$, endowed with its natural lift $\sigma$ of the absolute Frobenius of $k$. 

 Fix once and for all a connected reductive group $\bbf{G}$ over $\bb{Q}_p$ and let ${\rm Rep}(\bbf{G})$ be the category of finite dimensional algebraic representations of $\bbf{G}$.
  Fix a maximal torus $\bbf{T}$ inside a Borel subgroup $\bbf{B}$ of $\bbf{G}_{\overline{\bb{Q}}_p}$. The set 
    $X_*(\bbf{T})^+$ of dominant (with respect to $\bbf{B}$) cocharacters of $\bbf{T}$ is endowed with 
   the Bruhat partial order.\footnote{Recall that
     $\nu\leq \nu'$ if $\nu'-\nu$ is a linear combination of positive coroots (with respect to $\bbf{B}$) with non-negative integral coefficients.}
     
     Let $X_*(\bbf{G})/\bbf{G}$ be the set of $\bbf{G}(\overline{\bb{Q}}_p)$-conjugacy classes $\{\mu\}$ of cocharacters $\mu: \bb{G}_{m, \overline{\bb{Q}}_p}\to \bbf{G}_{\overline{\bb{Q}}_p}$. Any conjugacy class $\{\mu\}\in X_*(\bbf{G})/\bbf{G}$ 
    intersects $X_*(\bbf{T})^+$ in a unique element, giving a bijection 
     $$X_*(\bbf{G})/\bbf{G}\simeq X_*(\bbf{T})^+.$$  
     The obvious action of $\Gamma$ on $X_*(\bbf{G})/\bbf{G}$ induces an action on $X_*(\bbf{T})^+$. 
          
      Let $\mu: \bb{G}_{m, \overline{\bb{Q}}_p}\to \bbf{G}_{\overline{\bb{Q}}_p}$ be a geometric cocharacter. It induces a \emph{decreasing} filtration 
      $F^{\bullet}_{\mu}(V)$ on any $V\in {\rm Rep}(\bbf{G}_{\overline{\bb{Q}}_p})$, defined by 
      $$F^p_{\mu}(V)=\bigoplus_{n\leq -p} V[n],$$
      where $V[n]$ is the subspace on which $\mu(t)$ acts by $t^{n}$. 
      Let $\bbf{P}_{\mu}$ be the parabolic subgroup of $\bbf{G}_{\overline{\bb{Q}}_p}$ which stabilizes 
      the filtration $F^{\bullet}_{\mu}(V)$ for all $V\in {\rm Rep}(\bbf{G}_{\overline{\bb{Q}}_p})$, we have 
      $$\bbf{P}_{\mu}=\{g\in \bbf{G}_{\overline{\bb{Q}}_p}|\, \lim_{t\to \infty} \mu(t) g\mu(t)^{-1}\,\, \text{exists}.\}$$
      Similarly, we can define an increasing $\mu$-filtration 
      \[
      F_{\mu,p}=\bigoplus_{n\geq -p} V[n]
      \] with stabilizing parabolic $\bbf{P}_{\mu^{-1}}$ given by those $g\in \bbf{G}$ such that $\lim_{t\to 0}\mu(t)g\mu^{-1}(t)$  exists.  Thus, if $\mu\in X_*(\bbf{T})^+$ is dominant,  then $\bbf{P}_{\mu^{-1}}$ is a standard parabolic, i.e. it contains $\bbf{B}$. 
      
\begin{remark}\label{RemarkConventionsGroups} Our conventions agree with those of \cite[Section 19.4]{ScholzeWeinspadicgeometry}. Our groups $\bbf{P}_{\mu}$ and $\bbf{P}_{\mu^{-1}}$ correspond to $\bbf{P}_{\mu}^{\ob{std}}$ and $\bbf{P}_{\mu}$ in the notation of \cite{CaraianiScholze2017} respectively. In particular, if $(\bbf{G},X)$ is a global Shimura datum over $\bb{Q}$ and $h\colon \mathrm{Res}^{\bb{C}}_{\bb{R}} \bb{G}_m\to \bbf{G}_{\bb{R}}$ is a fixed element in $X$, the Hodge cocharacter $\mu$ is the restriction to the first component of $h\otimes_{\bb{R}}\bb{C}\to \bbf{G}_{\bb{C}}$. By convention, a real $\mathrm{Res}^{\bb{C}}_{\bb{R}} \bb{G}_m$-representation $V$ has a Hodge  decomposition $V_{\bb{C}}=\bigoplus_{(p,q)} V_{p,q}$ where the character of $V_{p,q}$ is precisely $(z^{-p}, \overline{z}^{-q})$. Hence,  $F^p_{\mu} V_{\bb{C}}=\bigoplus_{p'\geq p} V_{p',\bullet}$ and $F_{\mu,p}V=\bigoplus_{p'\leq p} V_{p',\bullet}$. 
\end{remark}

      The variety 
      $\bbf{G}_{\overline{\bb{Q}}_p}/\bbf{P}_{\mu}$ has $\overline{\bb{Q}}_p$-points that parametrize filtrations 
      $\bbf{G}(\overline{\bb{Q}}_p)$-conjugate to $F^{\bullet}_{\mu}$, and admits a model 
      $\FL_{\mu}$ defined over $E=E(\{\mu\})$, the reflex field of the conjugacy class $\{\mu\}$ (i.e. its minimal field of definition), a finite extension of 
      $\bb{Q}_p$ contained in $\overline{\bb{Q}}_p$. We denote 
      $$\Fl_{\mu}=\FL_{\mu}^{\rm ad}\to {\rm Spa}(E)$$
      the associated adic space. Notice that the map $g\mapsto g^{-1}$ produces an isomorphism $\bbf{G}_{\overline{\bb{Q}}_p}/\bbf{P}_{\mu}\cong \bbf{P}_{\mu}\backslash \bbf{G}_{\overline{\bb{Q}}_p}$.

   \subsection{The $\bb{B}_{\dR}$-affine Grassmannian} 
   
   It is not difficult to see that the presheaf $L^{(+)} \bbf{G}=\bbf{G}(\bb{B}_{\dR}^{(+)})$ over $\bb{Q}_p^{\lozenge}$ sending 
   an affinoid perfectoid $S/\bb{Q}_p^{\lozenge}$ to $\bbf{G}(\bb{B}_{\dR}^{(+)}(S^{\sharp}))$ is a $v$-sheaf.\footnote{This holds if we replace $\bbf{G}$ by any affine algebraic variety over $\bb{Q}_p$.} Here and below
  $S^{\sharp}$ is the untilt of $S$ over $\bb{Q}_p$ induced by the map $S\to \bb{Q}_p^{\lozenge}$. There is a natural map of $v$-sheaves over 
  $\bb{Q}_p^{\lozenge}$
  $$\theta: L^+\bbf{G}\to \bbf{G}^{\diamond},$$
  induced by  $\theta: \bb{B}_{\dR}^{+}(S^{\sharp})\to  \mathcal{O}(S^{\sharp})$. 
   The 
\emph{$\bb{B}_{\dR}$-affine Grassmannian} for $\bbf{G}$ is the quotient 
$${\rm Gr}_{\bbf{G}}=L\bbf{G}/L^+\bbf{G},$$ taken either in the category of \'etale or of $v$-sheaves (see \cite[Proposition  20.2.2]{ScholzeWeinspadicgeometry}).

  If $S={\rm Spa}(R, R^+)$ is affinoid perfectoid over $\bb{Q}_p$ then ${\rm Gr}_{\bbf{G}}(S^{\lozenge})$ depends only on $R$ (thus ${\rm Gr}_{\bbf{G}}$ is partially proper) and 
  we will simply write ${\rm Gr}_{\bbf{G}}(R)$ instead of ${\rm Gr}_{\bbf{G}}(S^{\lozenge})$. By \cite[Definition 20.2.1 and Proposition 20.2.2]{ScholzeWeinspadicgeometry}, ${\rm Gr}_{\bbf{G}}(R)$ is identified with the 
  set of isomorphism classes of pairs consisting of a $\bbf{G}$-torsor on ${\rm Spec}(\bb{B}^+_{\dR}(S^{\sharp}))$ together with a trivialization of its restriction to 
  ${\rm Spec}(\bb{B}_{\dR}(S^{\sharp}))$.

     Suppose that $S={\rm Spa}(C,C^+)/ \Spa (\bb{Q}_p)$ 
      is a geometric point. Then $B^+:=\bb{B}_{\dR}^+(C^{\sharp})$ is a strictly henselian dvr
       and we let $\xi$ be a uniformizer and 
      $B=\bb{B}_{\dR}(C)=B^+[1/\xi]$.  We let $\xi^{\mu}$ be the image of $\mu(\xi)$ in ${\rm Gr}_{\bbf{G}}(C)$. 
               The Cartan decomposition yields a canonical decomposition
      $${\rm Gr}_{\bbf{G}}(C)=\coprod_{\lambda\in X_*(\bbf{T})^+} {\rm Gr}_{\bbf{G}, \lambda}(C),$$
      where ${\rm Gr}_{\bbf{G}, \lambda}(C)$ is the left $\bbf{G}(B^+)$-orbit of $\xi^{\mu}$  (it is independent of the choice of
      the uniformizer $\xi$).\footnote{Our conventions are those of \cite[Section 19.2]{ScholzeWeinspadicgeometry}.} 
      
           Consider the subsheaf ${\rm Gr}_{\bbf{G}, \leq \mu}$ (resp. ${\rm Gr}_{\bbf{G},\mu}$) of ${\rm Gr}_{\bbf{G},E}:={\rm Gr}_{\bbf{G}}\times_{\bb{Q}_p^{\lozenge}}E^{\lozenge}$ 
        with the property that a map $\tau\colon S\to {\rm Gr}_{\bbf{G},E}$ factors through ${\rm Gr}_{\bbf{G}, \leq \mu}$ (resp. ${\rm Gr}_{\bbf{G},\mu}$) if and only if 
        $x^*(\tau)\in \coprod_{\lambda\leq \mu} {\rm Gr}_{\bbf{G},\lambda}(C)\subset {\rm Gr}_{\bbf{G}}(C)$ (resp. $x^*(\tau)\in {\rm Gr}_{\bbf{G},\lambda}(C)$)
        for all geometric points $x={\rm Spa}(C,C^+)\to S$.

        \begin{theorem}\label{affgrass}
        
        \begin{enumerate}
          \item  ${\rm Gr}_{\bbf{G}, \leq \mu}$ is  closed subsheaf of ${\rm Gr}_{\bbf{G},E}$ and the natural map ${\rm Gr}_{\bbf{G}, \leq \mu}\to E^{\lozenge}$ is a proper map of spatial diamonds.
          
          \item  ${\rm Gr}_{\bbf{G},\mu}$ is an open subsheaf of ${\rm Gr}_{\bbf{G}, \leq \mu}$ and there is a natural Bialynicki-Birula morphism of locally spatial diamonds over $\Spd E$, which is an isomorphism if and only if $\mu$ is minuscule
            $${\rm BB}_{\mu}: {\rm Gr}_{\bbf{G}, \mu}\to \mathcal{F\ell}_{\mu}^{\lozenge}.$$
            
            \end{enumerate}
        \end{theorem}
        
        \begin{proof}
        
     See    \cite[Proposition 19.2.3 and  Theorem 19.2.4]{ScholzeWeinspadicgeometry} for (1), the construction of (2) and the fact that it is an isomorphism if $\mu$ is minuscule. If $\mu$ is not minuscule then $\mathrm{BB}_{\mu}$ is not an isomorphism by  \cite[Theorem 5.4 (4)]{HoweKlevdalII}.
         \end{proof}
         
          Let us recall the construction of ${\rm BB}_{\mu}$, since it will be needed later on. For simplicity assume that there is $\mu\in \{\mu\}$ defined over $E$ (else replace 
          $E$ by a finite Galois extension over which such representative exists, and use Galois descent). 
          By \cite[Proposition  VI.2.4]{FarguesScholze} the map ${\rm Gr}_{\bbf{G}, \mu}\to 
          E^{\lozenge}$ has a natural section $[\mu]$ (obtained by sending an affinoid perfectoid $S/E^{\lozenge}$ to the element of $\xi^{\mu}$, where 
          $\xi$ is any generator of $\theta: \bb{B}^+_{\dR}(S^{\sharp})\to \s{O}(S^{\sharp})$), the 
           natural left action of $L^+\bbf{G}$ on ${\rm Gr}_{\bbf{G}, \mu}$ is transitive  by definition, and the stabilizer $(L^+\bbf{G})_{\mu}$ of $\xi^{\mu}$ is sent by the natural map 
           $\theta: L^+\bbf{G}\to \bbf{G}^{\lozenge}$ onto $\bbf{P}_{\mu}^{\lozenge}$, that is, the stabilizer $(L^+\bbf{G})_{\mu}$ consists on those elements $g\in L^+\bbf{G}$ such that $\mu(\xi)^{-1}g\mu(\xi)\in L^+\bbf{G}$, and so the image $\overline{g}$ of $g$ into $\bbf{G}^{\lozenge}$ must be in the locus where $\lim_{t\to 0} \mu(t)^{-1} g \mu(t)$ converges which is precisely $\bbf{P}_{\mu}$.           Therefore ${\rm BB}_{\mu}$ is simply the map
           $${\rm Gr}_{\bbf{G}, \mu}\simeq  (L^+\bbf{G})_{\mu}\backslash L^+\bbf{G}\to  \bbf{P}_{\mu}^{\lozenge} \backslash \bbf{G}^{\lozenge}\simeq \mathcal{F\ell}_{\mu}^{\lozenge}.$$

\subsection{Geometry of ${\rm Bun}_{\bbf{G}}$} We briefly recall the main properties of the stack ${\rm Bun}_{\bbf{G}}$ of $\bbf{G}$-bundles on the Fargues-Fontaine curve.
 We keep the notations of the previous sections.     

\subsubsection{$\bbf{G}$-isocrystals}

 To any $b\in \bbf{G}(\breve{\bb{Q}}_p)$ we can attach the following objects:
 
\begin{itemize} 
 
\item  A $\bbf{G}$-isocrystal\footnote{Recall that this is an additive, exact tensor functor from ${\rm Rep}(\bbf{G})$ to the category of 
 isocrystals, i.e. finite dimensional $\breve{\bb{Q}}_p$-vector spaces with a $\sigma$-linear automorphism.}
  $D_b\in {\rm Isoc}_{\bbf{G}}$, sending $(V, \rho)\in {\rm Rep}(\bbf{G})$ to 
$$D_b((V,\rho))=(V\otimes_{\bb{Q}_p} \breve{\bb{Q}}_p, \rho(b)({\rm id}_V\otimes \sigma))$$
  for 
   $(V, \rho)\in {\rm Rep}(\bbf{G})$. 
  
\item  A reductive group $\bbf{G}_b$ over $\bb{Q}_p$ with $\bbf{G}_b(\bb{Q}_p)={\rm Aut}(D_b)$. More precisely, for any 
   $\bb{Q}_p$-algebra $R$ the group $\bbf{G}_b(R)$ is the subgroup of $g\in \bbf{G}(\breve{\bb{Q}}_p\otimes_{\bb{Q}_p} R)$ for which 
   $gb=b\sigma(g)$.   
   
\item  A Newton morphism 
        $\nu_b: \bb{D}_{\breve{\bb{Q}}_p}\to \bbf{G}_{\breve{\bb{Q}}_p}$, encoding the Dieudonn\'e-Manin decomposition of the various 
   $D_b((V, \rho))$. Here 
$\bb{D}$ is the slope pro-torus with character group $\bb{Q}$ 
and we refer to \cite[Theorem 1.8]{RapoportRichartz} for Kottwitz's characterization of 
$\nu_b$. The morphism $\nu_b$ factors through the center of $\bbf{G}_{\breve{\bb{Q}}_p}$ if and only if 
$\bbf{G}_b$ is an inner form of $\bbf{G}$ (see \cite[Proposition  1.12]{RapoportRichartz}), in which case we say that 
$b$ is \emph{basic}.

\end{itemize}

    Any $\bbf{G}$-isocrystal is isomorphic to $D_b$ for some $b\in \bbf{G}(\breve{\bb{Q}}_p)$, and 
   $D_{b_1}\simeq D_{b_2}$ if and only 
   $b_1,b_2$ are $\sigma$-conjugate, i.e. $b_2=gb_1\sigma(g)^{-1}$ for some 
    $g\in \bbf{G}(\breve{\bb{Q}}_p)$. The Kottwitz set 
    $B(\bbf{G})$ of $\sigma$-conjugacy classes $[b]$ of elements $b\in \bbf{G}(\breve{\bb{Q}}_p)$
   classifies therefore isomorphism classes of $\bbf{G}$-isocrystals. Moreover, we have 
$\nu_{gb\sigma(g)^{-1}}=g\nu_b g^{-1}$ and $\nu_b=b\sigma(\nu_b)b^{-1}$ for 
   all $g,b\in \bbf{G}_{\breve{\bb{Q}}_p}$, thus $b\mapsto \nu_b$ induces a map 
   $$\nu_{\bbf{G}}: B(\bbf{G})\to \mathcal{N}(\bbf{G}):=({\rm Hom}(\bb{D}_{\breve{\bb{Q}}_p}, \bbf{G}_{\breve{\bb{Q}}_p})/\bbf{G}(\breve{\bb{Q}}_p)-\text{conj})^{\sigma}.$$ 
    
Let $\pi_1(\bbf{G})$ be the Borovoi fundamental group, a discrete 
     $\Gamma$-module, quotient
      of $X_*(\bbf{T})$ by the coroot lattice\footnote{Note that $\Gamma$ acts on $X_*(\bbf{T})$ by $\gamma.\chi=g\gamma(\chi)g^{-1}$, where 
     $g\in \bbf{G}(\overline{\bb{Q}}_p)$ is such that $g(\gamma(\bbf{T},\bbf{B}))g^{-1}=(\bbf{T},\bbf{B})$. The resulting action of 
     $\Gamma$ on $\pi_1(\bbf{G})$ is independent of the choice of 
     $\bbf{B}$ and $\bbf{T}$, up to canonical isomorphism.}. Kottwitz proved (see \cite[Theorem 1.15]{RapoportRichartz}) the existence of a unique functorial map 
   (in $G$) $\kappa: B(\bbf{G})\to \pi_1(\bbf{G})_{\Gamma}$ which 
       is induced by the $p$-adic valuation map $v_p: \bbf{G}(\bb{Q}_p)\to \bb{Z}=\pi_1(\bbf{G})_{\Gamma}$ for $\bbf{G}=\bb{G}_m$. The map $\kappa$ induces a bijection between $B(\bbf{G})_{\rm basic}$ and $\pi_1(\bbf{G})_{\Gamma}$ and the map
       $(\nu, \kappa): B(\bbf{G})\to \mathcal{N}(\bbf{G})\times \pi_1(\bbf{G})_{\Gamma}$ is injective. 
   
   If $W$ is the Weyl group of $(\bbf{G}_{\overline{\bb{Q}}_p}, T)$, there is a natural bijection
   $\mathcal{N}(\bbf{G})\simeq (X_*(T)_{\bb{Q}}/W)^{\Gamma}$ and  
  the set $X_*(T)_{\bb{R}}/W$ has a natural partial order, compatible with the Bruhat order on $X_*(T)/W\simeq X_*(T)^+$, 
   where $\Omega_1\leq \Omega_2$ if $\Omega_1$ is contained in the convex hull of $\Omega_2$ (see \cite[Lemma 2.2]{RapoportRichartz}). 
       This induces a partial order on $B(\bbf{G})$, defined by $x_1\leq x_2$ if and only if 
       $\kappa(x_1)=\kappa(x_2)$ and $\nu_b(x_1)\leq \nu_b(x_2)$.

\subsubsection{Topology and stratification of ${\rm Bun}_{\bbf{G}}$}

Let ${\rm Bun}_{\bbf{G}}$ be the small $v$-stack sending $S\in \Perf_k$ to 
the grupoid of $\bbf{G}$-bundles on $X_S$ (see \cite[Proposition  III.1.3]{FarguesScholze}). 
Any $b\in \bbf{G}(\breve{\bb{Q}}_p)$ 
gives rise to a functor $S\in \Perf_k\mapsto \mathcal{E}_{b, S}\in {\rm Bun}_{\bbf{G}}(S)$ corresponding to 
a map 
   $$\mathcal{E}_b: *:=k^{\lozenge}\to {\rm Bun}_G.$$
   The $\bbf{G}$-torsor $\mathcal{E}_{b, S}$ on $X_{S}$ is obtained via descent from the trivial torsor $\bbf{G}\times Y_{S}$ with Frobenius $b\times \varphi$. 
Varying $b$, this induces a map $B(\bbf{G})\to |{\rm Bun}_{\bbf{G}}|$, and Fargues \cite{FargerTorseursCourbe}
proved that it is a bijection. Based on this and previous work of Fargues and Scholze \cite[Section  III.3]{FarguesScholze}, Viehmann \cite{ViehmannStrata} proved that this map is a homeomorphism, when we endow $B(\bbf{G})$ with the order topology induced by the \emph{opposite} of the partial order on $B(\bbf{G})$ described in the previous subsection.
    More concretely, 
the closure in $|{\rm Bun}_{\bbf{G}}|$ of a subset $S$ is $\bigcup_{\mathcal{E}\in S} \overline{\{\mathcal{E}\}}$ and for any 
$b\in \bbf{G}(\breve{\bb{Q}}_p)$ the closure of (the point corresponding to) $[b]$ consists of those $[b']$ with 
$[b]\leq [b']$ \cite[Theorem 6.7 and Corollary 6.8]{ViehmannStrata}.

  The automorphism group of $\n{E}_{b}$ is the $v$-sheaf of groups 
\[
\widetilde{G}_b= \underline{\Aut}_{\bbf{G}}(\n{E}_b): (S\in \Perf_k)\mapsto \Aut_{\bbf{G}\times X_{S}} (\n{E}_{b,S}).
\]
Let 
$\widetilde{G}_b^{>0}$ be the subgroup of $\widetilde{G}_b$ consisting of unipotent automorphisms with respect to the Harder-Narasimhan filtration of $\n{E}_{b,S}$. Fargues and Scholze proved in \cite[Proposition III.5.1]{FarguesScholze} that 
there is a canonical split exact sequence 
   $$0\to \tilde{G}_b^{>0}\to \tilde{G}_b\to \underline{\bbf{G}_b(\bb{Q}_p)}\to 0,$$
  the section $\underline{\bbf{G}_b(\bb{Q}_p)}\to  \tilde{G}_b$ being provided by the interpretation of 
  $\bbf{G}_b(\bb{Q}_p)$ as automorphism group of the isocrystal $D_b$ and the fact that 
  $H^0(X_S, \mathcal{O})=\underline{\bb{Q}_p}(S)$.
   Moreover $\tilde{G}_b^{>0}$ is a unipotent group diamond, successive extension of positive Banach-Colmez spaces, of dimension 
   $\langle \nu([b]), 2\rho\rangle$ (where $2\rho$ is the sum of positive roots). When $b$ is basic 
   the group $\tilde{G}_b^{>0}$ is trivial and so $ \tilde{G}_b=\underline{\bbf{G}_b(\bb{Q}_p)}$ in this case.

  For any $[b]\in B(\bbf{G})$ the 
 stratum
   $${\rm Bun}_{\bbf{G}}^b:={\rm Bun}_{\bbf{G}}\times_{|{\rm Bun}_{\bbf{G}}|} \{[\mathcal{E}_b]\}$$
  classifies those $\bbf{G}$-bundles on $X_S$ which are $v$-locally isomorphic to 
  $\mathcal{E}_b$. It is a locally closed substack of ${\rm Bun}_{\bbf{G}}$ and the map 
  $\mathcal{E}_b: *\to {\rm Bun}_{\bbf{G}}$ induces an isomorphism 
   $$[*/ \tilde{G}_b]\simeq {\rm Bun}_G^b$$
   between the classifying stack of $\tilde{G}_b$-torsors and ${\rm Bun}_G^b$. The substack 
   ${\rm Bun}_G^b$ is open in ${\rm Bun}_G$ when $b$ is basic.

\subsection{Moduli spaces of modifications}

Recall the following standard definition:

\begin{definition}
\label{DefinitioLocalShimuraDatum}
A local shtuka datum is a triple $(\bbf{G},[b],\{\mu\})$ consisting of a reductive group $\bbf{G}$ over $\bb{Q}_p$, a conjugacy class $\{\mu\}$ of cocharacters $\bb{G}_m\to \bbf{G}_{\overline{\bb{Q}}_p}$, and an element $[b]\in B(\bbf{G})$. If moreover $\{\mu\}$ is minuscule, we suppose that $[b]\in B(\bbf{G},\{\mu^{-1}\})$ and say that $(\bbf{G},[b],\{\mu\})$ is a local Shimura datum.\footnote{For $\mu$ minuscule, the local Shimura variety associated to  $(\bbf{G},[b],\{\mu\})$ is non-empty if and only if $[b]\in B(\bbf{G},\{\mu^{-1}\})$, see \cite[Proposition 24.1.2]{ScholzeWeinspadicgeometry}.}
\end{definition}

    We fix from now on a local shtuka datum $(\bbf{G},[b],\{\mu\})$ and we freely use the notations introduced in the previous sections. 

    Let $\mathcal{E}, \mathcal{E}'$ be $\bbf{G}$-bundles on a sous-perfectoid space $X$ with a closed Cartier divisor $D$. Recall that a \emph{modification at $D$} $\alpha: \mathcal{E}\dashrightarrow \mathcal{E}'$ is an isomorphism of $\bbf{G}$-torsors 
     $$\alpha: \mathcal{E}|_{X\setminus D}\simeq \mathcal{E}'|_{X\setminus D}$$ which is meromorphic along $D$.

      Let $S$ be a perfectoid space over $\breve{E}^{\lozenge}$ with associated untilt $S^{\sharp}$.    The map  $D=S^{\sharp}\to Y_{S}$ defines a degree $1$ closed Cartier divisor on the curve $Y_{S}$ and hence on the Fargues--Fontaine curve   $X_{S}$. Let $\mathcal{E}$ and $\mathcal{E}'$  be two $\mathbf{G}$-bundles on $X_{S}$ and let $\alpha: \mathcal{E}\dashrightarrow \mathcal{E}'$ be a modification. Suppose first that $S={\rm Spa}(C, C^+)$ is a geometric point. The completion $\widehat{\mathcal{E}}'$ along $D$ is isomorphic to a trivial $\bbf{G}$-torsor over ${\rm Spec}(\bb{B}_{\dR}(C^{\sharp}))$. 
   Fixing  a trivialisation  $\widehat{\mathcal{E}}_1\cong \widehat{\mathcal{E}}'$, we have a point $x_{\alpha}\in {\rm Gr}_{\bbf{G}}(C)=\mathbf{G}(\mathbb{B}_{\dR}(C))/\mathbf{G}(\mathbb{B}^+_{\dR}(C))$ corresponding to the
   $\bbf{G}$-torsor $\widehat{\mathcal{E}}$ over ${\rm Spec}(\bb{B}^+_{\dR}(C^{\sharp}))$ endowed with the trivialisation over ${\rm Spec}(\bb{B}_{\dR}(C^{\sharp}))$ induced by $\alpha$ and the chosen trivialisation of 
   $\widehat{\mathcal{E}}'$. We say that $\alpha$ is \emph{bounded by $\{\mu\}$}
   if $x_{\alpha}\in {\rm Gr}_{\bbf{G}, \leq \mu}(C)$. More explicitly, after fixing a trivialization of both $\widehat{\mathcal{E}}$ and $\widehat{\mathcal{E}}'$, the modification $\alpha$ is precisely the datum of an element $x_{\alpha}\in \mathbf{G}(\mathbb{B}_{\dR}(C))$, forgetting the trivialization of $\widehat{\mathcal{E}}$ produces then an element in the right quotient $\mathbf{G}(\mathbb{B}_{\dR}(C))/\mathbf{G}(\mathbb{B}^+_{\dR}(C))=\mathrm{Gr}_{\mathbf{C}}(C)$, if instead we only forget the trivialization of $\widehat{\mathcal{E}}'$ we get an element in the left quotient $\mathbf{G}(\mathbb{B}^+_{\dR}(C))\backslash \mathbf{G}(\mathbb{B}_{\dR}(C))$.  For general $S$, we say that 
   $\alpha$ is \emph{bounded by $\{\mu\}$} if for any geometric point $x: {\rm Spa}(C, C^+)\to S$ the induced modification $x^*(\alpha)$ between $\bbf{C}$-torsors on $X_{{\rm Spa}(C,C^+)}$ is bounded by $\mu$.

\begin{definition}
The space $\n{M}_{\infty}:=\n{M}_{\bbf{G},b,\mu,\infty}$ of modifications at infinite level associated to the triple $(\bbf{G}, [b], \{\mu\})$ is the sheaf on 
$\Perf_k/\breve{E}^{\lozenge}$ sending 
$S/\breve{E}^{\lozenge}$ to the set of isomorphism classes of modifications 
$\alpha: \n{E}_1\dashrightarrow \n{E}_b$ which are bounded by $\mu$. 

\end{definition}

\begin{remark}[Period maps]\label{RemarkConstructionsModifications}
  By construction, the $v$-sheaf $\n{M}_{\infty}$ has the following properties:

  \begin{enumerate}
  \item It is equipped with commuting right and left actions of $\underline{G(\bb{Q}_p)}$ and of 
  $\tilde{G}_b$ respectively, where $g\in \underline{G(\bb{Q}_p)}(S)=\tilde{G}_1$ acts by $\alpha\mapsto \alpha\circ g$, while $h\in \tilde{G}_b$ acts by $\alpha\mapsto h \circ \alpha$. 
  
 \item  There is a canonical map of $v$-sheaves over $\breve{E}^{\lozenge}$
  $$g_{\rm dR}: \n{M}_{\infty}\to \mathrm{Gr}_{\mathbf{G}}$$
induced by the formal completion of the modification $\alpha$ at the divisor $S^{\sharp}\to X_{S}$ given by the point $S\to \breve{E}^{\diamond}$, and the natural trivialization of the torsor $\n{E}_b$. By definition, this map even factors as a map 
\[
\pi_{\GM}\colon  \mathcal{M}_{\infty} \to \mathrm{Gr}_{\mathbf{G},\leq \mu, \breve{E}} 
\]
called the \textit{Grothendieck-Messing period map}.  This map is $\widetilde{G}_b$-equivariant, and factors through the quotient 
\[
\mathcal{M}_{\infty}\to \mathcal{M}_{\infty}/\underline{G(\bb{Q}_p)}\to \mathrm{Gr}_{\mathbf{G},\leq \mu,\breve{E}}. 
\]

   \item If we swap the roles of $\mathcal{E}_1$ and $\mathcal{E}_b$ and use the trivialization of $\mathcal{E}_1$, the modification $\alpha^{-1}\colon\mathcal{E}_b\dashrightarrow \mathcal{E}_1$ is bounded by $\mu^{-1}$ and produces a $\underline{G(\mathbb{Q}_p)}$-equivariant \textit{Hodge--Tate period map} 
   \[
   \pi_{\HT}\colon \mathcal{M}_{\infty} \to \widetilde{G}_{b}\backslash \mathcal{M}_{\infty}\to \mathrm{Gr}_{\mathbf{G},\leq \mu^{-1},\breve{E}}. 
   \] 
   In total, we get a $\underline{G(\bb{Q}_p)}\times \widetilde{G}_b$-equivariant correspondence 
   \begin{equation}\label{eqCorrespondenceGManddHT}
   \begin{tikzcd}
    & \mathcal{M}_{\infty} \ar[ld,"\pi_{\GM}"'] \ar[rd,"\pi_{\HT}"] & \\
 \mathrm{Gr}_{\mathbf{G},\leq \mu,\breve{E}}    & &  \mathrm{Gr}_{\mathbf{G},\leq \mu^{-1},\breve{E}}
   \end{tikzcd}
   \end{equation}
      
   \end{enumerate}
   \end{remark}

 Let $b\in B(\mathbf{G})$ and $S\in \Perf_k$. Given an $S$-point $x_{\alpha}$ of $\mathrm{Gr}_{\mathbf{G}}(S)$ we can use the Beauville-Laszlo gluing lemma  of vector bundles to construct  a modification of $\mathcal{E}_b$ in $X_{S}$. This functor produces a  \textit{Beauville-Laszlo map} of $v$-stacks  $\mathrm{BL}_b\colon \mathrm{Gr}_{\mathbf{G}}\to \mathrm{Bun}_{\mathbf{G}}$.  Following \cite[Definition 3.1]{ViehmannStrata}, for $b,b'\in B(|\mathbf{G})$ we denote the Newton stratum  $$\Gr_{\bbf{G},b}^{[b']}:={\rm BL}_b^{-1}({\rm Bun}^{b'}_{\bbf{G}}).$$ We also let $\mathrm{Gr}_{\bbf{G},b,\leq \mu,\breve{E}}^{[b']}$ be the $\mu$-bounded locus over $\breve{E}^{\lozenge}$.

Thus, we have the following cartesian diagrams, consequence of the Beauville-Laszlo gluing:
 
 \begin{equation}
\label{cart diag1}
   \begin{tikzcd}
 \n{M}_{\infty} \ar[r] \ar[d, "\pi_{\GM}"']& k^{\lozenge} \arrow[d] \\ 
\Gr_{\bbf{G},b,\leq \mu,\breve{E}}^{[1]} \arrow[r, "{\rm BL}_b"']  & {\rm Bun}_{\bbf{G}}^1 
\end{tikzcd},
\end{equation}
and 
 \begin{equation}
\label{cart diag2}
 \begin{tikzcd}
 \n{M}_{\infty} \ar[r] \ar[d, "\pi_{\HT}"']& k^{\lozenge} \arrow[d] \\ 
\Gr_{\bbf{G},1,\leq \mu^{-1},\breve{E}}^{[b]} \arrow[r, "{\rm BL}_1"']  & {\rm Bun}_{\bbf{G}}^b
\end{tikzcd}
\end{equation}

 In particular we obtain the well-known result:
 
 \begin{proposition}\label{class}
The map $$\pi_{\GM}:  \n{M}_{\infty} \to \Gr_{\bbf{G},b,\leq \mu,\breve{E}}^{[1]}$$ is a $\tilde{G}_b$-equivariant
 pro\'etale $\underline{\bbf{G}(\bb{Q}_p)}$-torsor and the map
 $$\pi_{\HT}: \n{M}_{\infty} \to\Gr_{\bbf{G},1,\leq \mu^{-1},\breve{E}}^{[b]}$$ is a $\underline{\bbf{G}(\bb{Q}_p)}$-equivariant pro\'etale 
 $\tilde{G}_b$-torsor.
 \end{proposition}

When $b$ is basic we have a duality for the diagram \eqref{eqCorrespondenceGManddHT}, cf. \cite[Corollary 23.3.2]{ScholzeWeinspadicgeometry}.

\begin{proposition}
\label{PropositionDualityLocalShtukas}
Suppose that $b$ is basic and set $\check{\bbf{G}}=\bbf{G}_b$, $\check{b}=b^{-1}\in \bbf{G}_b(\breve{\bb{Q}}_p)=\bbf{G}(\breve{\bb{Q}}_p)$ and $\check{\mu}=\mu^{-1}$ under the identification $\bbf{G}_{\overline{\bb{Q}}_p} \cong \check{\bbf{G}}_{\overline{\bb{Q}}_p}$. Then there is a natural $\bbf{G}(\bb{Q}_p)\times \check{\bbf{G}}(\bb{Q}_p)$-equivariant isomorphism
\begin{equation}
\label{eqisodualitymoduliShtukkas}
\n{M}_{\bbf{G},b,\mu,\infty}\cong \n{M}_{\check{\bbf{G}},\check{b},\check{\mu},\infty}
\end{equation}
interchanging the maps $\pi_{\GM}$ and $\pi_{\HT}$ of \eqref{eqCorrespondenceGManddHT}. 
\end{proposition}
\begin{proof}
By \cite[Corollary III.4.3]{FarguesScholze} there is an isomorphism of $v$-stacks 
$$\iota: {\rm Bun}_{\bbf{G}}\simeq {\rm Bun}_{\check{\bbf{G}}},\,\, \iota(\mathcal{E})={\rm Isom}_{\bbf{G}}(\mathcal{E}_b, \mathcal{E}),$$
such that $\iota(\mathcal{E}_{b'})=\mathcal{E}_{b'b^{-1}}$ for all 
$b'$, in particular $\iota(\mathcal{E}_1)=\mathcal{E}_{\check{b}}$ and $\iota(\mathcal{E}_ b)=\mathcal{E}_1$. Applying 
$\iota$ to a modification $\alpha: \mathcal{E}_1\dashrightarrow \n{E}_b$ yields a modification 
$\iota(\alpha): \mathcal{E}_{\check{b}}\dashrightarrow \n{E}_1$ of $\check{\bbf{G}}$-bundles and by definition we have 
$g_{\rm dR}(\iota(\alpha))=g_{\rm dR}(\alpha)^{-1}$ in $L\bbf{G}_{\breve{E}}\simeq L\bbf{\check{G}}_{\breve{E}}$. The result follows. 
\end{proof}

\section{Geometric Sen operators and the Hodge-Tate period map of local Shimura varieties}
\label{Section:SenOperators}

In this section we compute the geometric Sen operator of local Shimura varieties, proving the local analogue of \cite[Theorem 5.2.5]{RCLocAnCompleted}. We start  in \cref{ss:DeRhamLocalSystems} with  the Riemann-Hilbert functor of \cite{ScholzeHodgeTheory2013,LZ},  the notion of de Rham local systems and the Hodge--Tate filtration.  In  \Cref{ss:SetUpNotationShimuraGeoRep} we give the appropriate  set up of the relevant period maps and notation from geometric representation theory.    In \Cref{ss:RiemannHilbertLocalShimura} we discuss the  Riemann--Hilbert correspondence for local Shimura varieties; in contrast with the global situation this one is almost a tautology after a bookkeeping of the constructions. In \Cref{ss:PullbackHTSheavesEquivariant} we compute the pullbacks of equivariant vector bundles along the Hodge--Tate period map. Finally, in \Cref{SubsectionSenOperators} we compute the geometric Sen operators of local Shimura varieties.

\subsection{De Rham local systems and the Hodge-Tate filtration}\label{ss:DeRhamLocalSystems}
   
     Let $K$ be a complete discretely valued extension of $\bb{Q}_p$ with perfect residue field and let 
     $X$ be a smooth rigid analytic variety over $K$. We can associate to $X$ two categories:
     
\begin{itemize}     
     
     \item The category 
     ${\rm Loc}_{\bb{Q}_p}(X)$ of 
     $\bb{Q}_p$-local systems on $X$, i.e. sheaves of $\underline{\bb{Q}}_p$-modules on $X_{\proet}$ that are locally free of finite rank.

    \item  The category ${\rm FIC}_X$ of filtered integrable connections on $X$, i.e. vector bundles 
     $\mathcal{E}$ on 
     $X_{\rm an}$ endowed with a separated, exhaustive, decreasing filtration 
     $F^{\bullet}\mathcal{E}$ by sub $\s{O}_X$-modules that are locally on $X_{\rm an}$ direct summands, and an integrable connection $\nabla: \mathcal{E}\to \mathcal{E}\otimes_{\s{O}_X}\Omega^1_X$, satisfying Griffiths transversality, i.e. $\nabla(F^i  \mathcal{E})\subset F^{i-1}\mathcal{E}\otimes_{\mathcal{O}_X}\Omega^1_X$. 
     
     \end{itemize}
          
     Recall that on $X_{\proet}$ there is the big period sheaf $\OBdr$, an $\s{O}_X$-module endowed with an integrable connection and a 
     decreasing exhaustive separated filtration satisfying Griffiths transversality.\footnote{Strictly speaking, one needs to work with the pro\'etale site of \cite{ScholzeHodgeTheory2013} for defining the sheaf $\OBdr$ and realizing the Riemann-Hilbert correspondence. However, the notion of   $\mathbb{Q}_p$-pro\'etale system is independent of the version of pro\'etale site we choose to work with. } For any $\bbf{L}\in {\rm Loc}_{\bb{Q}_p}(X)$ define 
     $$D_{\rm dR}(\bbf{L}):=\lambda_*(\bbf{L}\otimes_{\underline{\bb{Q}}_p} \OBdr),$$
     where $\lambda: X_{\proet}\to X_{\et}$ is the natural projection of sites. This is a vector bundle on 
     $X_{\rm an}$ by \cite[Theorem 3.8]{LZ}, endowed with an integrable connection and a decreasing (Hodge) filtration, satisfying Griffiths transversality, induced by those on $\OBdr$.\footnote{The results of \cite{LZ} are stated for \'etale $\mathbb{Q}_p$ local systems, i.e. pro\'etale $\mathbb{Q}_p$-local systems that are the generic fiber of a pro\'etale $\mathbb{Z}_p$-local system. It is easy to show that any pro\'etale $\mathbb{Q}_p$-local system admits locally in the analytic topology of $X$ a $\mathbb{Z}_p$-lattice, therefore the conclusion of \cite[Theorem 3.8]{LZ} also hold for them.} If $\mathcal{E}$ is a vector bundle on 
    $X_{\rm an}$ (or equivalently on $X_{\et}$), we will simply write $\mathcal{E}$ for $\lambda^{-1}\mathcal{E}$, a  $\s{O}_X$-vector bundle on $X_{\proet}$. 
     By construction, there is a natural comparison map, compatible with filtrations and connection
     $${\rm comp}_{\bbf{L}}: D_{\rm dR}(\bbf{L})\otimes_{\s{O}_X} \OBdr\to \bbf{L}\otimes_{\underline{\bb{Q}}_p} \OBdr,$$
     and we let ${\rm Loc}^{\rm dR}_{\bb{Q}_p}(X)$ be the full subcategory of de Rham local systems, i.e. those for which the comparison map is an isomorphism.\footnote{The reader can easily check that this definition of de Rham local systems is the same as the one by Scholze in \cite{Scholze2013moduli}, using that
     for any vector bundle $\mathcal{E}$ on $X$ we have a canonical isomorphism 
     $\mathcal{E}\simeq \lambda_*(\mathcal{E}\otimes_{\mathcal{O}_X} \OBdr)$.}
      Theorem 3.8 of  \cite{LZ} ensures that the functor 
     $$D_{\rm dR}: {\rm Loc}^{\rm dR}_{\bb{Q}_p}(X)\to {\rm FIC}_X$$
     is well-defined and is a tensor-functor. It is also easy to see that it is exact.
     
     \begin{proposition}\label{resumeSch}
      Let $\bbf{L}\in {\rm Loc}^{\rm dR}_{\bb{Q}_p}(X)$ and set $\mathcal{E}=D_{\rm dR}(\bbf{L})$ and 
        $$M=\bbf{L}\otimes_{\underline{\bb{Q}}_p} \bb{B}_{\dR}^+, \,\, M_0=(\mathcal{E}\otimes_{\mathcal{O}_X} 
     \OBdr^+)^{\nabla=0}.$$
     Identify $M_0$ with a sub $\bb{B}_{\dR}^+$-module of $\bbf{L}\otimes_{\underline{\bb{Q}}_p} \bb{B}_{\dR}$ via the isomorphism induced by 
     ${\rm comp}_{\bbf{L}}$
    $$(\mathcal{E}\otimes_{\mathcal{O}_X} 
     \OBdr)^{\nabla=0}\simeq \bbf{L}\otimes_{\underline{\bb{Q}}_p} \bb{B}_{\dR}.$$
      
\begin{enumerate}      
      
      \item  The $\bb{B}_{\dR}^+$-modules 
     $M$ and $M_0$ are finite free locally on $X_{\proet}$
     and the inclusion $M_0\subset \bbf{L}\otimes_{\underline{\bb{Q}}_p} \bb{B}_{\dR}$ described above induces an isomorphism
     $$\bb{B}_{\dR}\otimes_{\bb{B}_{\dR}^+} M_0\simeq \bb{B}_{\dR}\otimes_{\underline{\bb{Q}}_p}  \bbf{L}.$$
     
     \item  There are natural isomorphisms of $ \widehat{\s{O}}_X$-modules 
     (for all $i\in \bb{Z}$)
     $$\frac{M\cap {\rm F}^i(M_0)}{M\cap {\rm F}^{i+1}(M_0)}\simeq {\rm F}^{-i}(\mathcal{E})\otimes_{\s{O}_X}
     \widehat{\s{O}}_X(i), \,\, \frac{{\rm F}^i(M_0)}{{\rm F}^{i+1}(M_0)}\simeq \mathcal{E}\otimes_{\s{O}_X}
     \widehat{\s{O}}_X(i).$$
In the previous equation, if $N$ is a finite free $ \bb{B}_{\dR}^+$-module on $X_{\proet}$, we endow $N$ with the filtration ${\rm F}^i(N)=(\ker \theta)^i N$, where 
     $\theta: \bb{B}_{\dR}^+\to \widehat{\s{O}}_X$ is the natural quotient map.

     \item The increasing relative Hodge-Tate filtration 
        $${\rm F}_{\rm HT, i}(\bbf{L}\otimes_{\underline{\bb{Q}}_p}\widehat{\s{O}}_X):=\frac{M\cap {\rm F}^{-i}(M_0)}{{\rm F}^1(M)\cap {\rm F}^{-i}(M_0)}$$
    on $\bbf{L}\otimes_{\underline{\bb{Q}}_p}\widehat{\s{O}}_X=M/{\rm F}^1(M)$ comes with a canonical isomorphism
     $${\rm gr}_{\rm HT, i}(\bbf{L}\otimes_{\underline{\bb{Q}}_p}\widehat{\s{O}}_X)\simeq {\rm gr}^{i}(\mathcal{E})\otimes_{\s{O}_X} 
     \widehat{\s{O}}_X(-i).$$
     
     \end{enumerate}
     
     \end{proposition}
     
     \begin{proof}
      Parts (1) and (2) follow from \cite[Theorem 7.6 and  Proposition 7.9]{Scholze2013moduli}. Part (3) is a direct consequence of part (2), once we observe that 
      $$\frac{{\rm F}^1(M)\cap {\rm F}^i(M_0)}{{\rm F}^1(M)\cap {\rm F}^{i+1}(M_0)}\simeq  \frac{M\cap {\rm F}^{i-1}(M_0)}{M\cap {\rm F}^{i}(M_0)}(1).$$     \end{proof}

      The category  ${\rm Loc}^{\rm dR}_{\bb{Q}_p}(X)$ is stable under subquotients and tensor products,  and for $\bbf{G}$ an affine algebraic group over $\bb{Q}_;$ we define the category of de Rham $\bbf{G}$-local systems
     ${\rm Loc}^{\rm dR}_{\bbf{G}}(X)$ using the Tannakian definition of $\bbf{G}$-torsors, that is a de Rham $\bbf{G}$-local system is an additive and exact tensor functor ${\rm Rep}(\bbf{G})\to {\rm Loc}^{\rm dR}_{\bb{Q}_p}(X)$.  The category 
     ${\rm Loc}_{\bbf{G}}(X)$ of 
     $\bbf{G}$-local systems on $X$ 
          consists of additive and exact tensor functors
     ${\rm Rep}(\bbf{G})\to {\rm Loc}_{\bb{Q}_p}(X)$, and 
     the category 
     ${\rm FIC}_{\bbf{G}, X}$ of filtered integral $\bbf{G}$-connections on $X$ is defined similarly, using the Tannakian formalism and 
     the category ${\rm FIC}_X$.  Since the functor $D_{\rm dR}$ is a tensor functor when restricted to de Rham local systems, it gives rise to a functor of $\bbf{G}$-local systems 
     \[
     D_{\rm dR}\colon {\rm Loc}^{\rm dR}_{\bbf{G}}(X)\to  {\rm FIC}_{\bbf{G}, X}. 
     \]

\subsection{The setup}\label{ss:SetUpNotationShimuraGeoRep}

We fix a local Shimura datum $(\bbf{G},[b],\{\mu\})$ (thus $\{\mu\}$ is minuscule and $b\in B(G,\{\mu^{-1}\})$) with reflex field $E/\bb{Q}_p$, as well as a compact open subgroup $K$ of $\bbf{G}(\bb{Q}_p)$.
We denote for simplicity 
$$\n{M}_{\infty}:=\n{M}_{\bbf{G},b,\{\mu\},\infty}, \quad \n{M}_{K}=\n{M}_{\infty}/\underline{K},$$
living over 
$\breve{E}^{\lozenge}$. The Bialynicki-Birula maps 
\[
{\rm BB}_{\mu}: \Gr_{\bbf{G},\leq \mu,\breve{E}} = \Gr_{\bbf{G},\breve{E} ,\mu} \to \Fl_{\mu, \breve{E}}^{\lozenge},\quad {\rm BB}_{\mu^{-1}}: \Gr_{\bbf{G},\leq \mu^{-1},\breve{E}} = \Gr_{\bbf{G} ,\mu^{-1},\breve{E}} \to \Fl_{\mu^{-1}, \breve{E}}^{\lozenge} \]
are isomorphisms (\Cref{affgrass}), yielding a $\bbf{G}(\bb{Q}_p)\times \widetilde{G}_b$-equivariant diagram of period maps
\begin{equation}
\label{eqHodgeTatePeriodMapv2}
\begin{tikzcd}
 & \n{M}_{\infty} \ar[rd,"\pi_{\HT}"] \ar[ld, "\pi_{\GM}"']& \\ 
 \Fl_{\mu,\breve{E}}^{\lozenge}  & &  \Fl^{\lozenge}_{\mu^{-1},\breve{E}}. 
\end{tikzcd}
\end{equation}

   Let 
   $\Fl_{\mu,\breve{E}}^{a}$ be the open adic subspace of $\Fl_{\mu,\breve{E}}$ corresponding under the isomorphism ${\rm BB}_{\mu}$ to the open sub-diamond 
   $\Gr_{\bbf{G},b,\leq \mu,\breve{E}}^{[1]}$ of $\Gr_{\bbf{G},\leq \mu,\breve{E}}$. It follows from  \Cref{class} that the map 
   $$\pi_{\GM}:  \n{M}_{\infty}\to \Fl_{\mu,\breve{E}}^{a, \lozenge}$$
   is a $\tilde{G}_b$-equivariant $\underline{\bbf{G}(\bb{Q}_p)}$-torsor for the pro\'etale topology, hence the resulting map 

 \begin{equation}
 \label{eqEtaleGMMap}
\pi_{\GM}: \n{M}_{K}\to \Fl_{\mu,\breve{E}}^{\lozenge}
\end{equation}
 is \'etale,\footnote{Since the \'etale sites of $ \Fl_{\mu}^{\lozenge}$ and $ \Fl_{\mu}$ are equivalent, via the diamond functor.}  and there is a unique \'etale map $\n{M}_K^{\rm rig}=\n{M}_{\bbf{G},b,\mu,K}^{\rm rig} \to \Fl_{\mu, \breve{E}}$ whose diamond is $\pi_{\GM}$, and $\n{M}_K^{\rm rig}$ is a smooth rigid analytic variety over $\breve{E}$. For simplicity we denote the structural sheaf and sheaf of differential forms of $\mathcal{M}_K^{\rm rig}$ by 
 $$\s{O}_{\n{M}}=\s{O}_{\n{M}_K^{\rm rig}}, \quad \Omega^1_{\n{M}}=\Omega^1_{\n{M}_K^{\rm rig}},$$
and consider them as \'etale sheaves.

We fix  a finite extension  $E'/E$ which splits $\bbf{G}$ and for which we can find a representative 
 Hodge-cocharacter $\mu:\bb{G}_{m,E'}\to \bbf{G}_{E'}$ in $\{\mu\}$ defined over $E'$.\footnote{Taking this finite extension is unnecessary for the forthcoming discussion but it allows us to use the dictionary between representations of the chosen parabolic $\bbf{P}_{\mu}$ and $\bbf{G}$-equivariant quasi-coherent sheaves on the flag variety.   We leave to the reader the cocharacter-free formulation of the statements in terms of filtered $\bbf{G}$-representations.}   We will work systematically over 
 $\breve{E'}$, so we will base change all sheaves, spaces, etc. to $\breve{E'}$ or its diamond.  We let  $\bbf{P}_{\mu}$, $\bbf{P}_{\mu^{-1}}$ the associated  parabolic subgroups of $\bbf{G}_{E'}$, so that 
 $$\Fl_{\mu, \breve{E'}}=(\bbf{G}_{\breve{E'}}/ \bbf{P}_{\mu})^{\rm ad}, \quad \Fl_{\mu^{-1}, \breve{E'}}=(\bbf{G}_{\breve{E'}}/ \bbf{P}_{\mu^{-1}})^{\rm ad}.$$
 Let $\bbf{N}_{\mu}$ and $\bbf{N}_{\mu^{-1}}$ be the unipotent radicals of $\bbf{P}_{\mu}$ and $\bbf{P}_{\mu^{-1}}$ and let 
 $\bbf{M}$ be their common Levi subgroup. Note that all these objects depend on our choice of the representative 
 $\mu$ (for instance $\bbf{M}$ is the centralizer of $\mu$). 
 
The equivalence of algebraic stacks 
\[
[\Spec \breve{E'}/\mathbf{P}_{\mu}]= [\mathbf{G}_{ \breve{E'}}\backslash \mathbf{G}_{\breve{E'}}/\mathbf{P}_{\mu}]=[\mathbf{G}_{\breve{E'}}\backslash \FL_{\mu,E'}]
\] 
 depicted by the element $1\in \FL_{\mu, \breve{E'}}$ induces an equivalence
 \begin{equation}\label{eqEquivalenceQuasiCohSheavesFlagParabolic}
 \mathcal{W}_{\bbf{G},\mu}\colon {\rm Rep}(\mathbf{P}_{\mu})\simeq \mathrm{QCoh}([\Spec \breve{E'}/\mathbf{P}_{\mu}])\xrightarrow{\sim} \mathrm{QCoh}([\mathbf{G}_{\breve{E'}}\backslash \FL_{\mu, \breve{E'}}])
 \end{equation}
between algebraic $\mathbf{P}_{\mu}$-representations and $\mathbf{G}_{\breve{E'}}$-equivariant quasi-coherent sheaves on $\FL_{\mu, \breve{E'}}$.
 More concretely, this sends an algebraic $\mathbf{P}_{\mu}$-representation $W$ to $G_{\breve{E'}}\times^{\mathbf{P}_{\mu}} W$, and its inverse sends an equivariant quasi-coherent sheaf to its fibre at $1$.

  For any $V\in \Rep(\bbf{G})$ the $\bbf{P}_{\mu, \breve{E}'}$-representation $V_{\breve{E'}}:=V\otimes_{\bb{Q}_p} \breve{E'}$ endowed with its natural 
  decreasing $\bbf{P}_{\mu}$-stable filtration induces an object
  $$\n{F}_{V}^{\dR}= \s{O}_{\Fl_{\mu, \breve{E'}}} \otimes_{\breve{E'}} V_{\breve{E'}}\in  {\rm FIC}_{\Fl_{\mu, \breve{E'}}}$$
   By an abuse of notation we will also write $V_{\dR}$ for the restriction to the admissible locus 
   $\Fl_{\mu, \breve{E'}}^a\subset \Fl_{\mu, \breve{E'}}$. The object $\n{F}_{V}^{\dR}$ is endowed with a natural flat connection with basis given by horizontal sections $V_{\breve{E}'}$.

 We introduce the standard Lie algebroid notation of the localization theory of Beilinson-Bernstein, see \cite{BeilinsonBernstein}.  We let $$\f{g}^0_{\GM}:=\n{F}_{\f{g}}^{\dR}.$$   Let $\f{n}_{\mu}\subset \f{p}_{\mu} $ denote the Lie algebras of $\mathbf{N}_{\mu}\subset \mathbf{P}_{\mu}$ endowed with the adjoint action, and let $\f{n}_{\mu}^0\subset \f{p}_{\mu}^0$ be the associated equivariant sheaves over $\Fl_{\mu, \breve{E'}}$, similarly, we let $$\f{m}=\Lie \mathbf{M}_{\breve{E}'}$$ and let $\f{m}^0_{\mu}=\f{p}^0_{\mu}/\f{n}^0_{\mu}$ denote its associated sheaf over $\Fl_{\mu,\breve{E}'}$.  The natural anchor map  $\alpha\colon \f{g}^0_{\GM}\to \mathcal{T}_{\Fl_{\mu}}$  induces an isomorphism $$\f{g}^0_{\GM}/\f{p}^0_{\mu}\xrightarrow{\sim} \mathcal{T}_{\Fl_{\mu},\breve{E}'}.$$
 
  We have analogue objects $\f{g}^0_{\HT},\f{n}^0_{\mu^{-1}}, \f{p}^{0}_{\mu^{-1}}$ and $\f{m}^0_{\mu^{-1}}$ living over the Hodge-Tate flag variety $\Fl_{\mu^{-1},\breve{E}'}$, as well as an anchor map $\f{g}^0_{\HT}\to \mathcal{T}_{\Fl_{\mu^{-1}}}$ inducing an isomorphism $$\f{g}^0_{\HT}/\f{p}^0_{\mu^{-1}}\xrightarrow{\sim} \mathcal{T}_{\Fl_{\mu^{-1}}}.$$

\subsection{The Riemann-Hilbert correspondence on local Shimura varieties} \label{ss:RiemannHilbertLocalShimura}

In this section we will translate the diagram \eqref{eqHodgeTatePeriodMapv2} in terms of $p$-adic Hodge theory of local Shimura varieties and deduce a Riemann-Hilbert correspondence for pro\'etale local systems arising from algebraic $\bbf{G}$-representations. Contrary to the case of  global Shimura varieties, this is mostly formal and comes by unwinding definitions.

  Fix an algebraic representation $V\in \Rep(\mathbf{G})$ and consider the $v$-sheaf $V_{\proet}$ on $ \Fl_{\mu, \breve{E'}}^a$ obtained 
from the $\underline{\bbf{G}(\bb{Q}_p)}$-torsor $\pi_{\rm GM}: \n{M}_{\infty, \breve{E}'}\to \Fl_{\mu, \breve{E'}}^a$ (see  \Cref{DefinitionProetaleSystemsShimura}). By restriction we obtain a sheaf 
    denoted $\n{F}_{V}$ on $X:=\Fl_{\mu, \breve{E'}, \proet}^a$.

For any $S\in \Perf_k/\breve{\bb{Q}}_p^{\lozenge}$ the pullback of $\n{E}_{b}$ to the 
formal completion $X_S^{\wedge,\iota}$ of $X_S$ along the closed Cartier divisor $\iota: S^{\sharp}\to X_{S}$ is a trivial $\bbf{G}$-torsor. Since the automorphism group of the trivial $\bbf{G}$-torsor over $X_S^{\wedge,\iota}$  is $L^+\bbf{G}(S)$, pullback along the formal completion gives rise to a group homomorphism of $v$-sheaves  over $\breve{\bb{Q}}_p^{\diamond}$
\begin{equation}
\label{eqPullbackAutomorphismstoBdR}
\widetilde{G}_b\times \breve{\bb{Q}}_p^{\lozenge} \to L^+\bbf{G}. 
\end{equation}

   The Riemann-Hilbert correspondence for local Shimura varieties is then:

\begin{theorem}
\label{PropositionRiemannHilbert}
  The pro\'etale local system $\n{F}_{V}$ is de Rham, with associated filtered flat connection $$D_{\rm dR}(\n{F}_{V})\simeq \n{F}_{V}^{\dR}.$$
    More precisely, there is a $\widetilde{G}_b$-equivariant isomorphism of sheaves of filtered $\bb{B}_{\dR}$-modules on $ X_{\proet}$
\begin{equation}
\label{eqRiemannHilbertLocal}
\n{F}_{V}\otimes_{\bb{Q}_p} \bb{B}_{\dR} \cong  \underline{V}\otimes_{\bb{Q}_p} \bb{B}_{\dR} =(\n{F}^{\dR}_{V}\otimes_{\s{O}} \OBdR{X})^{\nabla=0},
\end{equation}
where the $\bb{B}_{\dR}^+$-filtration in the left-hand side is the trivial one, and the filtration in the right-hand side is given by 
\[
\Fil^i( \underline{V}\otimes_{\bb{Q}_p} \bb{B}_{\dR}):=(\Fil^i (\n{F}_{V}^{\dR}\otimes_{\s{O}_{X}} \OBdr))^{\nabla=0}.
\]
The action of $\widetilde{G}_b$ is trivial on $\n{F}_{V}$ in the left-hand side  and it factors through $\widetilde{G}_b\to L\bbf{G}$ and the natural action on $\underline{V}\otimes_{\bb{Q}_p} \bb{B}_{\dR}$ in the right-hand side. 
\end{theorem}

\begin{remark}
By \cite[Corollary 17.1.9]{ScholzeWeinspadicgeometry}, for a smooth rigid variety $Y$ one can identify filtered $\bb{B}^+_{\dR}$-vector bundles on the pro\'etale site $Y_{\proet}$ and filtered $\bb{B}_{\dR}^+$-vector bundles on $Y_v$. Thus, the equivariant isomorphism \eqref{eqRiemannHilbertLocal} can also be stated as a $\widetilde{G}_b$-equivariant isomorphism of filtered $\bb{B}_{\dR}$ sheaves on  the $v$-site
\[
\n{F}_{V}\otimes_{\bb{Q}_p} \bb{B}_{\dR} \cong \underline{V}\otimes_{\bb{Q}_p} \bb{B}_{\dR}.
\]
\end{remark}

\begin{proof}[Proof of \cref{PropositionRiemannHilbert}]
 Since $\n{F}_{V}^{\dR}$ has horizontal sections $V$, we have an isomorphism 
\[
\underline{V}\otimes_{\bb{Q}_p} \bb{B}_{\dR}=(\n{F}_{V}^{\dR}\otimes_{\s{O}_X} \OBdr)^{\nabla=0}
\]
in  $X_{\proet}$. By \cite[Theorem 7.6]{ScholzeHodgeTheory2013},  
$$\bb{M}':=(\Fil^0(\n{F}_{V}^{\dR}\otimes_{\s{O}_X} \OBdr))^{\nabla=0}$$ is a $\bb{B}_{\dR}^+$-lattice inside $\underline{V}\otimes_{\bb{Q}_p} \bb{B}_{\dR}$ in the pro\'etale site of $X$, and by \cite[Corollary 17.1.9]{ScholzeWeinspadicgeometry} we can view $\bb{M}'$ as a $\bb{B}_{\dR}^+$-lattice in the $v$-site of $X$. Thus, we will view $\underline{V}\otimes_{\bb{Q}_p} \bb{B}_{\dR}$ as a filtered $\bb{B}_{\dR}^+$-module in the $v$-site of $X$ with $\Fil^i(\underline{V}\otimes_{\bb{Q}_p} \bb{B}_{\dR})= \xi^i \bb{M}'$ for $\xi$ a local generator of the kernel of $\theta: \bb{B}_{\dR}^+\to \widehat{\s{O}}_{X}$. 

Now consider a map $S\to \n{M}_{\infty, \breve{E'}}$ with $S\in {\rm Perf}_k$, corresponding to a pair $(S^{\sharp}, \alpha)$. Let 
$\n{V}_1$ and $\n{V}_{b}$ be the vector bundles over the Fargues--Fontaine curve $X_{S}$ defined by the torsors $\n{E}_1$ and $\n{E}_b$ respectively, evaluated at the fixed representation $V$. Since $S$ lives over $\n{M}_{\infty, \breve{E'}}$ and $\n{F}_{V}$ is obtained by descent of 
$\underline{V}$ along the $\bbf{G}(\qp)$-torsor $\n{M}_{\infty, \breve{E'}}$, we have an isomorphism 
$$(\n{F}_{V}\otimes_{\bb{Q}_p} \bb{B}_{\dR})(S)\simeq V\otimes_{\qp} \bb{B}_{\dR}(S^{\sharp})\simeq  \n{V}_1\otimes_{\s{O}_{X_S}} \bb{B}_{\dR}(S^{\sharp}).$$
The modification $\alpha$ induces by completion 
an isomorphism 
\[
 \n{V}_1\otimes_{\s{O}_{X_S}} \bb{B}_{\dR}(S^{\sharp}) \cong \n{V}_{b}\otimes_{\s{O}_{X_S}} \bb{B}_{\dR}(S^{\sharp}) = (\underline{V}\otimes_{\bb{Q}_p} \bb{B}_{\dR})(S)
\]
Combining these yields 
a $\bbf{G}(\bb{Q}_p)\times \widetilde{G}_b$-equivariant isomorphism
$$(\n{F}_{V} \otimes_{\bb{Q}_p} \bb{B}_{\dR})(S)\simeq (\underline{V}\otimes_{\bb{Q}_p} \bb{B}_{\dR})(S),$$
where $\bbf{G}(\bb{Q}_p)$  acts trivially on  $\underline{V}$ and via the projection $\bbf{G}(\bb{Q}_p)\to L\bbf{G}$ on $\n{F}_{V} \otimes_{\bb{Q}_p} \bb{B}_{\dR}$, and $\widetilde{G}_b$  acts trivially on $\n{F}_{V}$ and via the projection $\widetilde{G}_b\to L\bbf{G}$ on $\underline{V}\otimes_{\bb{Q}_p} \bb{B}_{\dR}$. The fact that the isomorphism \eqref{eqRiemannHilbertLocal} is compatible with the filtration follows from the definition of the Bialynicki-Birula map and  \cite[Proposition 19.4.2]{ScholzeWeinspadicgeometry}. 
\end{proof}

  Using \cref{PropositionRiemannHilbert} and the definition of the Bialynicki-Birula map  
we obtain the following corollary.

\begin{corollary}
\label{CorollaryHodgeTateFiltration}
 There is a natural $\bbf{G}(\bb{Q}_p)\times \widetilde{G}_b$-equivariant isomorphism of $\widehat{\s{O}}$-modules over $\n{M}_{\infty, \breve{E}'}$ 
\[
{\rm F}_{{\rm HT}, n}(\n{F}_V) \cong \pi_{\HT}^*(\n{W}_{\bbf{G},\mu^{-1}}(\Fil_{n}(V_{\breve{E'}}))),
\]
where ${\rm F}_{{\rm HT}, \bullet}(\n{F}_V)$ is the relative Hodge-Tate filtration of $\n{F}_{V,\widehat{\s{O}}}=\n{F}_{V}\otimes_{\bb{Q}_p} \widehat{\s{O}}$ and
$\Fil_{\bullet}(V_{\breve{E'}})$ is the (increasing) $\bbf{P}_{\mu^{-1}}$-filtration of $V_{\breve{E'}}$ and $\n{W}_{\bbf{G},\mu^{-1}}$ is the equivalence in \eqref{eqEquivalenceQuasiCohSheavesFlagParabolic} (with $\mu$ replaced by $\mu^{-1}$).
\end{corollary}

 Taking graded pieces in the above corollary and using  \Cref{resumeSch}
one deduces the isomorphism of $\bbf{M}$-torsors on the infinite level Shimura variety. See \cite[Theorem 2.1.3]{CaraianiScholze2017} and \cite[Theorem 4.2.1]{RCLocAnCompleted} for the global (much harder) analogue.

\begin{corollary}
\label{CorollaryIsoTorsorsShimuraVariety}
Let $W\in \Rep(\bbf{M})$ be an irreducible algebraic representation over $\breve{E'}$, with highest weight $\mu(W)\in \bb{Z}$ with respect to $\mu$. There is a  natural $\bbf{G}(\bb{Q}_p)\times \widetilde{G}_b$-equivariant $\otimes$-isomorphism of $\widehat{\s{O}}$-modules over $\n{M}_{\infty, \breve{E}'}$  
\[
  \pi_{\HT}^*(\n{W}_{\bbf{G},\mu^{-1}}(W)) \cong   \pi_{\GM}^*(\n{W}_{\bbf{G},\mu}(W)) \otimes_{\widehat{\s{O}} } \widehat{\s{O}} (-\mu(W)).
\]
 In particular, if $\bbf{M}_{\GM}$ and $\bbf{M}_{\HT}$ denote the natural $\bbf{M}$-torsors living over $\Fl_{\mu,\breve{E}'}$ and $\Fl_{\mu^{-1}, \breve{E}'}$ respectively,  we have a $\bbf{G}(\bb{Q}_p)\times\widetilde{G}_b$-equivariant isomorphism of $\bbf{M}$-torsors over the ringed site $(\n{M}_{\infty,\breve{E}',v},\widehat{\s{O}})$
\[
 \pi_{\HT}^*(\bbf{M}_{\HT}) \cong \pi_{\GM}^*(\bbf{M}_{\GM})\times^{\bb{G}_m,\mu} \bb{G}_m(-1),
\]
where $\bb{G}_m$ injects into the center of $\bbf{M}$ via $\mu$, and $\bb{G}_m(-1)$ is the $\bb{G}_m$-torsor of trivializations of the Tate twist $\widehat{\s{O}}(-1)$. 
\end{corollary}

\subsection{Pullbacks of equivariant vector bundles along $\pi_{\rm HT}$} \label{ss:PullbackHTSheavesEquivariant}

In this section we compute the Faltings extension of the local Shimura varieties in terms of the representation theory of the Hodge-Tate flag variety.  We let $\mathcal{M}\subset \Fl_{\mu,\breve{E}}$ denote the admissible locus of the flag variety, and let $\mathcal{M}_{\infty}\to \mathcal{M}$ the $G(\bb{Q}_p)$-pro\'etale torsor given by the infinite level local Shimura variety.   Recall that $\gr^1 \OBdr^+$ is a $\widehat{\s{O}}$-vector bundle in the pro\'etale site, sitting in an  
  canonical exact sequence, called the Faltings extension 
 $$0\to \widehat{\s{O}}\to \gr^1 \OBdr^+ (-1)\to \Omega^1_{\n{M}} \otimes_{\s{O}_{\n{M}}} \widehat{\s{O}}(-1)\to 0.$$
 In particular $\gr^1 \OBdr^+$ defines naturally a $v$-vector bundle that we denote in the same way. The canonical map $e: \widehat{\s{O}}\to \gr^1\OBdr^+ (-1)$ induces a natural isomorphism 
\[
\gr^0 \OBdr = \Sym_{\widehat{\s{O}}} (\gr^1 \OBdr^+ (-1))/ (1-e(1)),
\]
where $1$ is the unit in the symmetric algebra, allowing us to consider $\gr^0 \OBdr$ as a $v$-sheaf which is a filtered colimit of $v$-vector bundles, see also \cite[Remark 2.1.2]{RCLocAnCompleted}.

  The algebra $\s{O}(\bbf{N}_{\mu^{-1}})$ of regular functions on
$\bbf{N}_{\mu^{-1}}$  has a left action 
of $\bbf{P}_{\mu^{-1}}=\bbf{N}_{\mu^{-1}}\rtimes \bbf{M}_{\mu^{-1}}$, whose restriction to $\bbf{N}_{\mu^{-1}}$ is the left regular action
and whose 
 restriction to $\bbf{M}_{\mu^{-1}}$ is the adjoint action $(m\cdot f)(n)= f(m^{-1}nm)$. Since 
 $\mu$ is minuscule, $\bbf{N}_{\mu^{-1}}$ is abelian and so can be identified as a group scheme with 
 its Lie algebra via the exponential map. This identifies $\s{O}(\bbf{N}_{\mu^{-1}})$ with the symmetric algebra of
 $\f{n}^{\vee}_{\mu^{-1}}$ and induces an increasing $\bbf{P}_{\mu^{-1}}$-filtration $\s{O}(\bbf{N}_{\mu^{-1}})^{\leq n}$ with graded pieces\footnote{See 
\cite[Proposition 3.3.1]{RCLocAnCompleted} and its previous discussion.} $$\gr_n(\s{O}(\bbf{N}_{\mu^{-1}}))\cong \Sym^n_{E'} (\f{n}^{\vee}_{\mu^{-1}}).$$ 

 The exact sequence 
$$0\to E'\to \s{O}(\bbf{N}_{\mu^{-1}})^{\leq 1}\to \f{n}^{\vee}_{\mu^{-1}}\to 0$$
induces an exact sequence 
$$0\to \widehat{\s{O}}\to \pi^*_{\HT}(\n{W}_{\bbf{G},\mu^{-1}}(\s{O}(\bbf{N}_{\mu^{-1}})^{\leq 1}))\to \pi^*_{\HT}(\f{n}^{0,\vee}_{\mu^{-1}})\to 0$$
of 
$\widehat{\s{O}}$-modules over $\mathcal{M}_{\infty, \breve{E'}, v}$. The next theorem shows that this is nothing but the  of the Faltings extension up to a sign.

\begin{theorem}
\label{TheoremOCCidentification}
\begin{enumerate}
\item The Kodaira-Spencer map of $\Fl_{\mu, \breve{E}'}$ induces a 
$\bbf{G}(\bb{Q}_p)\times \widetilde{G}_b$-equivariant isomorphism of $\widehat{\s{O}}$-modules over $\mathcal{M}_{\infty, \breve{E'}, v}$ 

\begin{equation}
\label{eqKodairaSpencer'}
\KS': \pi_{\HT}^*(\f{n}^{0,\vee}_{\mu^{-1}})\xrightarrow{\sim} \Omega^1_{\n{M}}\otimes_{\s{O}_{\n{M}}} \widehat{\s{O}}(-1)
\end{equation}
 
 \item  There is a natural $\bbf{G}(\bb{Q}_p)\times\widetilde{G}_b$-equivariant isomorphism of extensions\footnote{The bottom exact sequence is the Faltings extension.}
\[
\begin{tikzcd}
0\ar[r] & \widehat{\s{O}}\ar[r] \ar[d,"\id"'] & \pi^*_{\HT}(\n{W}_{\bbf{G},\mu^{-1}}(\s{O}(\bbf{N}_{\mu^{-1}})^{\leq 1}))\ar[d, "\alpha"'] \ar[r] & \pi^*_{\HT}(\f{n}^{0,\vee}_{\mu^{-1}}) \ar[r] \ar[d, "-\KS'"] & 0 \\ 
0 \ar[r] & \widehat{\s{O}}\ar[r] & \gr^1 \OBdr^+ (-1)  \ar[r] &\Omega^1_{\n{M}} \otimes_{\s{O}_{\n{M}_K}} \widehat{\s{O}}(-1) \ar[r] & 0
\end{tikzcd}
\]

inducing a $\bbf{G}(\bb{Q}_p)\times \widetilde{G}_b$-equivariant isomorphism of $\widehat{\s{O}}$-algebras over $\mathcal{M}_{\infty, \breve{E'}, v}$ 
\[
\gr^0(\OBdr) \cong \pi^*_{\HT}(\n{W}_{\bbf{G},\mu^{-1}}(\s{O}(\bbf{N}_{\mu^{-1}}))).
\]

\end{enumerate}
 
\end{theorem}

\begin{proof}
The proof follows fully and faithfully the arguments in \cite[Theorem 5.1.4]{RCLocAnCompleted}, where the key inputs are the Riemann-Hilbert correspondence of \cref{PropositionRiemannHilbert} and the Kodaira-Spencer isomorphism that we will explain below. 
 Note that in \textit{loc. cit.} $\gr^0(\OBdr)$ is denoted $\OCC{\n{S}h}$, and the Killing form of the derived Lie algebra $\f{g}^{\mathrm{der}}$ is used to obtain identifications
 $$\f{g}\cong \f{g}^{\vee},\,\, \f{m}\cong \f{m}^{\vee},\,\, \f{n}_{\mu^{-1}}^{\vee}\cong \f{g}_{E'}/\f{p}_{\mu^{-1}}.$$

Let $\f{g}$ be the adjoint representation of $\bbf{G}$ over $\bb{Q}_p$, and let $\n{F}_{\f{g}^{\vee}}^{\dR}= \s{O}_{\Fl_{\mu, \breve{E}'}}\otimes_{\bb{Q}_p} \f{g}^{\vee}=\f{g}^{0,\vee}_{\GM}$ be its associated $\bbf{G}$-equivariant vector bundle with flat connection over $\Fl_{\mu,\breve{E}'}$.  Since $\mu$ is minuscule,  $\f{g}^{0,\vee}_{\GM}$ has Hodge filtration concentrated in degrees $[-1,1]$ given by 
\[
(\f{g}^0_{\GM}/\f{p}^0_{\mu})^{\vee} \subset (\f{g}^0_{\GM}/\f{n}^0_{\mu})^{\vee} \subset \f{g}^{0,\vee}_{\GM}
\]
such that 
\[
\gr^i \f{g}^{0,\vee}_{\GM}=\begin{cases}
(\f{g}^0_{\GM}/\f{p}^0_{\mu})^{\vee} & \mbox{ if }i=1,\\
\f{m}^{0,\vee}_{\mu} & \mbox{ if }i=0,\\
\f{n}^{0,\vee}_{\mu} & \mbox{ if }i= -1. 
\end{cases}
\]
The flat connection $\nabla: \f{g}_{\GM}^{0,\vee}\to \f{g}_{\GM}^{0,\vee}\otimes_{\s{O}_{\Fl_{\mu, \breve{E}'}}} \Omega^1_{\Fl_{\mu, \breve{E}'}}$ induces a map between $\gr^1$-pieces
\[
(\f{g}_{\GM}^0/\f{p}^0_{\mu})^{\vee} \to \f{m}^{0,\vee}_{\mu}\otimes_{\s{O}_{\Fl_{\mu, \breve{E}'}}} \Omega^1_{\Fl_{\mu, \breve{E}'}}.
\]
Taking adjoints we get a $\bbf{G}$-equivariant map 
\begin{equation}
\label{eqKodairaSpencer1}
\widetilde{\KS}: (\f{g}^0_{\GM}/\f{p}^0_{\mu})^{\vee}\otimes_{\s{O}_{\Fl_{\mu, \breve{E}'}}}\f{m}^0_{\mu}\to \Omega^1_{\Fl_{\mu, \breve{E}'}}.
\end{equation}

Looking at the fiber at $[1]\in \Fl_{\mu, \breve{E}'}$ the map \eqref{eqKodairaSpencer1} is nothing but the natural adjoint action of $\f{m}_{\mu}$ on $(\f{g}_{E'}/\f{p}_{\mu})^{\vee}$ with $\f{g}_{E'}=\f{g}\otimes_{\bb{Q}_p} E'$:
\[
(\f{g}_{E'}/\f{p}_{\mu})^{\vee}\otimes_{E'}\f{m}_{\mu} \xrightarrow{\ad} (\f{g}_{E'}/\f{p}_{\mu})^{\vee} \cong \Omega^1_{\Fl_{\mu, \breve{E}'}}|_{[1]}.
\]

Therefore, the map $\widetilde{\KS}$ induces the $\bbf{G}$-equivariant  Kodaira-Spencer isomorphism over the flag variety 
\begin{equation}
\label{eqKodairaSpencer2}
\KS:(\f{g}^0_{\GM}/\f{p}^0_{\mu})^{\vee}\xrightarrow{\sim}  \Omega^1_{\Fl_{\mu, \breve{E}'}}.
\end{equation}
which is the inverse of the dual of the anchor map $\overline{\alpha}\colon \f{g}^0_{\GM}/\f{p}_{\mu}^0\xrightarrow{\sim} \mathcal{T}_{\Fl_{\mu, \breve{E}'}}$.\footnote{  We observe that $\KS$ is already defined over $E$ as the anchor map is so.}

We deduce the following lemma.

\begin{lemma}
\label{PropKodairaSpencer}
Let $K\subset \bbf{G}(\bb{Q}_p)$ be an open compact subgroup. The Kodaira-Spencer map of $\n{M}_{K,\breve{E}'}$
\[
\widetilde{\KS}: \gr^1(\f{g}^{0,\vee}_{\GM})\otimes_{\s{O}_{\n{M}_K}} \gr^0(\f{g}^0_{\GM})\to \Omega^1_{\n{M}_K}
\]
 constructed in analogue fashion as in \eqref{eqKodairaSpencer1}  factors through an isomorphism 
 \[
 \KS: \gr^1(\f{g}^{0,\vee}_{\GM})\xrightarrow{\sim} \Omega^1_{\n{M}_K}
 \]
 which is noting but the pullback along $\n{M}_{K,\breve{E}'}\to \Fl_{\mu, \breve{E}'}$ of the Kodaira-Spencer isomorphism \eqref{eqKodairaSpencer2}. 
\end{lemma}
\begin{proof}
This follows from the Kodaira-Spencer isomorphism \eqref{eqKodairaSpencer2} and the fact that the map $\n{M}_{K}\to \Fl_{\mu,\breve{E}}$ is \'etale, namely, the filtered vector bundle with flat connection $\f{g}^{\vee}_{\dR}$ over  $\n{M}_{K}$ is the pullback of the analogue filtered vector bundle with flat connection over the flag variety. 
\end{proof}

By  \cref{PropositionRiemannHilbert} the local system  $\n{F}_{\f{g}^{\vee}}$ on the admissible locus  $\mathcal{M}=\Fl_{\mu,\breve{E}}^a$ (see \cref{DefinitionProetaleSystemsShimura}) is de Rham with associated filtered flat connection $\n{F}_{\f{g}^{\vee}}^{\dR}$.  Then,  \Cref{CorollaryHodgeTateFiltration,CorollaryIsoTorsorsShimuraVariety}  give rise $\bbf{G}(\bb{Q}_p)\times \widetilde{G}_b$-equivariant $\widehat{\s{O}}$-linear isomorphisms on $\n{M}_{\infty,\breve{E}'}$
\begin{equation}
\label{EquationComparisonGradedAdjoint}
\pi_{\HT}^*(\f{n}^{0,\vee}_{\mu^{-1}})\cong \gr_1(\n{F}_{\f{g}^{\vee},\widehat{\s{O}}}) \cong \gr^1(\n{F}_{\f{g}^{\vee}}^{\dR})\otimes_{\s{O}_{\n{M}}} \widehat{\s{O}}(-1).
\end{equation}
Composing \eqref{eqKodairaSpencer2} and \eqref{EquationComparisonGradedAdjoint} we get the following incarnation of the Kodaira-Spencer map
\begin{equation}
\label{eqKodairaSpencer'}
\KS': \pi_{\HT}^*(\f{n}^{0,\vee}_{\mu^{-1}})\xrightarrow{\sim} \Omega^1_{\n{M}}\otimes_{\s{O}_{\n{M}}} \widehat{\s{O}}(-1)
\end{equation}
as $\bbf{G}(\bb{Q}_p)\times\widetilde{G}_b$-equivariant $\widehat{\s{O}}$-sheaves on $\n{M}_{\infty,\breve{E}',v}$.

 \end{proof}

\subsection{Computation of the geometric Sen operators}
\label{SubsectionSenOperators}

   For the pro\'etale $K$-torsor $\pi_K: \n{M}_{\infty, \breve{E}'}\to \n{M}_{K,\breve{E}'}$ the associated geometric Sen operator
\begin{equation}
\label{eqGeoSenOperatorlocalShimura}
\theta_{\n{M}_K}: \n{F}_{\f{g}^{\vee},\widehat{\s{O}}} \to \Omega^1_{\n{M}_K}\otimes_{\s{O}_{\n{M}_K}} \widehat{\s{O}}(-1)
\end{equation}
is described by following theorem.

\begin{theorem}
\label{TheoremComputationGeoSenOperator}
The geometric Sen operator \eqref{eqGeoSenOperatorlocalShimura} is the descent along the $K$-torsor $\pi_K: \n{M}_{\infty, \breve{E}'}\to \n{M}_{K,\breve{E}'}$ of the $\bbf{G}(\bb{Q}_p)\times \widetilde{G}_b$-equivariant map of $\widehat{\s{O}}$-vector bundles obtained by pulling back along $\pi_{\HT}$ the map
\[
\f{g}^{0,\vee}_{\HT}\to \f{n}^{0,\vee}_{\mu^{-1}}
\]
over $\Fl_{\mu^{-1}}$, under the identification $\pi_{\HT}^*(\f{n}^{0,\vee}_{\mu^{-1}})\cong \Omega^1_{\n{M}_K}\otimes_{\s{O}_{\n{M}_K}} \widehat{\s{O}}(-1)$ given by the Kodaira-Spencer isomorphism \eqref{eqKodairaSpencer'}. 
In particular $\theta_{\n{M}_K}$ is a surjectiver $\bbf{G}(\bb{Q}_p)\times \widetilde{G}_b$-equivariant map of $\widehat{\s{O}}$-sheaves over $\n{M}_{\infty, \breve{E}'}$.

\end{theorem}
\begin{proof}
The proof is the same as the one of \cite[Theorem 5.2.5]{RCLocAnCompleted} where  \cref{TheoremOCCidentification} replaces \cite[Theorem 5.1.4]{RCLocAnCompleted}. For the convenience of the reader we recall the key steps. First, by Remark 3.2.5 in \cite{RCGeoSen} it suffices to compute the geometric Sen operator on a faithful representation $V$ of the Lie algebra of $\bbf{G}$.  
 Using  \Cref{CorollaryHodgeTateFiltration}, we are reduced to computing the geometric Sen operator of the pullback under $\pi_{\HT}$ of the equivariant vector bundle attached to a $\mathbf{P}_{\mu^{-1}}$-representation $W$ (and apply this to the various pieces in the $\mathbf{P}_{\mu^{-1}}$-stable filtration of 
 $V_{\breve{E'}}$). We can reduce to the universal case $W=\mathcal{O}(\mathbf{P}_{\mu^{-1}})=\mathcal{O}(\mathbf{N}_{\mu^{-1}})\otimes_{\breve{E'}} 
 \mathcal{O}(\mathbf{M})$. By  \Cref{CorollaryIsoTorsorsShimuraVariety} and  \Cref{gsen} we are further reduced to the case 
 $W=\mathcal{O}(\mathbf{N}_{\mu^{-1}})$ (using the Peter-Weyl decomposition of the representation $ \mathcal{O}(\mathbf{M})$). But the $\pi_{\HT}$-pullback of the equivariant sheaf attached to this representation is 
 the sheaf $\OC$, whose geometric Sen operator is the negative of the residual connexion induced by $\OBdr$, and this corresponds to the natural connexion on 
 $\mathcal{O}(\mathbf{N}_{\mu^{-1}})$ induced by the derivative of the action of $\mathbf{N}_{\mu^{-1}}$.

\end{proof}

\section{Locally analytic vectors of infinite level Shimura varieties and de Rham cohomology}
\label{SectionLocAnInfiniteLevel}

  Let $(\bbf{G}, b, \{\mu\})$ be a local Shimura datum with reflex field $E$ and let 
  $E'/E$ and $\mu\in \{\mu\}$ defined over $E'$, such that $E'$ splits $\bbf{G}$. Let $C/E'$ be the $p$-adic completion of an algebraic closure of $E'$.    The purpose of this section is to study the decompletions by locally analytic vectors of the sheaves of periods at infinite level.
  
 \subsection{Independence of locally analytic vectors}\label{ss:IndependenceLocAn}
  
   We being with an immediate consequence of \Cref{theo:LocAnTwoTowers}. Suppose that $b$ is basic, and let    $(\bbf{\check{G}}, \check{b}, \check{\mu})$  be the dual local Shimura datum as in \Cref{PropositionDualityLocalShtukas}.  We let $\n{M}_{\infty}=\n{M}_{\mathbf{G},b,\mu,\infty}\cong \n{M}_{\check{\mathbf{G}}, \check{b},\check{\mu}, \infty}$ be the infinite level Shimura variety over $\breve{E}$.  The following theorem generalizes \cite[Corollary 5.3.9]{LuePanII}. 
  
  \begin{theorem}\label{TheoMainComparisonLocAn}
  Let $I\subset (0,\infty)$ denote a closed interval, and let $U\subset \n{M}_{C}$ be a quasi-compact open subspace. Let $K_U\subset \mathbf{G}(\bb{Q}_p)$ and $\check{K}_U\subset \check{\mathbf{G}}(\bb{Q}_p)$ be the stabilizers of $U$.   Then the natural maps 
  \[
  R\Gamma_v(U,\bb{B}_I)^{RK_U-{\rm la}}\xleftarrow{\sim} R\Gamma_v(U,\bb{B}_I)^{R (K_U\times \check{K}_U)-{\rm la}} \xrightarrow{\sim} R\Gamma_v(U, \bb{B}_I)^{R\check{K}_U-{\rm la}}
  \]
are equivalences of solid $K_U\times \widetilde{K}_{U}$-representations.   The same holds for the sheaf $\widehat{\s{O}}$. 
  \end{theorem}
  \begin{proof}
  This follows form \Cref{theo:LocAnTwoTowers} and the fact that both $U^{\diamond}/\underline{K_U}$ and $U^{\diamond}/\underline{\check{K}_{U}}$ are smooth quasi-compact and separated rigid spaces. 
  \end{proof}

\subsection{Vanishing of higher locally analytic vectors for $\widehat{\s{O}}$ and the arithmetic Sen operator}

 We assume that  $b$ is general. The following theorem describes the main properties of the locally analytic decompletion of the sheaf $\widehat{\s{O}}$ at infinite level.

\begin{theorem}\label{TheoVanishingLocAnVectors}
Let $U\subset \n{M}_{\infty,C}$  be a quasi-compact open subspace obtained as pullback from an affinoid $U_K\subset \n{M}_{\infty,C}$ at finite level $K\subset \mathbf{G}(\bb{Q}_p)$. Suppose in addition that $U_K$ admits toric coordinates\footnote{i.e. it admits an \'etale map to a torus $\bb{T}^d_C$, which factors as a composition of rational localizations and finite \'etale maps}. The following hold:

\begin{enumerate}
\item  The complex $R\Gamma_v(U,\widehat{\s{O}})^{RK-{\rm la}}$ 
 sits in degree $0$ and is equal to the locally analytic (with respect to the action of $K\subset \bbf{G}(\bb{Q}_p)$) vectors of $\widehat{\s{O}}(U)$.

\item The Lie algebroid $\f{n}^0_{\mu^{-1},C}\subset \f{g}^0_{\HT,C}=\s{O}_{\Fl_{\mu^{-1},C}}\otimes_{\bb{Q}_p} \f{g}$ kills $\widehat{\s{O}}(U_{\infty})^{\bbf{G}(\bb{Q}_p)-\rm la}$, hence we obtain a horizontal action of  $\f{m}^0_{\mu^{-1},C}$ on $\widehat{\s{O}}(U_{\infty})^{K-\rm la}$. 

\item  The space $\widehat{\s{O}}(U_{\infty})^{K-\rm la}$ has an arithmetic Sen operator as in \cite[Definition 7.1.2]{RCLocAnCompleted}, given by the opposite of the derivative of the Hodge cocharacter $-\theta_{\mu}=\theta_{\mu^{-1}}\in \f{m}^0_{\mu^{-1},C}$.

\end{enumerate}

\end{theorem}

\begin{proof} Parts (1) and (2) are direct consequences of the surjectivity of the geometric Sen operator ( \Cref{TheoremComputationGeoSenOperator}) and of  \Cref{geoSentorsor}, see also the proofs of Proposition 6.2.8 (1) and Corollary 6.2.13 in \cite{RCLocAnCompleted}.  The existence and computation of the arithmetic Sen operator follows the same path as the proof of Theorem  7.2.1 of \textit{loc. cit.}. 
\end{proof}

We let $\s{O}^{{\rm la}}_{\n{M}}$ denote the sheaf of locally analytic vectors of $\widehat{\s{O}}$ restricted to the underlying  topological space of $\n{M}_{\infty,C}$.

\subsection{The horizontal Cartan action}

 We keep the notation of Lie algebroids of \Cref{ss:SetUpNotationShimuraGeoRep}. Recall that at infinite level we have the period maps $\pi_{\GM}\colon \n{M}_{\infty}\to \Fl_{\mu,\breve{E}}$ and $\pi_{\HT}\colon \n{M}_{\infty}\to \Fl_{\mu^{-1},\breve{E}}$.    The following is the generalization of \cite[Corollary 5.3.13]{LuePanII}.

\begin{theorem}\label{TheoCenterBothSides}
The actions of $\f{m}^0_{\mu}$ and $\f{m}^0_{\mu^{-1}}$ on $\s{O}^{{\rm la}}_{\n{M}}$ by derivations are identified via the pullback 
\begin{equation}\label{eqomoaepmapwfnawrf}
\f{m}^0_{\mu^{-1}} \otimes_{\s{O}_{\Fl_{\mu^{-1}}}} \s{O}^{{\rm la}}_{\n{M}} =\f{m}^{0, {\rm la}} = \s{O}^{{\rm la}}_{\n{M}}\otimes_{\s{O}_{\Fl_{\mu}}} \f{m}^0_{\mu}.
\end{equation}
In particular, the central character of the actions of $\f{m}^0_{\mu}$ and $\f{m}^0_{\mu^{-1}}$ on $\s{O}^{{\rm la}}_{\n{M}}$ agree under the natural isomorphism of  the center of the enveloping algebras $\n{Z}(\f{m}_{\mu})_C \cong \n{Z}(\f{m}_{\mu^{-1}})_C$. 
\end{theorem}
\begin{proof}
In the following we forget about the action of the Galois group of $E$ and fix a trivialization of the Tate twist $\bb{Z}_p(1)\cong \bb{Z}_p$ obtained by fixing a sequence of $p$-th power roots of unit $(\zeta_{p^n})_{n}$.  In the following all completed tensor products are solid.

 Let $\n{M}^{{\rm la}}_{\infty,C}$ be the  ringed space whose underlying topological space is $|\n{M}_{\infty,C}|$ and sheaf of functions given by the algebra  $\s{O}^{{\rm la}}_{\n{M}}$. We have locally analytic Hodge-Tate period maps 
\[
\begin{tikzcd}
 & \n{M}^{{\rm la}}_{\infty,C} \ar[rd,"\pi_{\HT}^{{\rm la}}"] \ar[ld,"\pi_{\GM}^{{\rm la}}"'] & \\
\Fl_{\mu} & & \Fl_{\mu^{-1}}.
\end{tikzcd}
\]
Let $W$ be an algebraic representation of the Levi $\bbf{M}$, taking locally analytic vectors in \cref{CorollaryIsoTorsorsShimuraVariety} we get $\bbf{G}(\bb{Q}_p)\times \check{\bbf{G}}(\bb{Q}_p)$-equivariant $\otimes$-isomorphisms of vector bundles over $\n{M}^{{\rm la}}_{\infty,C}$
\begin{equation}\label{eqLocAnLeviTorsors}
\n{W}_{\bbf{G},\mu^{-1}}(W) \otimes_{\s{O}_{\Fl_{\bbf{G},\mu^{-1}}}}  \s{O}^{{\rm la}}_{\n{M}} =  \s{O}^{{\rm la}}_{\n{M}}  \otimes_{\s{O}_{\Fl_{\mu}}} \n{W}_{\bbf{G},\mu}(W).
\end{equation}
Let $\bbf{M}_{\mu}\to \Fl_{\mu}$ and  $\bbf{M}_{\mu^{-1}}\to \Fl_{\mu^{-1}}$  be the natural $\bbf{M}$-torsors, the $\otimes$-equivalence \eqref{eqLocAnLeviTorsors} gives rise to a natural isomorphism of $\bbf{M}$-torsors over $\n{M}^{{\rm la}}_{\infty,C}$
\[
\pi_{\HT}^{{\rm la},*}(\bbf{M}_{\mu^{-1}})\cong \pi_{\GM}^{{\rm la},*}(\bbf{M}_{\mu}). 
\]
Thus, if $\bbf{M}^{{\rm an}}_{\mu}$ and $\bbf{M}_{\mu}^{{\rm an}}$ denote the analytification of the algebraic torsors over the flag varieties, the period maps refine to a mixed period map 
\begin{equation}\label{eqwojoanflasasf}
\pi^{{\rm la}}_{\GM,\HT}: \n{M}^{{\rm la}}_{\infty,C} \to \bbf{M}^{{\rm an}}_{\mu}\times^{\bbf{M}^{{\rm an}}} \bbf{M}^{{\rm an}}_{\mu^{-1}}. 
\end{equation}
Note that the Lie algebra $\f{g}\times \f{g}_b$ acts on $\bbf{M}^{{\rm an}}_{\mu}\times^{\bbf{M}^{{\rm an}}} \bbf{M}^{ {\rm an}}_{\mu^{-1}}$ by derivations, and by construction both horizontal actions $\f{m}^0_{\mu}$ and $\f{m}_{\mu^{-1}}^0$ are identified after pullback (similarly for the infinitesimal actions of $\n{Z}(\f{m}_{\mu})_{C}$ and $\n{Z}(\f{m}_{\mu^{-1}})_{C}$).

Let $K_p\subset \bbf{G}(\bb{Q}_p)$  be a compact open subgroup. Consider $\s{O}^{{\rm la}}_{\n{M}}$ as an analytic sheaf on $\n{M}_{K_p,C}$. By  equation (6.25) in the proof of \cite[Corollary 6.2.14]{RCLocAnCompleted},  there is a natural $K_p$-equivariant isomorphism of pro\'etale $\widehat{\s{O}}$-modules on $\n{M}_{K_p,C}$
\begin{equation}\label{eqoojqonmaelfae}
\s{O}^{{\rm la}}_{\n{M}}\widehat{\otimes}_{\s{O}_{\n{M}_{K_p}}} \widehat{\s{O}} \cong C^{{\rm la}}(K_p,\widehat{\s{O}})^{\f{n}^0_{\mu^{-1},\star_1}=0}
\end{equation}
where $\star_1$ refers to the left action by derivations, and $K_p$ acts on the sheaf $\s{O}^{{\rm la}}_{\n{M}}$ on the left and by the right regular action on the right. The equivalence is induced by the orbit map.   

Now, the action of $\f{m}^{0}_{\mu}$ by derivations on $\s{O}^{{\rm la}}_{\n{M}}$ is $\s{O}_{\n{M}_{K_p}}$-linear and  $K_p$-invariant,  thus by \eqref{eqoojqonmaelfae} it must factors through a morphism of Lie algebroids
\[
\f{m}^0_{\mu}\to \f{g}^0_{\mu^{-1}}/\f{n}^0_{\mu^{-1}}\otimes_{\s{O}_{\Fl_{\mu^{-1}}}} \s{O}^{{\rm la}}_{\n{M}}.
\]
By looking at Hodge-Tate weights, we even must have a factorization $\f{m}^0_{\mu}\to \f{m}^0_{\mu^{-1}}\otimes_{\s{O}_{\Fl_{\mu^{-1}}}} \s{O}^{{\rm la}}_{\n{M}}$. Therefore, to determine the action of $\f{m}^0_{\mu}$ it suffices to see how it acts on $\bbf{M}$-representations via the torsor $\bbf{M}_{\mu^{-1}}$, but the map  \eqref{eqwojoanflasasf} implies that this action is the same as that arising from the torsor $\bbf{M}_{\mu}$, proving that both horizontal actions are are compatible in the sense of \eqref{eqomoaepmapwfnawrf} proving what we wanted. 
\end{proof}

\subsection{De Rham cohomology of infinite level local Shimura varieties}\label{SubSec:deRhamCoho}

In this section we show that the sheaf $\s{O}^{\rm la}_{\n{M}}$ produces an isomorphism between the de Rham cohomology (with compact supports) of the two towers for a duality of local Shimura varieties.  Similar results have been obtained independently by  Bosco-Dospinescu-Niziol \cite{FlipFlop} and Benchao Su \cite{BenchaoSudeRham}. In order to state the theorem, we keep the notation of  \Cref{ss:SetUpNotationShimuraGeoRep}. Let $\bbf{M}_{\mu}$ and $\bbf{M}_{\mu^{-1}}$ be the natural $\bbf{M}$-torsors over $\Fl_{\mu}$ and $\Fl_{\mu^{-1}}$ respectively. 

 Let $$X=\bbf{M}^{\rm an}_{\mu}\times^{ \bbf{M}^{\rm an}} \bbf{M}^{\rm an}_{\mu^{-1}}$$ and consider the mixed Lie algebroid over $\s{O}^{\rm la}_{\n{M}}$ 
\[
\n{T}^{\rm la}= \n{T}_{X}\otimes_{\s{O}_X} \s{O}^{\rm la}_{\n{M}} 
\]
obtained as the pullback of the tangent bundle of $X$ via the map $\pi^{\rm la}_{\GM,\HT}$ of \eqref{eqwojoanflasasf}.  The next result is a re-interpretation of  \Cref{TheoCenterBothSides}.

\begin{prop}
The Lie algebroid $\n{T}^{\rm la}$ acts by derivations on $\s{O}^{\rm la}_{\n{M}}$, compatible with the derivations on $X$.

\end{prop}

\begin{proof}

Let $\f{n}_{\mu^{-1}}^{0,\rm la}\subset \f{p}^{0,\rm la}_{\mu^{-1}}\subset \f{g}^{0,\rm la}_{\HT}$ be the base change of $\f{n}^0_{\mu^{-1}}\subset \f{p}^{0}_{\mu^{-1}}\subset \f{g}^{0}_{\HT}$ from $\s{O}_{\Fl_{\mu^{-1}}}$  to $\s{O}^{\rm la}_{\n{M}}$ and let $\f{n}_{\mu}^{0,\rm la}\subset \f{p}^{0,\rm la}_{\mu}\subset \f{g}^{0,\rm la}_{\GM}$  be the base change of $\f{n}^0_{\mu}\subset \f{p}^{0}_{\mu}\subset \f{g}^{0}_{\GM}$ from $\s{O}_{\Fl_{\bbf{G},\mu}}$  to $\s{O}^{\rm la}_{\n{M}}$. \cref{TheoCenterBothSides} provides an isomorphism 
\[
\f{m}^0_{\mu^{-1}}\otimes_{ \s{O}_{\Fl_{\mu^{-1}}}} \s{O}^{\rm la}_{\n{M}} \cong \f{m}^{0,\rm la}\cong  \s{O}^{\rm la}_{\n{M}}  \otimes_{ \s{O}_{\Fl_{\mu}}} \f{m}^0_{ \mu}. 
\]
On the other hand, $\n{T}^{\rm la}$ is the quotient of $\f{g}^{0,\rm la}_{\GM}\oplus \f{g}_{\HT}^{0,\rm la}$ by the Lie algebroid $\widetilde{\f{p}}^{0,\rm la}$ sitting in the cartesian square
\[
\begin{tikzcd}
\widetilde{\f{p}}^{0,\rm la} \ar[r] \ar[d] & \f{m}^{0,\rm la} \ar[d, "(\iota{,}-\iota)"] \\
\f{p}^{0,\rm la}_{\mu} \oplus \f{p}^{0,\rm la}_{\mu^{-1}} \ar[r] & \f{m}^{0,\rm la}_{\mu} \oplus \f{m}^{0,\rm la}_{\mu^{-1}}
\end{tikzcd}
\] 
where $(\iota,-\iota)$ is the anti-diagonal map.  Since $\widetilde{\f{p}}^{0,\rm la}$ acts trivially on $\s{O}^{\rm la}_{\n{M}}$, the action by derivations of $\f{g}^{0,\rm la}_{\GM}\oplus \f{g}_{\HT}^{0,\rm la}$ on $\s{O}^{\rm la}_{\n{M}}$ descends to $\n{T}^{\rm la}$. 
\end{proof}

\begin{remark}
With some additional effort one can prove that $\s{O}^{\rm la}_{\n{M}}$ is formally smooth over $C$ and that its tangent space is given by  $\n{T}^{\rm la}$ but we will not need this fact for the applications in  this paper. 
\end{remark}
 
 Let $\s{DR}(\s{O}^{\bbf{G}(\bb{Q}_p)-\rm sm}_{\n{M}})$ (resp. $\s{DR}(\s{O}^{G_b-\rm sm}_{\n{M}})$) be the colimit of de Rham complexes 
 of the finite level local Shimura varieties $\n{M}_{\bbf{G},b,\mu,K_p}$ (resp. $\n{M}_{\Check{\bbf{G}},\Check{b},\Check{\mu},K_{b,p}}$) over the compact open subgroups $K_p$ (resp. $K_{b,p}$) of $\bbf{G}(\bb{Q}_p)$ (resp. $\Check{\bbf{G}}(\bb{Q}_p)=G_b$).
 
 \begin{theorem}\label{TheoComparisonCoho}

\begin{enumerate}

\item  There are natural $\bbf{G}(\bb{Q}_p)\times G_b$-equivariant  quasi-isomorphisms of complexes over the topological space $|\n{M}_{\infty,C}|$
\begin{equation}\label{eqComparisondeRhamCohomologies}
\s{DR}(\s{O}^{\bbf{G}(\bb{Q}_p)-\rm sm}_{\n{M}}) \cong R\Gamma(\n{T}^{\rm la}, \s{O}^{\rm la}_{\n{M}}) \cong \s{DR}(\s{O}^{\check{\bbf{G}}(\bb{Q}_p)-\rm sm}_{\n{M}}),
\end{equation}
where $R\Gamma(\n{T}^{\rm la}, \s{O}^{\rm la}_{\n{M}})$ is the de Rham cohomology of $\s{O}^{\rm la}_{\n{M}}$ with respect to the action by derivations of the Lie algebroid $\n{T}^{\rm la}$.

\item  There is a natural $\bbf{G}(\bb{Q}_p)\times  \check{\bbf{G}}(\bb{Q}_p)$-equivariant isomorphism of de Rham cohomologies with compact supports
\[
\varinjlim_{K_p} H^i_{\rm dR,c}(\n{M}_{\bbf{G},b,\mu,K_p,C}) \cong \varinjlim_{K_{b,p}} H^i_{\rm dR,c}(\n{M}_{\Check{\bbf{G}},\Check{b},\Check{\mu},K_{b,p},C}).
\]

\end{enumerate}

\end{theorem}

\begin{remark}
Using the theory of Gelfand stacks  of \cite{deRhamStacksFF},  \cref{TheoComparisonCoho} can be promoted to an equivalence of analytic stacks
\[
\varprojlim_{K_p} \n{M}_{\bbf{G},b,\mu,K_p}^{\rm dR} \cong  \n{M}_{\bbf{G},b,\mu,\infty}^{\rm dR} \cong  \varprojlim_{K_{b,p}} \n{M}_{\Check{\bbf{G}},\Check{b},\Check{\mu},K_{b,p}}^{\rm dR}.
\]
After taking quotients by the smooth groups $\bbf{G}(\bb{Q}_p)^{\rm sm}$ and $\Check{\bbf{G}}(\bb{Q}_p)^{sm}$ such an equivalence  produces an equivalence of analytic stacks 
\[
\Fl^{a,\rm dR}_{\mu,\breve{E}}/ \Check{\bbf{G}}(\bb{Q}_p)^{\rm sm} = \Fl^{a,\rm dR}_{\mu^{-1},\breve{E}}/ \bbf{G}(\bb{Q}_p)^{\rm sm}
\]
between the analytic de Rham stacks of the quotients of the admissible locus of the flag varieties. This gives rise to a \textit{Jacquet-Langlands equivalence} of equivariant analytic $D$-modules. 
\end{remark}

\begin{proof}[Proof of \cref{TheoComparisonCoho}]
 Let $\n{T}_{\mu}$ and $\n{T}_{\mu^{-1}}$ be the tangent spaces of $\Fl_{\mu}$ and $\Fl_{\mu^{-1}}$ respectively. We have identifications $\n{T}_{\mu}= \f{g}^{0}_{\GM}/\f{p}^0_{\mu}$ and $\n{T}_{\mu^{-1}}=\f{g}^0_{\HT}/\f{p}^0_{\mu^{-1}}$ via the anchor map.  By construction of the Lie algebroid $\n{T}^{\rm la}$, we have a short exact sequence
\[
0\to  \f{m}^{0,\rm la}\to \n{T}^{\rm la}\to \f{g}^{0,\rm la}_{\GM}/\f{p}^{0,\rm la}_{\mu} \oplus \f{g}^{0,\rm la}_{\HT}/\f{p}^{0,\rm la}_{\mu^{-1}} \to 0. 
\]
The pullback along the inclusion of  $ \f{g}^{0,\rm la}_{\GM}/\f{p}^{0,\rm la}_{\mu}$ in the direct sum corresponds to the Lie algebroid $\f{g}^{0,la}_{\GM}/\f{n}^{0,\rm la}_{\mu}$ (similarly the pullback for the inclusion of $\f{g}^{0,\rm la}_{\HT}/\f{p}^{0,\rm la}_{\mu^{-1}}$ is $\f{g}^{0,\rm la}_{\HT}/\f{n}^{0,\rm la}_{\mu^{-1}}$). Thus,  we can write the $\n{T}^{\rm la}$-de Rham complex as the composite 
\[
 R\Gamma( \f{g}^{0,\rm la}_{\HT}/\f{p}^{0,\rm la}_{\mu^{-1}} , R\Gamma(\f{g}^{0,\rm la}_{\GM}/\f{n}^{0,\rm la}_{\mu}, \s{O}^{\rm la}_{\n{M}})) \cong  R\Gamma(\n{T}^{\rm la}, \s{O}^{\rm la}_{\n{M}}) \cong R\Gamma(\f{g}^{0,\rm la}_{\GM}/\f{p}^{0,\rm la}_{\mu} ,R\Gamma(\f{g}^{0,\rm la}_{\HT}/\f{n}^{0,\rm la}_{\mu^{-1}},  \s{O}^{\rm la}_{\n{M}})).
\]

 Part (1) of the theorem follows then from \Cref{PropositionComparisondeRhamSmooth}.  
The claim about the cohomology comparisons for the de Rham cohomology with compact supports follows for example by using the definition of compactly supported de Rham cohomology  arising from the six functor formalism of analytic $D$-modules of \cite{camargo2024analytic}. One can also argue by using the ad-hoc definition of \cite{zbMATH01414388}. Indeed, the compactly supported cohomology of the de Rham complex of  \textit{loc. cit.}  is nothing but the compactly supported cohomology of the de Rham complex seen as a sheaf on the underlying Berkovich space of $\n{M}_{\bbf{G},b,\mu,\infty,C}$. To see that this cohomology with compact supports is well defined one can argue as follows: the map $\n{M}_{\bbf{G},b,\mu,\infty,C}\to \n{M}_{\bbf{G},b,\mu,K_p,C}$ gives rise to a $K_p$-torsor of Berkovich spaces
\begin{equation}\label{eqtorsorBerkovich}
\n{M}_{\bbf{G},b,\mu,\infty,C}^B\to \n{M}_{\bbf{G},b,\mu,K_p,C}^B.
\end{equation}
The space $\n{M}_{\bbf{G},b,\mu,K_p,C}^B$ is a   locally finite dimensional  Hausdorff space (being the Berkovich space of a rigid space) and \cite[Theorem 4.8.9 (i)]{heyer20246functorformalismssmoothrepresentations} implies that  $\n{M}_{\bbf{G},b,\mu,K_p,C}^B$ has a well define functor of cohomology with compact supports for  sheaves over $\bb{Q}_p$ (in the language of \textit{loc. cit.}  it is \textit{$\bb{Q}_p$-fine}). Since \eqref{eqtorsorBerkovich} is represented in profinite sets, $\n{M}_{\bbf{G},b,\mu,\infty,C}^B$ is also a $\bb{Q}_p$-fine map (this follows from \cite[Theorem 3.4.11 (ii)]{heyer20246functorformalismssmoothrepresentations} since any map between profinite sets is $\bb{Q}_p$-fine by construction, see Section 3.5.16 in \textit{loc. cit.}), i.e. it has a well defined functor of cohomology with compact supports.
\end{proof}

\begin{prop}\label{PropositionComparisondeRhamSmooth}

 The natural maps $\s{O}_{\n{M}}^{\bbf{G}(\bb{Q}_p)-\rm sm}\to R\Gamma(\f{g}^0_{\HT}/\f{n}_{\mu^{-1}}^0, \s{O}^{\rm la}_{\n{M}})$ and 
  $\s{O}_{\n{M}}^{\Check{\bbf{G}}(\bb{Q}_p)-\rm sm}\to R\Gamma(\f{g}^0_{\GM}/\f{n}_{\mu}^0, \s{O}^{\rm la}_{\n{M}})$
 are equivalences.

\end{prop}

\begin{proof} These claims are symmetric  with respect to the period maps, so it suffices to prove the first.  Let 
$$\mathcal{E}=\n{W}_{\bbf{G},\mu^{-1}} ( \s{O}(\bbf{N}_{\mu^{-1}}))$$
be the $\bbf{G}$-equivariant sheaf on $\Fl_{\mu^{-1}}$ attached to the $\bbf{P}_{\mu^{-1}}$-representation
$ \s{O}(\bbf{N}_{\mu^{-1}})$, so that 
by \cref{TheoremOCCidentification}
\begin{equation}\label{eqIsoFE2}
\gr^0(\OBdr)=\pi_{\HT}^* (\mathcal{E}). 
\end{equation}
Let $\nu_{K_p}: \n{M}_{K_p , C,\proet}\to \n{M}_{K_p,C,\an}$ be the natural projection of sites. By \cite[Proposition 6.16]{ScholzeHodgeTheory2013} one has a natural quasi-isomorphism
\[
R\nu_{K_p,*} \gr^0(\OBdR{\n{M}}) = \s{O}_{\n{M}_{K_p,C}} .
\]

  Let now $V\subset \n{M}_{K_p,C}$ be an open affinoid and let $V_{\infty}=\n{M}_{\infty,C}\times_{\n{M}_{K_p,C}} V$ be its pullback at infinite level. 
  We have the following chain of quasi-isomorphisms
  \[
\begin{aligned}
\s{O}_{\n{M}_{K_p,C}}(V) & = R\Gamma_{\proet}(V, \gr^0(\OBdR{\n{M}}))  \\ & = R\Gamma(K_p, R\Gamma_{\proet}(V_{\infty},\gr^0(\OBdR{\n{M}}) ) ) \\ 
& = R\Gamma(K_p, R\Gamma_{\proet}(V_{\infty},\gr^0(\OBdR{\n{M}}) )^{RK_p-\rm la} ) \\  
& = R\Gamma(K_p, R\Gamma_{\proet}(V_{\infty},\widehat{\s{O}}_{\n{M}} )^{RK_p-\rm la} \otimes_{\s{O}_{\Fl_{\bbf{G},\mu^{-1}}}}  \mathcal{E} )\\
& = R\Gamma(K_p, \s{O}^{\rm la}_{\n{M}}(V_{\infty})\otimes_{\s{O}_{\Fl_{\mu^{-1}}}}   \mathcal{E}) \\ 
& = R\Gamma^{\rm sm}(K_p, R\Gamma(\f{g},\s{O}^{\rm la}_{\n{M}}(V_{\infty})\otimes_{\s{O}_{\Fl_{\mu^{-1}}}}  \mathcal{E})),
\end{aligned}
\]
where:

\begin{itemize}

\item The second equality is decent along the $K_p$-torsor $V_{\infty}\to V$.

\item  The third equality is the comparison between solid and locally analytic group cohomology  \cite[Theorem 6.3.4]{RJRCSolidLocAn2}.

\item  The fourth equality follows from the projection formula of locally analytic vectors \cite[Corollary 3.1.15 (3)]{RJRCSolidLocAn2}, the isomorphism \eqref{eqIsoFE2}, and projection formula for vector bundles.

\item  The fifth equality follows from the vanishing of higher locally analytic vectors in \cref{TheoVanishingLocAnVectors}.

\item  The sixth equality is the Lie algebra/smooth vs locally analytic cohomology comparison in Theorem \cite[Theorem 6.3.4]{RJRCSolidLocAn2}

\end{itemize}

Taking colimits as $K_p\to 1$, we deduce that 
\begin{equation}\label{eqLieAlgebra1}
\s{O}^{\bbf{G}(\bb{Q}_p)-\rm sm}_{\n{M}}=\varinjlim_{K_p} R\nu_{K_p,*} \gr^0(\OBdr) = R\Gamma(\f{g},\s{O}^{\rm la}_{\n{M}}\otimes_{\s{O}_{\Fl_{\mu^{-1}}}}   \mathcal{E}) .
\end{equation}
Since $\bbf{N}_{\mu^{-1}}$ is an affine space, we have
\[
R\Gamma(\f{n}_{\mu^{-1}}, \s{O}(\bbf{N}_{\mu^{-1}})_{C})=C.
\]
Taking the associated equivariant vector bundles over the flag variety and taking pullbacks along $\pi_{\HT}$ one deduces that 
\[
R\Gamma(\f{n}^{0,\rm la}_{\mu^{-1}},\s{O}^{\rm la}_{\n{M}}\otimes_{\s{O}_{\Fl_{\mu^{-1}}}}  \mathcal{E}  ) = \s{O}^{\rm la}_{\n{M}}.
\]
Combining this with \eqref{eqLieAlgebra1}, and by computing $\f{g}^0_{\HT}$-lie algebra cohomology in two steps,  one gets that 

$$
\s{O}^{\bbf{G}(\bb{Q}_p)-\rm sm}_{\n{M}} = R\Gamma(\f{g}^{0,\rm la}_{\HT}/\f{n}^0_{\mu^{-1}}, R\Gamma(\f{n}^{0,\rm la}_{\mu^{-1}},\s{O}^{\rm la}_{\n{M}}\otimes_{\s{O}_{\Fl_{\mu^{-1}}}}  \mathcal{E} )) 
= R\Gamma(\f{g}^{0,\rm la}_{\HT}/\f{n}^{0,\rm la}_{\mu^{-1}}, \s{O}^{\rm la}_{\n{M}}),
$$
finishing the proof. 
\end{proof}

\section{The Jacquet-Langlands functor for admissible locally analytic  representations}
\label{ssJLFunctor}

After recalling the definition of the Jacquet-Langlands functor of \cite{ScholzeLubinTate} for admissible Banach representations, we prove that this functor is compatible with the passage to locally analytic vectors and that it has a natural compatibility with the action of the center of the enveloping algebra.

\subsection{Scholze's Jacquet-Langlands functor}

Let $F/\bb{Q}_p$ be a finite extension, with ring of integers $\n{O}\subset F$, uniformizer $\varpi\in \n{O}$ and residue field $\bb{F}=\bb{F}_q$, and let $\bb{C}_p$ be the $p$-completion of an algebraic closure of $F$. Let $n$ be a positive integer and consider the group $\bbf{G}=\GL_{n,F}$ and the cocharacter $\mu$ given by $(1,0,\ldots, 0)$, with $n-1$ occurrences of $0$. Let $[b]\in B(\bbf{G})$ correspond to a formal $p$-divisible $\n{O}$-module $\bb{X}_b$ over $\overline{\bb{F}}$ of dimension $1$ and $F$-height $n$. We have $\widetilde{G}_b=D^{\times}$, where 
$D$ is the division algebra over $F$ of invariant $1/n$.  
Let $\n{M}_{\infty}=\n{M}_{\GL_{n,F},[b],\{\mu\},\infty}$.

By \cite{RZ96} there is a formal scheme $\f{M}_{\bb{X}}$ over $\Spf \breve{\n{O}}$, formally smooth and locally formally of finite type representing the functor on formal schemes over $\breve{\n{O}}$ sending $S$ to the set of isomorphism classes of pairs $(X,\rho)$, where $X/S$ is a formal $\n{O}$-module, and $\rho: X\times_{S}\overline{S}\xrightarrow{\sim} \bb{X}\times_{\overline{\bb{F}}} \overline{S}$ is a quasi-isogeny of formal $\n{O}$-modules, where $\overline{S}=S\times_{\Spf \breve{\n{O}}} \Spec \overline{\bb{F}}$. Let $\n{M}_{\bb{X}}$ denote the generic fiber of $\f{M}_{\bb{X}}$ as a rigid space. This is the Lubin-Tate rigid space over $\breve{F}$, an infinite disjoint union of open polydiscs of dimension $n-1$.

\begin{theorem}[{\cite[Corollary 24.3.5]{ScholzeWeinspadicgeometry}}]
There is a natural equivalence of diamonds $\n{M}^{\lozenge}_{\bb{X}}\cong \n{M}_{K}$ with $K= \GL_n(\n{O})$. 
\end{theorem}

 In this situation the $\GL_n(F)\times D^{\times}$-equivariant period maps \eqref{eqHodgeTatePeriodMapv2} restrict to a diagram 
\begin{equation*}
\begin{tikzcd}
 & \n{M}_{\infty} \ar[rd,"\pi_{\HT}"] \ar[ld, "\pi_{\GM}"']& \\ 
\bb{P}^{n-1}_{\breve{F}}  & &  \Omega_{\breve{F}}
\end{tikzcd}
\end{equation*}
where
\begin{itemize}

\item $\pi_{\GM}$ is a pro\'etale $\GL_{n}(F)$-torsor and $D^{\times}$ acts on $\bb{P}^{n-1}_{\breve{F}}$ via the natural inclusion of the map $D^{\times}\subset \GL_{n}(\breve{F})$. 

\item  $\pi_{\HT}$ is a pro\'etale $D^{\times}$-torsor and $\Omega_{\breve{F}}\subset \bb{P}^{n-1}_{\breve{F}}$ is the $\GL_{n}(F)$-stable open Drinfeld space obtained by removing all $F$-rational hyperplanes. 

\end{itemize}

Thus, we have an equivalence of $v$-stacks
\[
[\bb{P}^{n-1,\lozenge}_{\breve{F}}/\underline{D}^{\times}] \cong [\Omega_{\breve{F}}^{\lozenge}/\underline{\GL_n(F)}]. 
\]

The Jacquet-Langlands functor is defined as follows.

\begin{definition}\label{Def:AdmissibleTorsion}
 Let $\pi$ be a $p^{\infty}$-torsion admissible representation of $\GL_{n}(F)$ over $\bb{Z}_p$ and let $\n{F}_{\pi}$ be the \'etale sheaf over $\bb{P}^{n-1}_{\breve{F}}$ obtained by descent along $\pi_{\GM}$. The Jacquet-Langlands functor $\n{JL}$ is the functor mapping such $\pi$ to the complex
\[
\n{JL}(\pi)= R\Gamma_{\et}(\bb{P}^{n-1}_{\bb{C}_p}, \n{F}_{\pi}). 
\]
\end{definition}

\begin{theorem}[{\cite[Theorem 1.1]{ScholzeLubinTate}}]\label{Theo:AdmissibleLTScholze}
For any $p^{\infty}$-torsion admissible smooth representation $\pi$ of $\GL_n(F)$ over $\bb{Z}_p$ and any
$i\in \mathbf{Z}$ 
\[
\n{JL}^i(\pi):= H^i_{\et}(\bb{P}^{n-1}_{\bb{C}_p}, \n{F}_{\pi})
\]
is an admissible smooth representation of $D^{\times}$

\end{theorem}

For convenience we shall consider the $p$-completed analogue of \cref{Theo:AdmissibleLTScholze}. Let $\pi$ be a $p$-adically complete admissible representation of $\GL_n(F)$ over $\bb{Z}_p$, and let $\n{F}_{\pi}$ be the pro\'etale sheaf over $\bb{P}^{n-1}_{\breve{F}}$ given by the limit of \'etale sheaves $\n{F}_{\pi} = \varprojlim_{s} \n{F}_{\pi/p^s}$. Finally, we denote by $\n{JL}(\pi)$ the  $p$-adically complete $D^{\times}$-representation 
\[
\n{JL}(\pi):=R\Gamma_{\proet}(\bb{P}^{n-1}_{\bb{C}_p}, \n{F}_{\pi})= R\varprojlim_{s} R\Gamma_{\et} (\bb{P}^{n-1}_{\bb{C}_p}, \n{F}_{\pi/p^s}). 
\]
\begin{corollary}\label{CoroBanachAdmissible}
If $\pi$ is a $p$-adically complete admissible representation of $\GL_{n}(F)$ over $\bb{Z}_p$, then $\n{JL}(\pi)$ is a complex of $p$-adically complete admissible representations of $D^{\times}$, i.e. the cohomology groups 
\[
\n{JL}^i(\pi) =H^i_{\proet}(\bb{P}^{n-1}_{\bb{C}_p}, \n{F}_{\pi})
\]
are $p$-adically complete admissible representations of $D^{\times}$. Moreover, we have
\[
\n{JL}^i(\pi)=\varprojlim_{s} \n{JL}^i(\pi/p^s).
\]
\end{corollary}
\begin{proof}
Let $K_D\subset D^{\times}$ be a uniform pro-$p$ open subgroup. 
The complex $\n{JL}(\pi)^{\vee}=R\Hom(\n{JL}(\pi),\bb{Z}_p)$ is a $p$-adically complete module over the Iwasawa algebra $\bb{Z}_{p,\square}[K_D]$ whose reduction modulo $p$ is a perfect $\bb{F}_{p,\square}[K_D]$-complex by  \Cref{Theo:AdmissibleLTScholze}. This implies that $\n{JL}(\pi)^{\vee} $ is itself a perfect complex of $\bb{Z}_{p,\square}[K_D]$-modules and so $\n{JL}(\pi)$ can be represented by a complex of admissible representations of $D^{\times}$. The rest of the statements are classical and left to the reader, see for example \cite[Proposition 1.2.12]{EmertonInterpolation}. 
\end{proof}

\subsection{Locally analytic Jacquet-Langlands functor}

Next we show that the Jacquet-Langlands functor of \cref{Def:AdmissibleTorsion} is compatible with locally analytic vectors. Let $\Pi$ be an admissible locally analytic representation of $\GL_{n}(F)$ over $\bb{Q}_p$, we let $\n{F}_{\Pi}$ be the pro\'etale sheaf over $\bb{P}^{n-1}_{\breve{F}}$ sending an affinoid perfectoid $S\to \bb{P}^{n-1,\lozenge}_{\breve{F}}$ to 
\[
\n{F}_{\Pi}(S)= (C(|\n{M}_{\infty}\times_{\bb{P}^{n-1,\lozenge}_{\breve{F}}} S|, \bb{Q}_p)\widehat{\otimes}_{\bb{Q}_p} \Pi )^{\GL_n(F)}
\]
where $\GL_{n}(F)$ acts via the diagonal action and 
the completed tensor product is a tensor product of LB representations (equivalently a solid tensor product). 
This is the same as the pro\'etale solid sheaf on $\bb{P}^{n-1}_{\breve{F}}$  obtained by descent from the constant sheaf on $\n{M}_{\infty}$ via \cite[Corollary 4.5]{anschutz2024descent}.

\begin{theorem}\label{TheoJLLocAn}
Let $\pi$ be a $p$-adically complete admissible representation of $\GL_n(F)$ and let $\Pi=(\pi[\frac{1}{p}])^{\GL_n(F)-\rm la}$. There is a natural equivalence
\begin{equation}\label{eqDerivedlaJL}
\left(\n{JL}(\pi)[\frac{1}{p}]\right)^{RD^{\times}-{\rm la}} \cong  R\Gamma_{\proet}(\bb{P}^{n-1}_{\bb{C}_p}, \n{F}_{\Pi}).
\end{equation}
Moreover, for all $i\in \bb{Z}$ we have an isomorphism of locally analytic admissible $D^{\times}$-representations 
\begin{equation}\label{eqNonDerivedJL}
\left(\n{JL}^i(\pi)[\frac{1}{p}] \right)^{D^{\times}-{\rm la}}\cong H^{i}_{\proet}(\bb{P}^{n-1}_{\bb{C}_p}, \n{F}_{\Pi}).
\end{equation}
\end{theorem}
\begin{proof} In the following proof we work with the derived  $\infty$-categories of solid  sheaves of diamonds as in \cite[\S 4]{anschutz2024descent}. 

\textit{Step 0.}  The equivalence in \eqref{eqNonDerivedJL} follows from \eqref{eqDerivedlaJL}. Indeed, the object $\n{JL}(\pi)[\frac{1}{p}]$ is a complex with cohomologies given by admissible Banach representations of $D^{\times}$. By \cite[Proposition  4.48]{RRLocallyAnalytic} (see also \cite[Proposition 2.3.1]{RCLocAnCompleted}) the higher locally analytic vectors of a Banach admissible representation vanish, then by the spectral sequence of \cite[Theorem 1.5]{RCLocAnCompleted} one deduces that 
\[
H^i((\n{JL}(\pi)[\frac{1}{p}])^{RD^{\times}-{\rm la}}) = (\n{JL}^i(\pi)[\frac{1}{p}] )^{D^{\times}-{\rm la}}.
\]

\textit{Step 1.} We first reinterpret the problem using the period sheaves.  By \cite[Proposition II.2.5]{FarguesScholze} we have a short exact sequence of pro\'etale sheaves
\[
0\to \bb{Q}_p\to \bb{B}_{[1,p]} \xrightarrow{\varphi-1} \bb{B}_{[1,1]}\to 0.
\]
Taking solid (equivalently $p$-complete in this case) tensor products with the sheaf $\n{F}_{\pi}$ we get a short exact sequence
\[
0 \to \n{F}_{\pi}[\frac{1}{p}]\to \bb{B}_{[1,p]}\widehat{\otimes}_{\bb{Z}_p} \n{F}_{\pi} \to \bb{B}_{[1,1]}\widehat{\otimes}_{\bb{Z}_p} \n{F}_{\pi}\to 0. 
\]
Taking pro\'etale cohomology we get an exact triangle
\[
\n{JL}(\pi)[\frac{1}{p}]\to R\Gamma_{\proet}( \bb{P}^{n-1}_{\bb{C}_p}, \bb{B}_{[1,p]}\widehat{\otimes}_{\bb{Q}_p} \n{F}_{\pi}) \to R\Gamma_{\proet}( \bb{P}^{n-1}_{\bb{C}_p}, \bb{B}_{[1,1]}\widehat{\otimes}_{\bb{Q}_p} \n{F}_{\pi})  \xrightarrow{+}. 
\]
On the other hand, taking LB-completed tensor products we get a short exact sequence of pro\'etale sheaves
\[
0\to \n{F}_{\Pi}\to \bb{B}_{[1,p]}\widehat{\otimes}_{\bb{Q}_p} \n{F}_{\Pi} \to \bb{B}_{[1,1]}\widehat{\otimes}_{\bb{Q}_p} \n{F}_{\Pi}\to 0. 
\]
Therefore, to prove the theorem it suffices to show that for all $I\subset (0,\infty)$ compact interval, we have a natural equivalence of  representations of $D^{\times}$
\begin{equation}\label{eqEquivalenceLocAnCoho}
R\Gamma_{\proet}(\bb{P}^{n-1}_{\bb{C}_p}, \bb{B}_I\widehat{\otimes}_{\bb{Z}_p} \n{F}_{\pi})^{RD^{\times}-{\rm la}} \cong R\Gamma_{\proet}(\bb{P}^{n-1}_{\bb{C}_p}, \bb{B}_I\widehat{\otimes}_{\bb{Z}_p} \n{F}_{\Pi}).
\end{equation}

\textit{Step 2.} We now reduce the proof of \eqref{eqEquivalenceLocAnCoho} to affinoid subspaces of $\bb{P}^{n-1}_{\bb{C}_p}$. Let $\nu:\bb{P}^{n-1}_{\bb{C}_p,\proet}\to \bb{P}^{n-1}_{\bb{C}_p,\an}$ be the projection of sites and let   $\f{U}=\{U_i\}_{i\in I}$ be a finite rational open cover of $\bb{P}^{n-1}_{\bb{C}_p}$. Then, for any pro\'etale sheaf $\s{F}$ over $\bb{P}^{n-1}_{\bb{C}_p}$ we have equivalences of complexes 
\[
R\check{\Gamma}_{\an}(\f{U}, R\nu_* \s{F}) \cong R\Gamma_{\proet}(\bb{P}^{n-1}_{\bb{C}_p} ,\s{F}).
\]
functorial on $\s{F}$, where the left hand side is the \v{C}ech cohomology given by 
\[
R\check{\Gamma}_{\an}(\f{U}, R\nu_* \s{F}) = \varprojlim_{V\in \mathrm{Int}(\f{U})} R\Gamma_{\an}(V,  R\nu_* \s{F})
\]
with $\mathrm{Int}(\f{U})$ the poset of finite intersections of elements in $\f{U}$.

 Therefore, in order to show \eqref{eqEquivalenceLocAnCoho} it suffices to prove that for $U\subset \bb{P}^{n-1}_{\bb{C}_p}$ a rational open subspace we have a natural equivalence
\begin{equation}\label{eqLocalEquivalenceRationalCover}
R\Gamma_{\proet}(U, \bb{B}_I\widehat{\otimes}_{\bb{Z}_p} \n{F}_{\pi})^{RD^{\times}-{\rm la}} \cong R\Gamma_{\proet}(U, \bb{B}_I\widehat{\otimes}_{\bb{Z}_p} \n{F}_{\Pi}). 
\end{equation}

\textit{Step 3.}  Finally, we prove \eqref{eqLocalEquivalenceRationalCover}. We can assume without loss of generality that $U$ is a rational subspace admitting a section $U\subset \n{M}_{\bb{X},\bb{C}_p}$.  Let $K_{D,U}\subset D^{\times}$ be a compact open subgroup stabilizing $U$ and let $U_{\infty}= \n{M}_{\infty}\times_{\n{M}_{\bb{X}}} U$.    Then

\[
\begin{aligned}
R\Gamma_{\proet}(U, \bb{B}_{I}\widehat{\otimes}_{\bb{Z}_p} \n{F}_{\pi})^{RD^{\times}-{\rm la}} & \cong R\Gamma(K_{D,U}, R\Gamma(\GL_n(\n{O}), R\Gamma_{\proet}(U_{\infty}, \bb{B}_{I}\widehat{\otimes}_{\bb{Z}_p} \n{F}_{\pi}) \widehat{\otimes}^L_{\bb{Q}_p} C^{\rm la}(K_{D,U},\bb{Q}_p))\\ 
& \cong R\Gamma(K_{D,U}\times \GL_{n}(\n{O}), R\Gamma_{\proet}(U_{\infty}, \bb{B}_I)\widehat{\otimes}_{\bb{Q}_p}^L \pi [1/p] \widehat{\otimes}_{\bb{Q}_p}^LC^{\rm la}(K_{D,U},\bb{Q}_p)) \\ 
 & \cong R\Gamma(\GL_{n}(\n{O}), R\Gamma_{\proet}(U_{\infty}, \bb{B}_I)^{RK_{H}-\rm la}\widehat{\otimes}_{\bb{Q}_p}^L \pi[1/p]) \\ 
 & \cong R\Gamma(\GL_n(\n{O}), R\Gamma_{\proet}(U_{\infty}, \bb{B}_I)^{RG-\rm la}\widehat{\otimes}_{\bb{Q}_p}^L \pi[1/p]) \\ 
 & \cong R\Gamma(\GL_n(\n{O}), R\Gamma_{\proet}(U_{\infty}, \bb{B}_I)^{RG-\rm la}\widehat{\otimes}_{\bb{Q}_p}^L \Pi )\\  
 & \cong R\Gamma(\GL_n(\n{O}), R\Gamma_{\proet}(U_{\infty}, \bb{B}_I)\widehat{\otimes}_{\bb{Q}_p}^L \Pi ) \\ 
 & \cong R\Gamma_{\proet}(U, \n{F}_{\Pi}).
\end{aligned}
\]
In the first equivalence we use descent along the $\GL_n(\n{O})$-torsor $U_{\infty}\to U$ and write explicitly the definition of $K_{D,U}$-locally analytic vectors. The second equivalence is clear as $U_{\infty}$ is qcqs and $\pi$ is $p$-complete. The third equivalence follows from projection formula of locally analytic vectors \cite[Corollary 3.1.15 (3)]{RJRCSolidLocAn2} and the fact that $\pi$ is a trivial $K_{D,U}$-representation. The fourth equivalence is \cref{TheoMainComparisonLocAn}.  The fifth equivalence follows from the projection formula of locally analytic vectors and the fact that $(\pi)[\frac{1}{p}])^{RG-{\rm la}}=\Pi$ as $\pi$ is an admissible representation. The sixth equivalence is the projection formula again. The last equivalence is  descent along the torsor $U_{\infty}\to U$.  This finishes the proof of the theorem.
\end{proof}

As a corollary we can prove that the Jacquet-Langlands functor for Banach admissible locally analytic representations preserves central characters. See 
 \cite{MR4940362} and \cite{MR4621880} for a completely different approach. 

\begin{theorem}\label{CoroCentralChar}
Let $\pi$ be an admissible unitary Banach representation of $G:=\GL_n(F)$ over $\bb{Q}_p$ such that $\pi^{G-\rm la}$ has infinitesimal character $\chi$. For all $i\in \bb{Z}$ the locally analytic $D^{\times}$-representation $\n{JL}^i(\pi)^{D^{\times}-\rm la}$ has infinitesimal character $\chi$ under the natural identification $\n{Z}(\Lie D^{\times})\otimes \bb{C}_p \cong \n{Z}(\Lie G)\otimes \bb{C}_p$. 
\end{theorem}
\begin{proof}
The statement can be proved after base change to $\bb{C}_p$. By \cite[Theorem 3.2]{ScholzeLubinTate}  we have a natural equivalence 
\[
\n{JL}(\pi)\widehat{\otimes}_{\bb{Q}_p} \bb{C}_p=R\Gamma_{\proet}(\bb{P}^{n-1}_{\bb{C}_p}, \n{F}_{\pi}\widehat{\otimes}_{\bb{Q}_p} \widehat{\s{O}}). 
\]
Then \eqref{eqEquivalenceLocAnCoho} yields a $D^{\times}$-equivariant equivalence 
\[
(\n{JL}(\pi)\widehat{\otimes}_{\bb{Q}_p} \bb{C}_p)^{RD^{\times}-\rm la} \cong R\Gamma_{\proet}(\bb{P}^{n-1}_{\bb{C}_p}, \n{F}_{\Pi,\widehat{\s{O}}}),
\]
thus it suffices to show that the RHS term has central character given by $\chi$. By picking a suitable affinoid cover $\{U_i\}_{i}$ of $\bb{P}^{n-1}_{\bb{C}_p}$ as in Steps 2 and 3 of the proof of \cref{TheoJLLocAn}, we are reduced to show that for any small enough open affinoid $U\subset \n{M}_{\bb{X}}$ with stabilizer $K_{D,U}\subset D^{\times}$, the central character of $R\Gamma_{\proet}(U, \n{F}_{\Pi,\widehat{\s{O}}})$ for the action of $K_{D,U}$ is $\chi$. Let $U_{\infty}\subset \n{M}_{\infty}$ be the pullback of $U$ to infinite level, by Step 3 of the proof of \cref{TheoJLLocAn} we have that 
\[
R\Gamma_{\proet}(U, \n{F}_{\Pi,\widehat{\s{O}}})\cong R\Gamma(\GL_n(\n{O}),  R\Gamma_{\proet}(U_{\infty}, \widehat{\s{O}})^{RK_{D,U}-\rm la}\widehat{\otimes}_{\bb{Q}_p} \Pi ),
\]
but by taking $U$ small enough, the vanishing of higher locally analytic vectors of \cref{TheoVanishingLocAnVectors} implies that 
\[
R\Gamma_{\proet}(U, \n{F}_{\Pi,\widehat{\s{O}}})\cong R\Gamma(\GL_n(\n{O}),  \s{O}^{\rm la}_{\n{M}}(U_{\infty})\widehat{\otimes}_{\bb{Q}_p} \Pi).
\]
The corollary follows from the identification of the central horizontal actions $\n{Z}(\f{m}_{\mu})_{\bb{C}_p}\cong \n{Z}(\f{m}_{\mu^{-1}})_{\bb{C}_p}$ on $\s{O}^{\rm la}_{\n{M}}(U_{\infty})$ of \cref{TheoVanishingLocAnVectors}   and the fact that the central actions of $\n{Z}(\Lie D^{\times})\otimes \bb{C}_p \cong \n{Z}(\Lie G)\otimes \bb{C}_p$ factor through the horizontal actions. 
\end{proof}

\bibliographystyle{alpha}
\bibliography{JLlocan}

\end{document}